\documentclass[12pt]{scrartcl}

\usepackage{array}
\usepackage{supertabular}
\usepackage{amsmath,amsfonts,amssymb,amsthm} 
\usepackage{bm}
\usepackage[dvipdfmx]{graphicx}
\usepackage{here}

\newcommand{\secffe}[2]{\overset{\;e}{A}\hspace{-2pt}\begin{smallmatrix}
{#2}\\{#1}\end{smallmatrix}}

\newcommand{\secffm}[2]{\overset{\;\scalebox{0.6}{$m$}}{A}\hspace{-2pt}\begin{smallmatrix}
{#2}\\{#1}\end{smallmatrix}}

\allowdisplaybreaks[4]

\newtheorem{prop}{Proposition}

\makeatletter
\newcommand{\oset}[2]{%
  {\mathop{#2}\limits^{\vbox to -.5\ex@{\kern-\tw@\ex@
   \hbox{\scriptsize #1}\vss}}}}
\makeatother

\usepackage{microtype} 

\usepackage[utf8]{inputenc} 
\usepackage[T1]{fontenc} 

\title {Asymptotic Expansion of Risk for a Regression Model with respect to $\alpha$-Divergence with an Application to the Sample Size Problem}
\author{Yo Sheena\thanks{Faculty of Economics and Law, Shinshu University. Faculty of Data Science, Shiga University}}

\date{April 2017}

\begin{document}
\maketitle

\begin{abstract}
For a regression model,  we consider the risk of the maximum likelihood estimator with respect to $\alpha$-divergence, which includes the special cases of Kullback-Leibler divergence, Hellinger distance and  $\chi^2$ divergence. The asymptotic expansion of the risk with respect to the sample size $n$ is given  up to the order $n^{-2}$.  We observed how the risk convergence speed (to zero) is affected by the error term distributions and the magnitude of the joint moments of the standardized explanatory variables under three concrete error term distributions: a normal distribution, a t-distribution and a skew-normal distribution. We try to use the (approximated) risk of m.l.e. as a measure of the difficulty of estimation for the regression model. Especially comparing the value of the (approximated) risk with that of a binomial distribution, we can give a certain standard for the sample size required to estimate the regression model.
\end{abstract}
\noindent
MSC(2010) \textit{ Subject Classification}: Primary 60F99; Secondary 62F12
\section{Introduction}
\label{section:int}
We consider the following regression model;
\begin{equation}
\label{regression_model}
y=\beta'\tilde{x}+\sigma \epsilon,
\end{equation}
where 
$$
\beta'=(\beta_0,\beta_1,\ldots,\beta_p)
$$
is the $p+1$-dimensional parameter vector, while
$$
\tilde{x}'=(x_0, x'),\quad x_0 \equiv 1,\qquad x'=(x_1,\ldots, x_p), 
$$
and $x'=(x_1,\ldots,x_p)$ is a $p$-dimensional explanatory random vector. $\epsilon$ is the error term. We assume that the distributions of  $\epsilon$ is known, but the distribution of $x$ is unknown. The unknown parameters to be estimated are $\beta \in R^{p+1}$ and $\sigma (>0)$. Without loss of generality, we can assume that $x'=(x_1,\ldots, x_p)$ is standardized, i.e.
\begin{equation}
\label{stand_x}
E[x_i]=0,\quad i=1,\ldots,p,\qquad E[x_i x_j]=
\begin{cases}
1, \text{ if $1\leq i=j \leq p,$} \\ 
0, \text{ if $1\leq i\ne j \leq p.$}
\end{cases}
\end{equation}

Let $f(\epsilon)$ and $h(x)$ respectively be the p.d.f.'s of $\epsilon$ and $x$, then the p.d.f. of $(y,x)$ is given by
\begin{equation}
\label{yx_pdf}
f(y,x \,|\, \beta, \sigma)\triangleq f_x(y|\:\beta,\sigma)\,h(x),
\end{equation}
where 
$$
f_x(y|\:\beta,\sigma)\triangleq \frac{1}{\sigma}f\biggl(\frac{y-\beta'\tilde{x}}{\sigma}\biggr).
$$
We assume that $f(\epsilon) $ is positive and differentiable three times over the real line.

Let's consider the maximum likelihood estimators (say $\hat{\beta}, \hat{\sigma}$) of $\beta, \sigma$. One way to evaluate the performance of m.l.e. is the closeness of  the predictive distribution designated by the p.d.f 
\begin{equation}
\label{yx_pdf_mle}
f(y, x \,|\,\hat{\beta}, \hat{\sigma})= f_x(y|\:\hat\beta,\hat\sigma)\,h(x)=\frac{1}{\hat{\sigma}}f\biggl(\frac{y-\hat{\beta}'\tilde{x}}{\hat{\sigma}}\biggr)\,h(x)
\end{equation}
to the true distribution given by \eqref{yx_pdf}.
 
We adopt divergences as the measure of closeness between two given distributions. A divergence is a premetric. Namely a divergence function $D[ d_1 : d_2]$ evaluated at two distributions $d_1$ and $d_2$  on a same sigma field  $\mathcal{X}$ satisfies 
$$
D[ d_1 : d_2] \geq 0 \text{ for any distributions $d_1$ and $d_2$}
$$
with equality iff $d_1=d_2$, but it is asymmetric,  and the "triangular inequality" does not always hold. 

Among possible divergences, \textit{f-divergence} is natural in dealing with probability distributions. (See Amari and Nagaoka \cite{Amari&Nagaoka}, Vajda \cite{Vajda}.)
First $f$-divergence is  parameter-free. If we change the way of parametrization of a parametric model, $f$-divergence is invariant in the following sense.  Suppose  a distribution $d$ on $\mathcal{X}$ can be designated  by a parameter $\theta$ in a parametric model $P_\theta=\{ (d|\theta) \: | \:\theta \in \Theta \}$, while it is expressed in another parametrization as $(d|\eta)$ in $P_\eta=\{ (d|\eta) \: | \:\eta \in H \}$. If $(d|\theta_i)$ and $(d|\eta_i)$ is the same distribution for $i=1,2$, 
$$
D[ (d|\theta_1) : (d|\theta_2)] = D[ (d|\eta_1) :  (d|\eta_2) ].
$$

Second it is invariant with respect to the transformation between the random variables that retains information.  Let $Y(X)$ be a sufficient statistic for the parametric model of a random object $X$, then $f$-divergence satisfies
$$
D[ (X|\theta_1) : (X|\theta_2) ] = D[ (Y|\theta_1) : (Y|\theta_2) ],
$$
where $(X|\theta_i)$ is the distribution of $X$ given by a parameter $\theta_i$ $(i=1,2)$ .  

In order to proceed a practical investigation of regression models, we need a more specific form of $f$-divergence. In this paper we focus on  an $\alpha$-divergence. It is an  important subclass of $f$-divergence.  Generally a divergence gives a geometrical structure on the manifold of a parametric distribution model, $P_\theta=\{ (d|\theta) \: | \:\theta \in \Theta \}$. (See  Eguchi \cite{Eguchi1}, Amari and Nagoka \cite{Amari&Nagaoka}.) The possible geometrical structures given by $f$-divergence can be realized by $\alpha$-divergences. Furthermore it is a basic divergence from the perspective of information geometry since it gives rise to a "dual" structure between  $\alpha$ and $-\alpha$ for the manifold of the given parametric model (see Eguchi \cite{Eguchi1}, Amari \cite{Amari3}, and Amari and Cichocki \cite{Amari&Cichocki}). Specifically $\alpha$-divergence $(-\infty < \alpha < \infty)$  between the two distributions, each of which is  given respectively by the p.d.f. $f(x; \theta_1)$ and $f(x; \theta_2)$, is defined as 
\begin{equation}
\label{alphadive}
\overset{\alpha}{D}[\theta_1: \theta_2]=
\begin{cases}
\frac{4}{1-\alpha^2}\Bigl\{ 1- \int_{\mathfrak X}f^{(1-\alpha)/2}(x; \theta_1)  f^{(1+\alpha)/2}(x; \theta_2) d\mu \Bigr\}, & \text{ if $\alpha \ne \pm 1,$}\\
\int_{\mathfrak X}f(x; \theta_2) \log \Bigl(f(x; \theta_2)/f(x; \theta_1) \Bigr)d\mu, & \text{ if  $\alpha=1$,}\\
\int_{\mathfrak X}f(x; \theta_1) \log \Bigl(f(x; \theta_1)/f(x; \theta_2) \Bigr) d\mu, & \text{ if $\alpha=-1$.}
\end{cases}
\end{equation}
$\alpha$-divergence is a broad class of divergences. Actually it includes Kullback--Leibler divergence ($\alpha=-1$), the Hellinger distance ($\alpha=0$) and  $\chi^2$ divergence ($\alpha=3$).

We will measure the performance of m.l.e. $\hat{\beta},\ \hat{\sigma}$ by the expected $\alpha$-divergence between two distributions \eqref{yx_pdf} and \eqref{yx_pdf_mle};
\begin{equation}
\label{def_alpha_div}
\overset{\alpha}{E\!D}(\beta,\sigma)\triangleq E\bigl[\overset{\alpha}{D}[(\hat\beta(\bm{y},\bm{x}),\hat\sigma(\bm{y},\bm{x})):(\beta, \sigma)]\bigr],
\end{equation}
where $(\bm{y, x})=\Bigl((y_1,x_1),\ldots,(y_n,x_n)\Bigr)$ are $n$ independent random samples from the true distribution \eqref{yx_pdf}. In other words, we evaluate the performance of m.l.e. using the risk of m.l.e. with respect to an $\alpha$-divergence. However, this risk of m.l.e. can not  be gained explicitly in many (most in a practical sense) cases, hence its asymptotic expansion with $n$ is useful since it gives a good approximation under a large size of samples. Sheena \cite{Sheena} gave the asymptotic expansion of $\overset{\alpha}{E\!D}$ up to the $n^{-2}$ order for a general parametric model. (Henceforth, we will call the truncated $\overset{\alpha}{E\!D}$ up to the $n^{-2}$ order by the name of  "the approximated $\overset{\alpha}{E\!D}$".) In this paper, we focused ourselves on the regression model \eqref{regression_model}, and derived the approximated $\overset{\alpha}{E\!D}$ for it. 

The result for a general regression model \eqref{regression_model} is still too lengthy to be out of use for a practical purpose. So we narrowed our scope further to some specific error distributions. (See Mathematica program in Appendix B which enables us to calculate the approximated  $\overset{\alpha}{E\!D}$, once the p.d.f. (and its derivatives) of an arbitrary error distribution is given.) This  paper is constructed as follows; In Section \ref{general_result}, we explained how the general result of \cite{Sheena} is applied to the regression model.  In Section \ref{homo_x},  we considered three specific error term distributions and observed an explicit form of the expansion of $\overset{\alpha}{E\!D}$: a normal distribution (Section \ref{Normal Distribution Error}), a $t$-distribution (Section \ref{t-Distribution Error}), a skew-normal distribution (Section \ref{Skew_Normal_Distribution_Error}).  In Section \ref{comp_error_dists}, we made a comparison among these three error distributions. Throughout Section \ref{homo_x}, we considered the case where the explanatory variable $x$ has a homogeneous distribution (i.e. invariant w.r.t. the permutations of the $x_i, i=1,\ldots,p$).  We combined the above error term distributions with various types of joint moments of $x$ to gain a concrete form of the approximated $\overset{\alpha}{E\!D}$ as the function of $n, p, \alpha.$ We observed how $n, p, \alpha$ affect $\overset{\alpha}{E\!D}$. In Section \ref{real_data}, we treated two real datasets, which give us examples for non-homogeneous distribution of $x$. 

As one of the possible applications of  $\overset{\alpha}{E\!D}$, we considered the sample size problem, that is, "how  large sample size is required to estimate the parameters of the regression model \eqref{regression_model} ?".
When a parametric distribution model is given, the difficulty of estimation (specification) of the parameter for that model could be measured in various ways. Sheena \cite{Sheena} proposed to measure it by the approximated  $\overset{\alpha}{E\!D}$.  In the paper, the author tried to use the approximated $\overset{\alpha}{E\!D}$ of a binomial distribution model $B(n, p)$ as a benchmark since it gives us an intuitive interpretation. For example, if a parametric distribution model has a similar value of $\overset{\alpha}{E\!D}(\theta)$ (at a given $\theta$) to $B(10, 0.01)$, we can understand that the task of the estimation is hard, since the value $0.01$ is too small to be estimated from as little as 10 samples. On the contrary, $\overset{\alpha}{E\!D}(\theta)$ of the model is close to that of $B(10, 0.5)$, it is a relatively easy task to estimate the parameter.  

In this paper, we formalized this idea and proposed two indicators (\textit{I.D.E.} and \textit{R.S.S.}) that could be used for a sample size problem. In Section \ref{general_result}, we gave the definition of the both indicators.  In Section \ref{homo_x} and \ref{real_data}, we calculated their concrete values under the given error distributions and the moments of $x$, and tried to give a solution to the sample size problem.
\section{Asymptotic risk of m.l.e. w.r.t. $\alpha$-divergence}
\label{general_result}
First we introduce a general result  of Sheena \cite{Sheena} on the asymptotic risk of m.l.e. with respect to $\alpha$-divergence. In order to improve readability, we use Einstein's summation convention, that is, the summation carried out as every pair of upper and lower index moves from 1 to $p$.

Let $\mathcal P$ be a parametric family of probability distributions on a space $\mathfrak{X}$, which is given by a family of positive-valued densities  $f(x; \theta)$ on $\mathfrak{X}$ with respect to a measure $\mu$:
\begin{equation}
\mathcal P=\{f(x; \theta)\: | \: \theta=(\theta^1, \dots, \theta^p) \in \Theta \},
\label{parametric_model}
\end{equation}
where $\Theta$ is an open set in $R^p$.

Consider the maximum likelihood estimator $\hat{\theta}(\bm{X})$ of $\theta$ based on  $n$ samples $\bm{X}=(X_1,\ldots,X_n)$ independently chosen from the distribution $f(x; \theta)$.  Closeness $\hat{\theta}$ and the true parameter $\theta$ is measured by \eqref{alphadive}, namely $\overset{\alpha}{D}[\hat{\theta}(\bm{X}):\theta]$. The risk is defined as the expectation of this random variable;
\begin{equation}
\label{def_risk}
\overset{\alpha}{E\!D}(\theta)\triangleq E_{\theta}\bigl[\overset{\alpha}{D}[\hat\theta(\bm{X}):\theta]\bigr].
\end{equation}

The asymptotic expansion of $\overset{\alpha}{E\!D}$ w.r.t. $n$ is given by 
\begin{align}
&\overset{\alpha}{E\!D} \nonumber\\
&=\frac{p}{2n}+\frac{1}{24n^{2}}\nonumber\\
&\times
\biggl[ 
(\alpha')^2\bigl\{3\overset{\:e}{F}+3T^{ijk}T_{ijk}-6\langle \secffe{i}{j}, (\secffm{j}{i}-\secffe{j}{i})\rangle-3\langle \secffe{i}{i}, (\secffm{j}{j}-\secffe{j}{j})\rangle+3p^2+6p\bigr\}\nonumber\\
&\hspace{8mm}+
\alpha'\bigl\{3\overset{\:e}{F}-5T^{ijk}T_{ijk}-6T_{is}^iT_j^{js}+6\langle \secffe{i}{j}, (\secffm{j}{i}-\secffe{j}{i})\rangle+3\langle \secffe{i}{i}, (\secffm{j}{j}-\secffe{j}{j})\rangle \nonumber\\
&\hspace{20mm}-3p^2-6p\bigr\}\nonumber\\
&\hspace{8mm}
+12\langle \secffe{j}{i}, \secffe{i}{j} \rangle -2\langle \secffe{j}{i}, \secffm{i}{j} \rangle -\langle \secffe{i}{i}, \secffm{j}{j} \rangle+T_{ijk}T^{ijk}+9T_{is}^i T^{js}_j+8\overset{\;e}{R}\hspace{-10pt}\begin{smallmatrix}{\ \  \ \ ij}\\{ij}\end{smallmatrix}-9\overset{\:e}{F}
\biggr]\nonumber\\
&+o(n^{-2}),\label{expan_ED_final}
\end{align}
where $\alpha'=(1-\alpha)/2.$
The main term equals $p/2n$. $p/n$ is the ratio of the number of the parameters to the sample size. (We will call this quantity  ``$p-n$ ratio'' hereafter.) The coefficient of $n^{-2}$, i.e. the terms inside the bracket have a geometrical meaning if we view $\mathcal{P}$ as a Riemannian manifold. We omit the geometrical explanation (see Sheena \cite{Sheena} ), and just describe their formal definitions. 

Define the following notations; for $1\leq i, j, k, l \leq p$, 
\begin{equation}
\label{def_L_{_}}
\begin{split}
&L_{(ij)}\triangleq E_\theta[l_{ij}],\quad L_{ij}\triangleq E_\theta[l_i l_j], \\
&L_{(ij)k}\triangleq E_\theta[l_{ij} l_k], \quad L_{ijk}\triangleq E_\theta[l_i l_j l_k] \\
&L_{(ij)(kl)}\triangleq E_\theta[l_{ij} l_{kl}], \quad L_{(ijk)l}\triangleq E_\theta[l_{ijk}l_l], \quad 
L_{(ij)kl}\triangleq E_\theta[l_{ij}l_k l_l],\quad L_{ijkl}\triangleq E_\theta[l_i l_j l_k l_l], 
\end{split}
\end{equation}
\begin{equation}
\label{def_L}
\begin{split}
&L11\triangleq g^{ij}g^{kl}L_{(il)jk},\quad L12\triangleq g^{ij}g^{kl}L_{(ij)kl},\quad L13\triangleq g^{ij}g^{kl}L_{ijkl},\\
&L14\triangleq g^{ij}g^{kl}L_{(ik)(jl)},\quad L15\triangleq g^{ij}g^{kl}L_{(ij)(kl)},\\
&L21\triangleq g^{ij}g^{kl}g^{su}L_{(ik)s}L_{jlu},\quad L22\triangleq g^{ij}g^{kl}g^{su}L_{(ij)k}L_{lsu}, \\
&L23\triangleq g^{ij}g^{kl}g^{su}L_{iks}L_{jlu},\quad  L24\triangleq g^{ij}g^{kl}g^{su}L_{ijk}L_{lsu},\\
&L25\triangleq g^{ij}g^{kl}g^{su}L_{(ik)s}L_{(jl)u},\quad L26\triangleq g^{ij}g^{kl}g^{su}L_{(ij)k}L_{(su)l},
\end{split}
\end{equation}
where $(g^{ij})$ is the inverse matrix of $(g_{ij})$ given by 
$$
g_{ij}\triangleq L_{ij}(\equiv -L_{(ij)}),\label{int_exp_g}
$$
and 
\begin{align*}
&l_i \triangleq l_i(x; \theta) \triangleq \frac{\partial}{\partial \theta_i} \log f(x ;\theta), \quad l_{ij}\triangleq l_{ij}(x; \theta) \triangleq \frac{\partial^2}{\partial\theta_i \partial\theta_j }\log f(x ;\theta), \quad \cdots, \\
&E_\theta[h(x;\theta)] \triangleq \int_{\mathfrak X} h(x; \theta) f(x; \theta) d\mu.
\end{align*}

Then each term of \eqref{expan_ED_final} is defined as follows.
\begin{align}
\overset{\:e}{F}&=g^{ij}g^{ks}\bigl(2L_{(is)jk}+L_{(ks)ij}+L_{ijks}\bigr)\nonumber\\
&\quad-g^{ks}g^{uj}g^{li}L_{ijk}\bigl(2L_{(su)l}+L_{sul}\bigr)\nonumber\\
&\quad -g^{ti}g^{uj}g^{ks}L_{(it)s}L_{juk}\nonumber\\
&=2L11+L12+L13-2L21-L23-L22,\label{int_exp_F}\\
T_{ijk}T^{ijk}&= L_{ijk} L_{stu} g^{is}g^{jt} g^{ku}=L23, \label{int_exp_T_ikkT^ijk}\\
T_{is}^i T_j^{js}&=L_{ijk} L_{stu} g^{ij} g^{st} g^{uk}=L24, \label{int_exp_T_is^iT_j^js}\\
\overset{\;e}{R}\hspace{-10pt}\begin{smallmatrix}{\ \  \ \ ij}\\{ij}\end{smallmatrix}
&=g^{ij}g^{sk}\bigl(L_{(ki)(js)}-L_{(ij)(ks)}+L_{(ki)js}-L_{(ij)ks}\bigr)\nonumber\\
&\quad +g^{sk}g^{ti}g^{uj}\bigl(-L_{(ki)j}L_{(st)u}+L_{(it)s}L_{(uj)k}+L_{sit}L_{(uj)k}-L_{stu}L_{(ij)k}\bigr)\nonumber\\
&=L14-L15+L11-L12-L25+L26+L22-L21,\label{int_exp_R}\\
\langle \secffe{i}{j}, \secffe{j}{i} \rangle &= g^{jk}g^{li}L_{(ik)(jl)}-g^{jk}g^{li}g^{st}L_{(ik)s}L_{(jl)t}-p\nonumber\\
&=L14-L25-p, \label{int_exp_secffe{i}{j}_secffe{i}{j}}\\
\langle \secffe{i}{i}, \secffe{j}{j} \rangle &= g^{ik}g^{jl}L_{(ik)(jl)}-g^{ik}g^{jl}g^{st}L_{(ik)s}L_{(jl)t}-p^2\nonumber\\
&=L15-L26-p^2,\label{int_exp_secffe{i}{i}_secffe{j}{j}}\\
\langle \secffe{i}{j}, \secffm{j}{i} \rangle &= g^{jk}g^{li}L_{(ik)jl}+g^{jk}g^{li}L_{(ik)(jl)}\nonumber\\
&\hspace{20mm}-g^{jk}g^{li}g^{st}L_{(ik)s}L_{(jl)t}-g^{jk}g^{li}g^{st}L_{(ik)s}L_{jlt}\nonumber\\
&=L11+L14-L25-L21, \label{int_exp_secffm{i}{j}_secffm{i}{j}}\\
\langle \secffe{i}{i}, \secffm{j}{j} \rangle &= g^{ik}g^{jl}L_{(ik)jl}+g^{ik}g^{jl}L_{(ik)(jl)}\nonumber\\
&\hspace{20mm}-g^{ik}g^{jl}g^{st}L_{(ik)s}L_{(jl)t}-g^{ik}g^{jl}g^{st}L_{(ik)s}L_{jlt}\nonumber\\
&=L12+L15-L26-L22.\label{int_exp_secffm{i}{i}_secffm{j}{j}}
\end{align}

Now we apply \eqref{expan_ED_final} to the case where $\mathcal{P}$ is given by
$$
\mathcal{P}=\{f(y,x \,|\, \beta, \sigma)\: | \: \beta \in R^{p+1},\: \sigma>0 \},
$$
where $f(y,x \,|\, \beta, \sigma)$ is given by \eqref{yx_pdf}.

Accordingly we define the following notations; for $i, j, k, l = 0, 1, \ldots, p, \sigma$
\begin{align*}
&L_{(ij)}\triangleq E_{\beta,\sigma}[l_{ij}],\quad L_{ij}\triangleq E_{\beta,\sigma}[l_i l_j], \\
&L_{(ij)k}\triangleq E_{\beta,\sigma}[l_{ij} l_k], \quad L_{ijk}\triangleq E_{\beta,\sigma}[l_i l_j l_k] \\
&L_{(ij)(kl)}\triangleq E_{\beta,\sigma}[l_{ij} l_{kl}], \quad L_{(ijk)l}\triangleq E_{\beta,\sigma}[l_{ijk}l_l], \quad 
L_{(ij)kl}\triangleq E_{\beta,\sigma}[l_{ij}l_k l_l],\quad L_{ijkl}\triangleq E_{\beta,\sigma}[l_i l_j l_k l_l],
\end{align*}
where
$$
l_i \triangleq \partial_i \log f(y,x | \beta, \sigma), \quad l_{ij} \triangleq \partial_i \partial_j \log f(y,x | \beta, \sigma), \quad l_{ijk} \triangleq \partial_i \partial_j \partial_k \log f(y,x | \beta, \sigma)
$$
with $\partial_i (i=0, 1, \ldots, p, \sigma)$ defined by
$$
\partial_i=
\begin{cases}
\frac{\partial}{\partial \beta_i} &\text{ if $0 \leq i \leq p$}, \\
\frac{\partial}{\partial \sigma} &\text{ if $i = \sigma$},
\end{cases}
$$
and
$$
E_{\beta,\sigma}[h(y, x; \beta, \sigma)]=\int_{R^p} \int_R h(y, x; \beta, \sigma) f(y, x | \beta, \sigma) dy dx.
$$
We also define other notations. \\
For $0 \leq i, j \leq p$, 
$$
\delta_{ij}=
\begin{cases}
1 &\text{ if $i=j$,}\\
0 &\text{ if $i\ne j$.}
\end{cases}
$$
For $0 \leq i, j, k, l \leq 4$,
\begin{equation}
\label{def_eta}
\eta[i,j,k,l]\triangleq \int_{-\infty}^\infty \Bigl(\frac{d^3 \log f(y)}{d y^3}\Bigr)^i 
\Bigl(\frac{d^2 \log f(y)}{d y^2}\Bigr)^j
\Bigl(\frac{d \log f(y)}{d y}\Bigr)^k
y^l f(y) dy.
\end{equation}
For $i, j, k, l \in \{0, 1,\ldots, p, \sigma\}$
\begin{equation}
\label{def_m[ ]}
\begin{split}
&m[i, j, k]\triangleq E[\dot{x}_i \dot{x}_j \dot{x}_k]=\int_{R^p} \dot{x}_i \dot{x}_j \dot{x}_k h(x) dx, \\
&m[i, j, k, l]\triangleq E[\dot{x}_i \dot{x}_j \dot{x}_k \dot{x}_l]=\int_{R^p} \dot{x}_i \dot{x}_j \dot{x}_k \dot{x}_l h(x) dx,
\end{split}
\end{equation}
where
$$
\dot{x}_i=
\begin{cases}
x_i &\text{ if $i \in \mathcal{I}\triangleq \{1,2,\ldots,p\}$},\\
1 & \text{ if $i \in \mathcal{S}\triangleq \{0,\sigma\}.$}
\end{cases}
$$

Straightforward calculation leads to the following results (see Appendix A for  the detailed calculation).
\begin{align}
g_{ij}&=
\delta_{ij}\sigma^{-2}\eta[0,0,2,0]=-\delta_{ij}\sigma^{-2}\eta[0,1,0,0], \quad 0 \leq i, j \leq p.\label{g_ij}\\
g_{i\sigma}&=
\begin{cases}
\sigma^{-2}\eta[0,0,2,1]=-\sigma^{-2} \eta[0,1,0,1] & \text{ if $i=0$,} \\
0 & \text{ if $1 \leq i \leq p.$}
\end{cases}
\label{g_is}
\\
g_{\sigma\sigma}&=\sigma^{-2}(1+2\eta[0,0,1,1]+\eta[0,0,2,2]).\nonumber\\
&=-\sigma^{-2}(1+\eta[0,1,0,2]+2\eta[0,0,1,1])\label{g_ss}\\
g^{ij}&=\delta_{ij}\sigma^{2}\eta^{-1}[0,0,2,0],\quad 1 \leq i, j \leq p.\label{g^ij}\\
g^{0i}&=g^{\sigma i}=0, \quad 1\leq i \leq p.\label{g^0i}\\
g^{00}&=\sigma^{2}\Delta^{-1}(1+2\eta[0,0,1,1]+\eta[0,0,2,2]).\label{g^00}\\
g^{0\sigma}&=\sigma^{2}\Delta^{-1}\eta[0,1,0,1].\label{g^0s}\\
g^{\sigma\sigma}&=\sigma^{2}\Delta^{-1}\eta[0,0,2,0].\label{g^ss}\\
&( \Delta=\eta[0,0,2,0](1+2\eta[0,0,1,1]+\eta[0,0,2,2])-\eta^2[0,1,0,1]) \nonumber
\end{align}
For $i, j, k, l =0, 1, \ldots, p, \sigma$, 
\begin{align}
L_{(ij)k}&=\sigma^{-3}m[i,j,k]\eta_{(ij)k}\label{L_(ij)k}\\
L_{ijk}&=\sigma^{-3}m[i,j,k]\eta_{ijk}\label{L_ijk}\\
L_{(ij)(kl)}&=\sigma^{-4} m[i,j,k,l]\eta_{(ij)(kl)}\label{L_(ij)(kl)}\\
L_{(ijk)l}&=\sigma^{-4}m[i,j,k,l]\eta_{(ijk)l}\label{L_(ijk)l}\\
L_{(ij)kl}&=\sigma^{-4}m[i,j,k,l]\eta_{(ij)kl}\label{L_(ij)kl}\\
L_{ijkl}&=\sigma^{-4}m[i,j,k,l]\eta_{ijkl}\label{L_ijkl},
\end{align}
where for $0 \leq i, j, k, l \leq p$,
\begin{align}
\eta_{(ij)k}&=-\eta[0,1,1,0]\label{eta_(ij)k}\\
\eta_{(i\sigma)k}&=-(\eta[0,1,1,1]+\eta[0,0,2,0])\label{eta_(is)k}\\
\eta_{(ij)\sigma}&=-(\eta[0,1,0,0]+\eta[0,1,1,1])\label{eta_(ij)s}\\
\eta_{(i\sigma)\sigma}&=-(\eta[0,1,0,1]+\eta[0,1,1,2]+\eta[0,0,2,1])\label{eta_(is)s}\\
\eta_{(\sigma\sigma)i}&=-(\eta[0,1,1,2]+2\eta[0,0,2,1])\label{eta_(ss)i}\\
\eta_{(\sigma\sigma)\sigma}&=-(1+3\eta[0,0,1,1]+\eta[0,1,0,2]+2\eta[0,0,2,2]+\eta[0,1,1,3])\label{eta_(ss)s}\\
\eta_{ijk}&=-\eta[0,0,3,0]\label{eta_ijk}\\
\eta_{ij\sigma}&=-(\eta[0,0,2,0]+\eta[0,0,3,1])\label{eta_ijs}\\
\eta_{i\sigma\sigma}&=- (2\eta[0,0,2,1]+\eta[0,0,3,2])\label{eta_iss}\\
\eta_{\sigma\sigma\sigma}&=-(1+3\eta[0,0,1,1]+3\eta[0,0,2,2]+\eta[0,0,3,3])\label{eta_sss}\\
\eta_{(ij)(kl)}&=\eta[0,2,0,0]\label{eta_(ij)(kl)}\\
\eta_{(i\sigma)(kl)}&=\eta[0,2,0,1]+\eta[0,1,1,0]\label{eta_(is)(kl)}\\
\eta_{(i\sigma)(j\sigma)}&=\eta[0,2,0,2]+2\eta[0,1,1,1]+\eta[0,0,2,0]\label{eta_(is)(js)}\\
\eta_{(ij)(\sigma\sigma)}&=\eta[0,1,0,0]+\eta[0,2,0,2]+2\eta[0,1,1,1]\label{eta_(ij)(ss)}\\
\eta_{(i\sigma)(\sigma\sigma)}&=\eta[0,1,0,1]+\eta[0,2,0,3]+3\eta[0,1,1,2]+2\eta[0,0,2,1]\label{eta_(is)(ss)}\\
\eta_{(\sigma\sigma)(\sigma\sigma)}
&=1+\eta[0,2,0,4]+4\eta[0,0,2,2]+2\eta[0,1,0,2]\nonumber\\
&\quad+4\eta[0,0,1,1]+4\eta[0,1,1,3]\label{eta_(ss)(ss)}\\
\eta_{(ijk)l}&=\eta[1,0,1,0]\label{eta_(ijk)l}\\
\eta_{(ijk)\sigma}&=\eta[1,0,0,0]+\eta[1,0,1,1]\label{eta_(ijk)s}\\
\eta_{(ij\sigma)k}&=2\eta[0,1,1,0]+\eta[1,0,1,1]\label{eta_(ijs)k}\\
\eta_{(i\sigma\sigma)j}&=4\eta[0,1,1,1]+2\eta[0,0,2,0]+\eta[1,0,1,2]\label{eta_(iss)j}\\
\eta_{(ij\sigma)\sigma}&=2\eta[0,1,0,0]+\eta[1,0,0,1]+2\eta[0,1,1,1]+\eta[1,0,1,2]\label{eta_(ijs)s}\\
\eta_{(i\sigma\sigma)\sigma}&=4\eta[0,1,0,1]+\eta[1,0,0,2]+4\eta[0,1,1,2]+2\eta[0,0,2,1]+\eta[1,0,1,3]\label{eta_(iss)s}\\
\eta_{(\sigma\sigma\sigma)i}&=6\eta[0,1,1,2]+6\eta[0,0,2,1]+\eta[1,0,1,3]\label{eta_(sss)i}\\
\eta_{(\sigma\sigma\sigma)\sigma}&=2+6\eta[0,1,0,2]+6\eta[0,0,1,1]+\eta[1,0,0,3]\nonumber\\
&\quad+2\eta[0,0,1,1]+6\eta[0,1,1,3]+6\eta[0,0,2,2]+\eta[1,0,1,4]\label{eta_(sss)s}\\
\eta_{(ij)kl}&=\eta[0,1,2,0]\label{eta_(ij)kl}\\
\eta_{(ij)k\sigma}&=\eta[0,1,1,0]+\eta[0,1,2,1]\label{eta_(ij)ks}\\
\eta_{(i\sigma)jk}&=\eta[0,1,2,1]+\eta[0,0,3,0]\label{eta_(is)jk}\\
\eta_{(ij)\sigma\sigma}&=\eta[0,1,0,0]+2\eta[0,1,1,1]+\eta[0,1,2,2]\label{eta_(ij)ss}\\
\eta_{(i\sigma)j\sigma}&=\eta[0,1,1,1]+\eta[0,0,2,0]+\eta[0,1,2,2]+\eta[0,0,3,1]\label{eta_(is)js}\\
\eta_{(\sigma\sigma)ij}&=\eta[0,0,2,0]+2\eta[0,0,3,1]+\eta[0,1,2,2]\label{eta_(ss)ij}\\
\eta_{(i\sigma)\sigma\sigma}&=\eta[0,1,0,1]+2\eta[0,1,1,2]+2\eta[0,0,2,1]
+\eta[0,1,2,3]+\eta[0,0,3,2]\label{eta_(is)ss}\\
\eta_{(\sigma\sigma)i\sigma}&=2\eta[0,0,2,1]+\eta[0,1,1,2]+\eta[0,0,2,1]+2\eta[0,0,3,2]+\eta[0,1,2,3]\label{eta_(ss)is}\\
\eta_{(\sigma\sigma)\sigma\sigma}&=1+4\eta[0,0,1,1]+\eta[0,1,0,2]+5\eta[0,0,2,2]\nonumber\\
&\quad+2\eta[0,1,1,3]+2\eta[0,0,3,3]+\eta[0,1,2,4]\label{eta_(ss)ss}\\
\eta_{ijkl}&=\eta[0,0,4,0]\label{eta_ijkl}\\
\eta_{ijk\sigma}&=\eta[0,0,3,0]+\eta[0,0,4,1]\label{eta_ijks}\\
\eta_{ij\sigma\sigma}&=\eta[0,0,2,0]+2\eta[0,0,3,1]+\eta[0,0,4,2]\label{eta_ijss}\\
\eta_{i\sigma\sigma\sigma}&=3\eta[0,0,2,1]+3\eta[0,0,3,2]+\eta[0,0,4,3]\label{eta_isss}\\
\eta_{\sigma\sigma\sigma\sigma}&=1+4\eta[0,0,1,1]+6\eta[0,0,2,2]+4\eta[0,0,3,3]+\eta[0,0,4,4].\label{eta_ssss}
\end{align}

If we insert these results \eqref{g^ij},...,\eqref{L_ijkl}
 into \eqref{def_L}, we can calculate the values of  \eqref{int_exp_F} to \eqref{int_exp_secffm{i}{i}_secffm{j}{j}}. Note that the summation (by Einstein's convention) in \eqref{def_L} to \eqref{int_exp_secffm{i}{i}_secffm{j}{j}} is carried over the range $0, 1, \ldots, p, \sigma$ for each index. The calculation process is so lengthy that we used Mathematica \cite{Mathematica}.  The general result expressed with abstract notations $\eta[i, j,k,l]$ (see \eqref{def_eta}) and $m[i,j,k], m[i, j,k, l]$ (see \eqref{def_m[ ]}) could be given, but it is too complicated to be out of use. Instead we put the Mathematica program in Appendix B so that we can easily calculate the approximated $\overset{\alpha}{E\!D}$ once the error term distribution and the moments of the explanatory variables are given, which  respectively determine  $\eta[i, j,k,l]$ and $m[i,j,k], m[i, j,k, l]$. 

Generally $\overset{\alpha}{E\!D}$ for the parametric model \eqref{parametric_model} depends on $\theta$. However $\overset{\alpha}{E\!D}$ for the regression model \eqref{regression_model} is independent of $\beta, \sigma$. This is obvious from the fact that \eqref{g^ij},...,\eqref{L_ijkl} include only $\sigma$, but it vanishes at \eqref{def_L}. We report that if the support of $f(\epsilon)$ is not the whole real line (e.g. $f(\epsilon)=0$ for negative values of $\epsilon$), $\eta[i, j, k, l]$, hence $\overset{\alpha}{E\!D}$ could be dependent on $(\beta, \sigma)$.

In the next section, we give the explicit result when an error distribution and the moments of $x$ are specified. We consider three specific cases where the error term distribution is respectively a normal distribution, a t-distribution and a skew-normal distribution. The different sets of the moments of $x$ are combined with these error distributions to give illustrating examples.

Now we mention one of the possible applications of the approximated $\overset{\alpha}{E\!D}$. For a parametric distribution model \eqref{parametric_model}, we naturally raise the following questions;
\begin{enumerate}
\item At which point $\theta$, is the parameter most difficult to be estimated ?
\item Compared with another model, this model is easier or harder to be estimated ?
\end{enumerate}
We propose to use the approximated  $\overset{\alpha}{E\!D}$ to give an answer to these questions. Maximum likelihood is the most common estimation method and intrinsic to the model, hence it is natural to measure ``the difficulty of estimating the model `` by its performance such as the risk w.r.t a certain loss function. As we mentioned in Introduction, the risk w.r.t. $\alpha$-divergence has favorable properties to answer to the above questions. In this paper we will use the approximated $\overset{\alpha}{E\!D}$ as a measure of the estimation difficulty.

In the case of  the regression model \eqref{regression_model}, the answer to the first question is obvious. Since $\overset{\alpha}{E\!D}$ is constant (independent of $\beta, \sigma^2$), the difficulty of estimation is same all over the parameter space. Concerning the second question, we take the binomial distribution model $B(n,m)$ ($n$: the sample size, $m$: the probability of an event) as the benchmark for comparison. 

\begin{figure}
\centering
\includegraphics[width=11cm,clip]{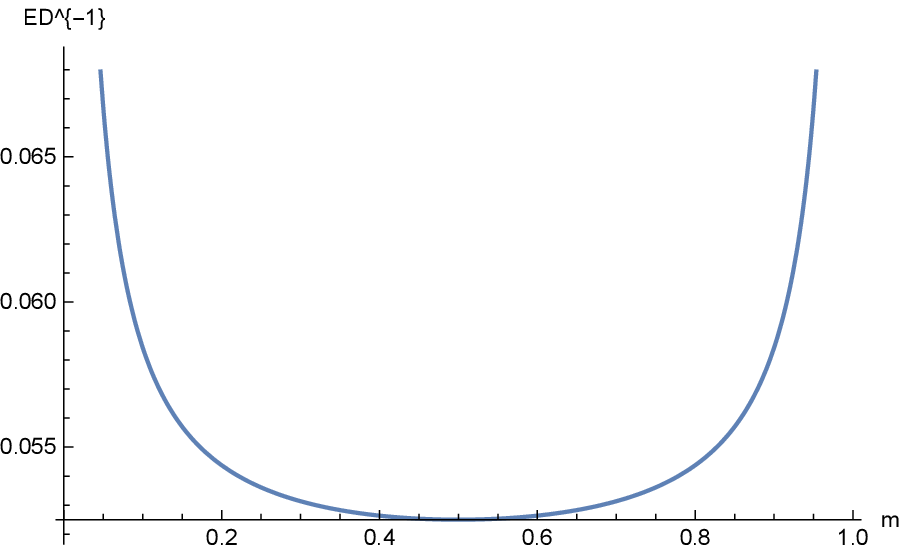}
\caption{$\overset{-1}{E\!D}$ of  $B(10,m)$}
\label{fig:binomial}
\end{figure}

The asymptotic expansion of $\overset{\alpha}{E\!D}$ for the  binomial distribution $B(n,m)$ is given by
\begin{align}
&\overset{\alpha}{E\!D}\nonumber\\
&=\frac{1}{2n}+\frac{1}{24n^{2}}\biggl[ 
(\alpha')^2(3M-9)+\alpha'(-11M+29)+10M-22
\biggr]+o(n^{-2}), \label{expan_binom}
\end{align}
where $\alpha'=(1-\alpha)/2$ and $M\triangleq 1/m+1/(1-m)$. (See the subsection 3.2 of Sheena \cite{Sheena}.)
For Kullback-Leibler divergence, put $\alpha=-1$, then we have
\begin{equation}
\overset{\scalebox{0.6}{\!\!$-1$}}{E\!D} =\frac{1}{2n}+\frac{1}{12n^{2}}(M-1)+o(n^{-2}).\label{expan_KL_binom}\\
\end{equation}
The graph of the approximated $\overset{\scalebox{0.6}{\!\!$-1$}}{E\!D}$ for $B(10, m)$ is given in Figure \ref{fig:binomial}.
We notice that the approximated $\overset{\scalebox{0.6}{\!\!$-1$}}{E\!D}$ is stable around the area $0.1 \leq m \leq 0.9$, however it rapidly increases outside this area. 

Let $\overset{\alpha}{E\!D}_{B(n,m)}$ denote the approximated $\overset{\alpha}{E\!D}$ for $B(n,m)$ and $\overset{\alpha}{E\!D}_R(n)$ denote that for a specific regression model where all the elements of the regression model ($p$, the error term distribution, the moments of $x$) are specified, hence $\overset{\alpha}{E\!D}_R$ is considered as the function of the sample size $n$.
Here we propose an indicator of the difficulty of estimation. 
\\
--\textit{Indicator of the Difficulty of Estimation (I.D.E.)}--
\begin{quotation}
Use a k times binomial experiment $B(k, m)$ as a benchmark. Solve the equation on $m$
\begin{equation}
\label{equation_difficulty_estimation}
\overset{\alpha}{E\!D}_{B(k,m)}=\overset{\alpha}{E\!D}_R((p+2)k)
\end{equation}
\end{quotation}
We easily notice the equation \eqref{equation_difficulty_estimation} is independent of $k$. Taking the sample size for the regression model as $(p+2)k$, we get the same $p-n$ ratio $1/k$ between the two models. Hence it  makes sense to compare the $n^{-2}$ order terms. The solution $m$ tells us intuitively how difficult the parameter estimation is for the regression model. For example if $m=0.001$, then we easily understand the estimation is difficult since it is difficult to estimate $m$ as small as 0.001 based on just 10 samples. On the contrary, if we have $m=0.8$, then the estimation from 10 samples seems not so hard unless we require high precision.

The above equation \eqref{equation_difficulty_estimation} might have no real roots, that is, the left-hand side of the equation is larger than the right-hand side for any $m$. In this case, we can conclude that the regression model could be estimated more easily than the binomial model with the same $p-n$ ratio.

In a reverse way, we can use  the approximated $\overset{\alpha}{E\!D}$ of the regression model for giving an answer to the sample size problem, that is, how large sample size is required to estimate the parameters of the regression model \eqref{regression_model}. 
\\
--\textit{Required Sample Size (R.S.S.)}--
\begin{quotation}
Use a 10 times normal coin toss $B(10, 0.5)$ as a benchmark. Solve the equation
\begin{equation}
\label{equation_sample_size}
\overset{\alpha}{E\!D}_{B(10,0.5)}=\overset{\alpha}{E\!D}_R(n).
\end{equation}
\end{quotation}
The solution $n$ indicates the sample size large enough to guarantee as easy estimation as 10 times normal coin toss. 

The equation \eqref{equation_sample_size} could have no real roots. This means that the left-hand side of the equation is larger than the right-hand side for any $n$. Since the equation is based on the ``approximated'' $\overset{\alpha}{E\!D}$, we must notice that this does not necessarily mean just a small sample (e.g. $p+2$ samples) is enough for the estimation of the regression model. For the approximation to work well, the appropriate sample size is needed. If  we want a concrete solution on the sample size problem, it could be gained by choosing an appropriately large $k$ of $B(k,0.5)$ instead of 10 on the left-hand side of \eqref{equation_sample_size}.
\section{Homogeneous  Explanatory Variables}
\label{homo_x}
\begin{figure}
\centering
\includegraphics[width=11cm,clip]{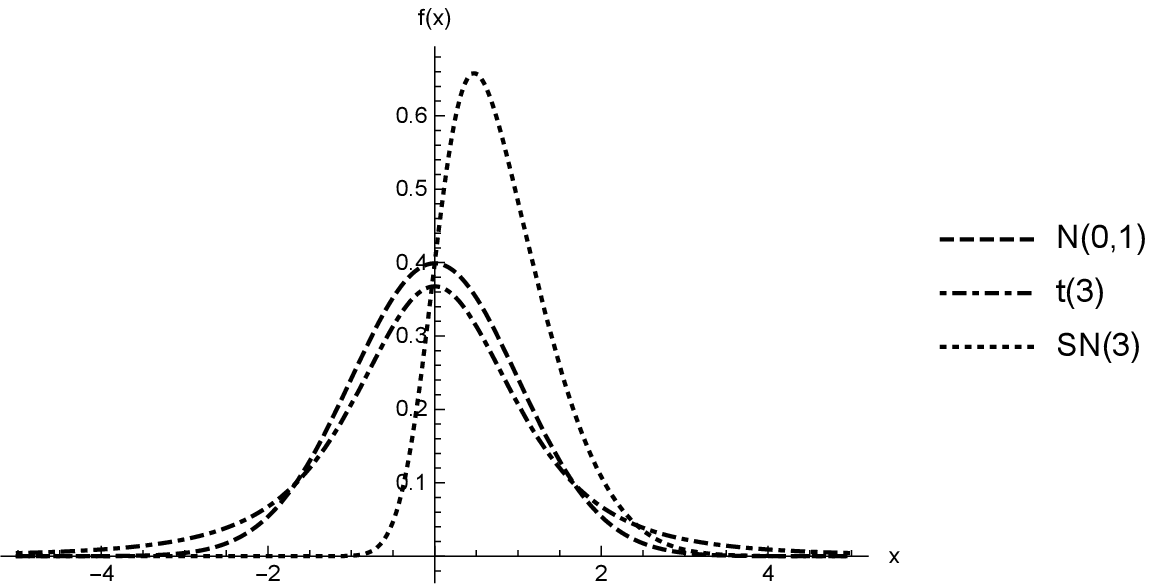}
\caption{p.d.f.'s of three error distributions.}
\label{fig:p.d.f_error}
\end{figure}

In this section, we consider three concrete forms of error distribution: a normal distribution, a t-distribution and a skew-normal distribution. 
A normal distribution is a theoretically basic error distribution. We are interested in how the fat tail property of  a $t$-distribution or the skewness of a skew-normal distribution affects the (approximated) $\overset{\alpha}{E\!D}$. For contrasting these properties, we choose 3 for the d.f. of the $t$ distribution and 3 for the shape parameter of the skew-normal distribution. Figure \ref{fig:p.d.f_error} is the graph of the p.d.f.'s of the three error distributions; the standard normal distribution ($N(0,1)$), the t-distribution with the d.f. of  3 ($t(3)$), the skew-normal distribution with the shape parameter of 3 ($SN(3)$).

$\overset{\alpha}{E\!D}$ also depends on the moments of $x$. As we can see from the definition \eqref{L_(ij)k}--\eqref{L_ijkl}, the maximum order of the joint moments of $x$ is four that appear in the expansion of $\overset{\alpha}{E\!D}$ up to the $n^{-2}$ order.
In this section, we consider the homogeneous case where the distribution of $x=(x_1, \ldots, x_p)$ is invariant w.r.t. any permutation of the elements. This is not practical but this case helps us observe the effect of the dimension $p$, so called "the curse of dimension". 

Here we define the notations of the homogeneous moments of $x$ as follows. For all distinguished $i, j, k, l\ (1 \leq i, j, k, l \leq p)$,
\begin{align}
&m_{4}\triangleq E[x_i^4],\ 
m_{31}\triangleq E[x_i^3 x_j],\ 
m_{22}\triangleq E[x_i^2 x_j^2],\nonumber\\ 
&m_{211}\triangleq E[x_i^2 x_j  x_k],\ 
m_{1111}\triangleq E[x_i  x_j  x_k x_l],\label{4th_moments}\\
&m_{3}\triangleq E[x_i^3],\ 
m_{21}\triangleq E[x_i^2 x_j],\ 
m_{111}\triangleq E[x_i x_j x_k],\label{3rd_moments}\\
&m_{2}\triangleq E[x_i^2]=1,\ 
m_{11}\triangleq E[x_i x_j]=0,\\
&m_{1}\triangleq E[x_i]=0,\\
&m_{0}\triangleq E[x_0]=1.
\label{def_m[i,j,k,l]}
\end{align} 
Under these homogeneous moments, we can state the approximated $\overset{\alpha}{E\!D}$ explicitly for each error distribution as a function of $n, p, \alpha$ and these moments . The result is given in the following subsections. 

We used the following four distributions of $x$ as specific examples of the moments of $x$ when we want to analyze the approximated $\overset{\alpha}{E\!D}$ in a more concrete form:
\begin{enumerate}
\item The standard $p$-dimensional normal distributions, $N_p(0, I_p)$
\begin{align*}
&m_{4}=3,\ m_{31}=0,\ m_{22}=1,\ m_{211=0},\ m_{1111}=0,\\
&m_{3}=0,\ m_{21}=0,\ m_{111}=0.
\end{align*}
\item The standard $p$-dimensional $t$-distribution, $t_p(0, I_p, \nu)$, that is, the $p$ dimensional multivariate $t$-distribution with zero mean vector, the unit matrix as the scale matrix and the degree of freedom $\nu$. Its p.d.f. is given by 
$$
h(x) \varpropto \Bigl(1+\nu^{-1}\sum_{i=1}^p x_i^2 \Bigr)^{-(\nu+p)/2}
$$
Note that $E[x_i]=0,\ i=1,\ldots,p$ and
$$
Cov(x_i, x_j)=E[x_i x_j]=
\begin{cases}
\nu/(\nu-2) &\text{ if  $i=j$,} \\
0 &\text{ if $i \ne j$.}
\end{cases}
$$
for $\nu >2$. Therefore after the normalization \eqref{stand_x}, we have 
$$
E[x_i^2 x_j^2]=
\begin{cases}
3(\nu-2)/(\nu-4) & \text{ if $i=j$,}\\
(\nu-2)/(\nu-4) & \text{ if $i \ne j$,}
\end{cases}
$$
under the condition $\nu>4.$ Notice that the effect of the fourth moment is enhanced by $(\nu-2)/(\nu-4)$ compared to the case $x \sim N_p(0, I_p)$.  We want to check the effect of  the fat tail property  of a $t$-distribution. Here we put $\nu$ as 4.2, then we have
\begin{align*}
&m_{4}=33,\ m_{31}=0,\ m_{22}=11,\ m_{211}=0,\ m_{1111}=0,\\
&m_{3}=0,\ m_{21}=0,\ m_{111}=0.
\end{align*}
\item A completely controlled distribution, where each $x_i, i=1,\ldots,p$ is independently and identically distributed as $P(x_i=1)=P(x_i=-1)=1/2$.
\begin{align*}
&m_{4}=1,\ m_{31}=0,\ m_{22}=1,\ m_{211}=0, \ m_{1111}=0,\\
&m_{3}=0,\ m_{21}=0,\ m_{111}=0.
\end{align*}
\item Pareto distributions, where  each $x_i, i=1,\ldots,p$ is independently and identically distributed as $P(b)$, Pareto distributions with Pareto index $b$.  Its p.d.f. is given by
$$
h(x)=
\begin{cases}
\prod_{i=1}^p bx_i^{-(b+1)} & \text{if $x_i>1$ for $i=1,\ldots,p,$} \\
0 &\text{otherwise.}
\end{cases}
$$
After the normalization \eqref{stand_x}, we have
\begin{align*}
m_{3}&=\text{Skewness of P(b)}=\frac{2 (b+1)}{b-3} \sqrt{\frac{b-2}{b}},\quad b>3,\\
m_{4}&=\text{Kurtosis  of P(b)}=\frac{6 \left(b^3+b^2-6 b-2\right)}{b (b-3) (b-4)}+3
\quad b>4.
\end{align*}
We are interested in the effect of the strong skewness and heavy tail of Pareto distribution. Here we put $b$ as 4.2.
Consequently 
\begin{align*}
&m_{4}=\frac{8129}{21},\ m_{31}=0,\ m_{22}=1,\ m_{211}=0,\ m_{1111}=0,\\
&m_{3}=\frac{26}{63}\sqrt{231},\ m_{21}=0,\ m_{111}=0.
\end{align*}
\end{enumerate}
\subsection{Normal Error Term Distribution}
\label{Normal Distribution Error}
Suppose that the distribution of $\epsilon$ is the  standard normal distribution, that is, 
$$
f(y)=\frac{1}{\sqrt{2\pi}}\exp\Bigl(-\frac{y^2}{2}\Bigr).
$$
Since 
\begin{align*}
&\log f(y)=-\frac{1}{2}\log(2\pi)-\frac{y^2}{2},\\
&\frac{d \log f(y)}{d y}=-y,\qquad\frac{d^2 \log f(y)}{d y^2}=-1,\qquad
\frac{d^3 \log f(y)}{d y^3}=0,
\end{align*}
we have for $0 \leq i, j, k, l \leq 4$,
\begin{align}
\eta[i,j,k,l]&=\frac{1}{\sqrt{2\pi}}\int_{-\infty}^\infty  0^i (-1)^{(j+k)} y^{k+l} \exp\Bigl(-\frac{y^2}{2}\Bigr)dy\nonumber\\
&=
\begin{cases}
0 &\text{ if $i \geq 1$,}\\
0 &\text{ if $i=0$ and $k+l$ is odd,}\\
(-1)^{(j+k)}(k+l-1)!! & \text{ if $i=0$ and $k+l$ is even}, \label{eta_normal}
\end{cases}
\end{align}
where $-1!!=1,\ 1!!=1,\ 3!!=3\cdot 1,\ 5!!=5\cdot 3 \cdot 1,\cdots.$ 

Skipping long calculation (See Appendix B for the calculation procedure using Mathematica \cite{Mathematica}. ), we give the final result in three different expressions  for $n^{-2}$ order term (each expression focuses on respectively $\alpha$, the moments of $x$, $p$).
\begin{align}
&\overset{\alpha}{E\!D}\nonumber\\
&=\frac{p+2}{2 n}\nonumber\\
&\quad+\frac{1}{96 n^2}\Bigl(\alpha^2 \bigl(84 + (48 - 9 m_{22} + 9 m_{4}) p + 9 m_{22} p^2\bigr)\nonumber\\
&\qquad- 8\alpha \bigl( -25 - 3 (6 + m_{22} - m_{4}) p + 3 (-1 + m_{22}) p^2\bigr)\nonumber\\
&\qquad+300 + 240 p + 81 m_{22} p - 81 m_{4} p + 48 p^2 - 81 m_{22} p^2\Bigr)\nonumber\\
&\qquad+o(n^{-2})\label{ED_normal_1}\\
&=\frac{p+2}{2 n}\nonumber\\
&\quad+\frac{1}{96 n^2}\Bigl(3p(-27- 8 \alpha + 3 \alpha^2) m_{4}+3p(p-1)(-27- 8 \alpha + 3 \alpha^2)m_{22}\nonumber\\
&\hspace{30mm}4\alpha^2 ( 12 p+21 ) + 4\alpha ( 6 p^2 + 36 p +50) \nonumber\\
&\hspace{30mm}+4(12 p^2+ 60 p+75)\Bigr)\nonumber\\
&\qquad+o(n^{-2})\label{ED_normal_2}\\
&=\frac{p+2}{2 n}\nonumber\\
&\quad+\frac{1}{96 n^2}\Bigl( (48 + 24 \alpha - 81 m_{22} - 24 \alpha m_{22} + 9 \alpha^2 m_{22})  p^2\nonumber\\
&\hspace{30mm}+ (240 + 144 \alpha + 48 \alpha^2 + 81 m_{22} + 24 \alpha m_{22} - 9 \alpha^2 m_{22}\nonumber\\
&\hspace{30mm}\quad -81 m_{4} - 24 \alpha m_{4} + 9 \alpha^2 m_{4}) p\nonumber\\
&\hspace{30mm}+300 + 200 \alpha + 84 \alpha^2\Bigr)\nonumber\\
&\qquad+o(n^{-2})\label{ED_normal_3}
\end{align}
The $n^{-2}$ order term has the following properties; 
\begin{enumerate}
\item The maximum dimension of  $p$ is two, hence $\overset{\alpha}{E\!D}$ is asymptotically determined by the $p-n$ ratio.
\item  Other moments than $m_4$ and $m_{22}$ do not appear.
\item The coefficients of $m_4$ and $m_{22}$ are non-positive when
 $3\alpha^2-8\alpha-27 \leq 0$, that is, 
\begin{equation}
\label{ED_normal_alpha_int}
-1.95\cdots \leq \alpha \leq 4.62\cdots.
\end{equation}
For $\alpha$ within this interval, the larger $m_{22}$ or $m_4$ gets, the less $\overset{\alpha}{E\!D}$ becomes. The divergences often used in statistical literature are all included in this interval: K-L divergence ($\alpha=-1$),  K-L dual divergence ($\alpha=1$),  Hellinger Divergence ($\alpha=0$), $\chi^2$ divergence ($\alpha=3$).
\end{enumerate}
Most commonly used $\alpha$-divergence is Kullback-Leibler divergence given by $\alpha=-1$;
\begin{align}
&\overset{\scalebox{0.6}{\!\!$-1$}}{E\!D}\nonumber\\
&=\frac{p+2}{2 n}+\frac{1}{12 n^2}\Bigl(23 + 6 (3 + m_{22} - m_4) p + (3 - 6 m_{22}) p^2\Bigr)\nonumber\\
&\quad+o(n^{-2})\label{ED_-1_n}
\end{align} 
Its "dual" divergence given by $\alpha=1$ satisfies the relationship
$$
\overset{\scalebox{0.6}{\!\!$1$}}{D}[\theta_1: \theta_2]=\overset{\scalebox{0.6}{\!\!$-1$}}{D}[\theta_2: \theta_1] . 
$$
The risk of m.l.e. with respect to this divergence is given by;
\begin{align}
&\overset{\scalebox{0.6}{\!\!$1$}}{E\!D}\nonumber\\
&=\frac{p+2}{2 n}+\frac{1}{12 n^2}\Bigl(73 + 6 (9 + 2 m_{22} - 2 m_4) p + (9 - 12 m_{22}) p^2\Bigr)\nonumber\\
&\quad+o(n^{-2})\label{ED_1_n}
\end{align} 
When $\alpha=0$, the $\alpha$-divergence becomes a distance, which is called Hellinger distance;
\begin{align}
&\overset{\scalebox{0.6}{\!\!$0$}}{E\!D}\nonumber\\
&=\frac{p+2}{2n}+\frac{1}{32 n^2}\Bigl(100 + (80 + 27 m_{22} - 27 m_4) p + (16 - 27 m_{22}) p^2\Bigr)\nonumber\\
&\quad+o(n^{-2})
\label{ED_0_n}
\end{align} 
When $\alpha=3$, it becomes $\chi^2$ divergence;
\begin{align}
&\overset{\scalebox{0.6}{\!\!$3$}}{E\!D}\nonumber\\
&=\frac{p+2}{2n}+\frac{1}{4 n^2}\Bigl(69 + (46 + 3 m_{22} - 3 m_4) p + (5 - 3 m_{22}) p^2\Bigr)\nonumber\\
&\quad+o(n^{-2})
\label{ED_2_n}
\end{align} 

For each distribution of $x$ introduced in the beginning of Section \ref{homo_x}, we have the following results.
\begin{enumerate}
\item  $x\sim N_p(0, I_p)$
\begin{align}
&\overset{\alpha}{E\!D} \nonumber\\
&=\frac{p+2}{2n}\nonumber\\
&\quad+\frac{1}{96 n^2}\Bigl(\alpha^2(84 + 66 p + 9 p^2) + 8 \alpha (25 + 12 p)+300 + 78 p - 33p^2\Bigr)\nonumber\\
&\quad+o(n^{-2})
\label{nncase_ED}
\end{align}
\item $x \sim t_p(0, I_p, 4.2)$
\begin{align}
&\overset{\alpha}{E\!D} \nonumber\\
&=\frac{p+2}{2n}\nonumber\\
&\quad+\frac{1}{96 n^2}\Bigl(3 \alpha^2 (28 + 82 p + 33 p^2)  + 8\alpha  (-25 + 48 p + 30 p^2)+300 - 1542 p - 843 p^2\Bigr)\nonumber\\
&\quad+o(n^{-2})
\label{ntcase_ED}
\end{align}
\item $x$ is controlled.
\begin{align}
&\overset{\alpha}{E\!D} \nonumber\\
&=\frac{p+2}{2n}\nonumber\\
&\quad+\frac{1}{96 n^2}
\Bigl(\alpha^2(84 + 48 p + 9 p^2)+8\alpha (25 + 18 p) +300 + 240 p - 33 p^2\Bigr)\nonumber\\
&\quad+o(n^{-2})
\label{nccase_ED}
\end{align}
\item $x_i$ is i.i.d. as $P(4.2)$ $(i=1,\ldots,p)$
\begin{align}
&\overset{\alpha}{E\!D} \nonumber\\
&=\frac{p+2}{2n}\nonumber\\
&\quad+\frac{1}{672 n^2}\Bigl(3 \alpha^2(196 + 8220 p + 21 p^2) -8 \alpha(-175 + 7982 p)\nonumber\\
&\qquad\qquad\qquad+3(700 - 72412 p - 77 p^2)\Bigr)\nonumber\\
&\quad+o(n^{-2})
\label{npcase_ED}
\end{align}
\end{enumerate}

\begin{figure}
\centering
\includegraphics[width=11cm,clip]{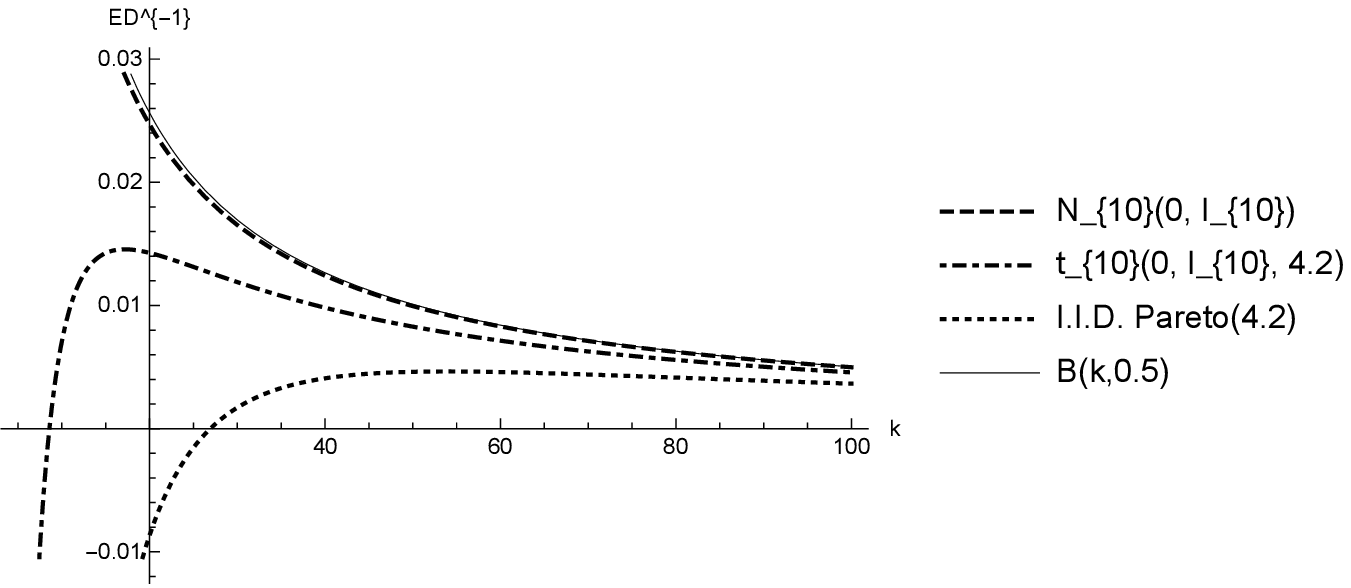}
\caption{$\overset{-1}{E\!D}$ when $\epsilon \sim N(0, 1)$}
\label{fig:normal_x_-1_1}
\end{figure}

\begin{figure}
\centering
\includegraphics[width=11cm,clip]{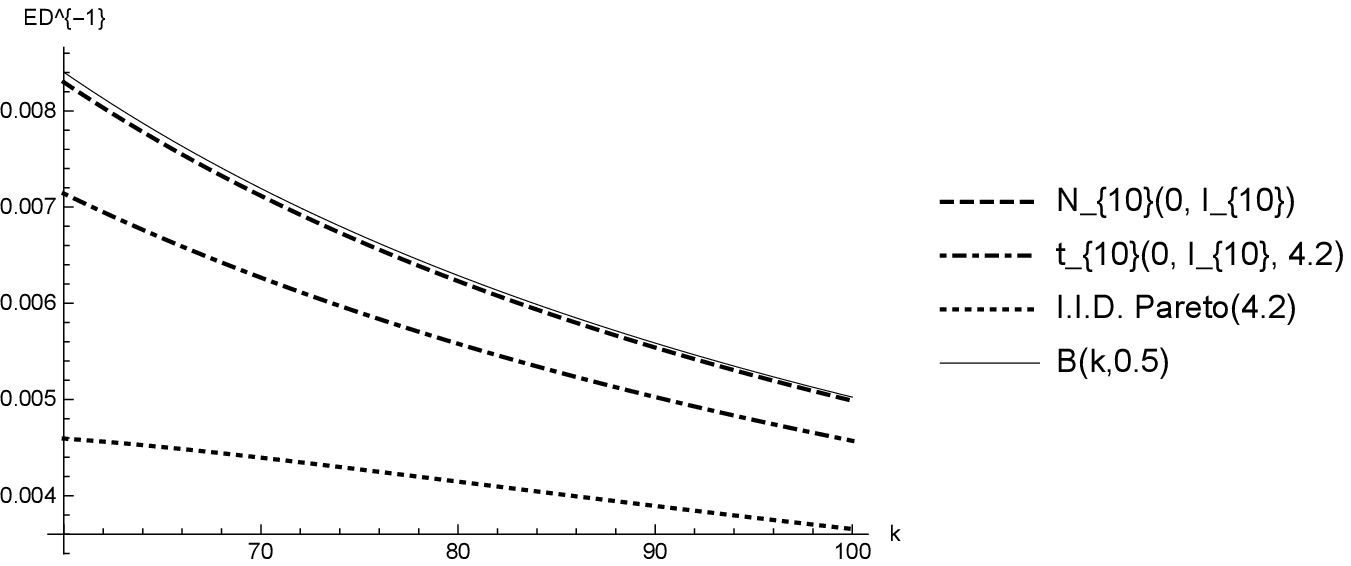}
\caption{$\overset{-1}{E\!D}$ when $\epsilon \sim N(0, 1)$}
\label{fig:normal_x_-1_2}
\end{figure}

We made a numerical comparison to see the effect of the joint moments of $x$. We set $p=10$ and $n=12k$, which means $p-n$ ratio equals $1/k$ since the number of the parameters of the regression model \eqref{regression_model} equals 12 when $p=10$. Figure \ref{fig:normal_x_-1_1} is the graph of the approximated $\overset{\scalebox{0.6}{\!\!$-1$}}{E\!D}$'s corresponding to each distribution of $x$ above-mentioned as $k$ varies from 5 to 100. (The graph for the controlled distribution is always quite similar to that for the normal distribution, hence for the clarity of the figures we will omit it in every figure hereafter.) Figure \ref{fig:normal_x_-1_2} magnifies the part $k=50, \ldots, 100.$  We put as the benchmark the approximated $\overset{\scalebox{0.6}{\!\!$-1$}}{E\!D}$ of the binomial model $B(k,0.5)$, that is, the $k$-times normal coin toss model. 
\begin{figure}
\centering
\includegraphics[width=11cm,clip]{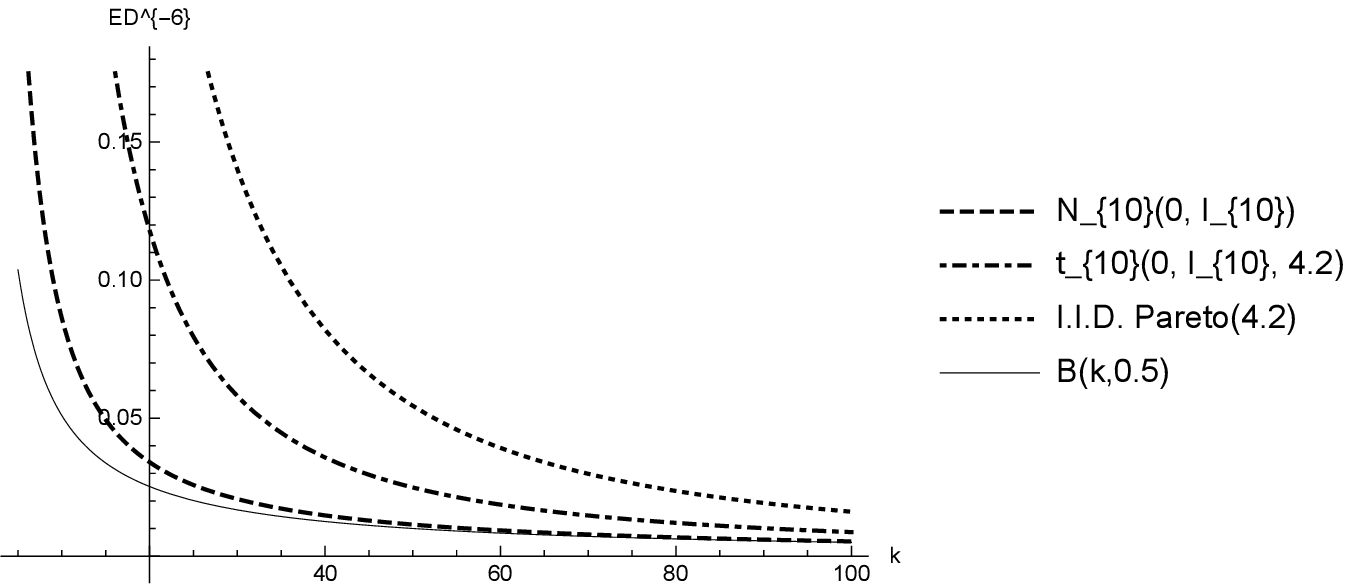}
\caption{$\overset{-6}{E\!D}$ when $\epsilon \sim N(0, 1)$}
\label{fig:normal_x_-6_1}
\end{figure}
\begin{table}
\centering
\caption{I.D.E. \& R.S.S. for $N(0,1)$ error distribution}
\label{ide_ssd_normal}
\begin{tabular}{|c|c|c|}
\hline 
      & I.D.E. &  R.S.S.\\
\hline
$x\sim N_{10}(0, I_{10})$ & *  & 111(10)\\
\hline
$x\sim t_{10}(0, I_{10}, 4.2)$& * & 322(40)\\
\hline
$x$ is controlled & * & 112(10) \\
\hline
$x$ is $i.i.d. P(4.2)$ &  *   & 741(110) \\
\hline
\end{tabular}
\end{table}

We notice that heavy tail property of Pareto distribution $P(4.2)$ or $t$-distribution $t(4.2)$ decreases difficulty in estimating the parameter, especially the large $m_4$ value 8129/21 of $P(4.2)$ makes the estimation easier. On the contrary, if $m_{4}$ and $m_{22}$ are as small as those of $N_{10}(0, I_{10})$ (or the controlled distribution), then the difficulty of estimation is close to the normal coin toss. 

Here we refer to the question how large sample size is required for the good approximation of $\overset{\alpha}{E\!D}$ by the expansion up to the $n^{-2}$ order term. It is very difficult to give a general answer to this question, but at least for a specific model, obviously we should not use the approximation unless it is positive or decreasing  with respect to $n$. For example, in Figure \ref {fig:normal_x_-1_1}, we see that the approximation for $t_{10}$ should be used for $k > 10$, namely $n>120$.

We observed that the effects of $m_{4}$ and $m_{22}$ depends on $\alpha$. If $\alpha$ is outside the interval \eqref{ED_normal_alpha_int}, the large value of $m_{44}$ or $m_{22}$ enhances the difficulty of the estimation. For example, if $\alpha=-6$, the order of various distributions of $x$ is  completely reversed to that for the case $\alpha=-1$ as we can see from Figure  \ref{fig:normal_x_-6_1}.

Now we consider I.D.E. and R.S.S. introduced in Section \ref{general_result}. We take Kullback-Leibler divergence ($\alpha=-1$)  as an example. Let $p$ be 10.  When $x$ is distributed as $N_p(0, I_p)$, we have 
$$
\overset{\scalebox{0.6}{\!\!$-1$}}{E\!D}(n)=
\frac{6}{n}-\frac{217}{12 n^2}.
$$
``Indicator of the difficulty of estimation'' is given as the solution of $m$ for the equation
$$
\frac{1}{2k}+\frac{1}{12k^{2}}(M-1)=\frac{1}{2k}-\frac{217}{12\times(12k)^2}
$$
(See \eqref{expan_KL_binom} for the left-hand side.)
Actually this quadratic equation of $m$  does not have the real roots. The left-hand side is always larger than the right-hand side. This means the estimation of the regression model is easier than the coin toss problem under the same $p-n$ ratio. 

Sample size determination is solving the next equation;
$$
\frac{1}{20}+\frac{1}{400}=\frac{6}{n}-\frac{217}{12 n^2}.
$$
where $n=111$ is  the rounded solution. 
For the other distributions (joint moments) of $x$, we can similarly calculate I.D.E. and R.S.S..The result is given in Table \ref{ide_ssd_normal}. "*" indicates that the equation has no solutions.  The number in the parenthesis in R.S.S. shows the sample size of the binomial model in the left-hand side of \eqref{equation_sample_size} (see the last paragraph of Section \ref{general_result}.) With the sample size given by R.S.S., the $p-n$ ratio of the regression model equals 12/R.S.S., while that of the coin toss model is equal to the reciprocal of the number in the parenthesis. Hence R.S.S divided by the number in the parenthesis could be another indicator. It is smaller for $t_{10}(0, I_{10}, 4.2)$ or $P(4.2)$ than that for $N_{10}(0, I_{10})$ or the controlled distribution. The  large joint moments $m_{4}$, $m_{22}$ for $t_{10}(0, I_{10},4.2)$ or $P(4.2)$ make estimation easier. We can guess that the large oscillation of $x$ is helpful to  estimate the values of $\beta$. Nevertheless of these differences, in general, the estimation for the regression model under the normal error distribution is not so troublesome, since 10 times as large sample size as the dimension of the parameter guarantees relatively easy estimation.
%
%
%
%
\subsection{$t$ Error Term Distribution}
\label{t-Distribution Error}
In this subsection we  investigate the case where $\epsilon$ has a $t$-distribution.  
Since 
$$
f(y)=c(\nu) (1+y^2/\nu)^{-(\nu+1)/2},\quad c(\nu)=\Gamma((\nu+1)/2)/(\sqrt{\pi \nu}\Gamma(\nu/2)),
$$
we have
\begin{align*}
\log f(y)&=\log c(\nu)-\frac{\nu+1}{2}\log(1+\frac{y^2}{\nu}), \\
\frac{d }{dy}\log f(y)&=-(\nu+1)\frac{y}{\nu+y^2},\\
\frac{d^2 }{dy^2}\log f(y)&=(\nu+1)\frac{y^2-\nu}{(\nu+y^2)^2},\\
\frac{d^3 }{dy^3}\log f(y)&=2(\nu+1)y \frac{3\nu-y^2}{(\nu+y^2)^3}.
\end{align*}
If we put $\tilde{y}$ as $\nu+y^2$, then
\begin{align*}
&\Bigl(\frac{d^3 }{dy^3}\log f(y)\Bigr)^i
\Bigl(\frac{d^2 }{dy^2}\log f(y)\Bigr)^j
\Bigl(\frac{d }{dy}\log f(y)\Bigr)^k
y^l\\
&=\Bigl((-1)^i 2^i (\nu+1)^i y^i (y^2-3\nu)^i \tilde{y}^{-3i}\Bigr)
\Bigl((\nu+1)^j (y^2-\nu)^j \tilde{y}^{-2j}\Bigr)
\Bigl((-1)^k (\nu+1)^k y^k \tilde{y}^{-k}\Bigr) y^l \\
&=2^i (-1)^{i+k} (\nu+1)^{i+j+k} y^{i+k+l} (y^2-3\nu)^i (y^2-\nu)^j \tilde{y}^{-(3i+2j+k)}\\
&=2^i (-1)^{i+k} (\nu+1)^{i+j+k} y^{i+k+l} \\
&\quad\times \Bigl(\sum_{s=0}^i y^{2s}(-3\nu)^{i-s} {}_i C_s\Bigr) 
\Bigl(\sum_{t=0}^j y^{2t}(-\nu)^{j-t} {}_j C_t \Bigr) \tilde{y}^{-(3i+2j+k)}\\
&=\sum_{s=0}^i \sum_{t=0}^j 2^i (-1)^{2i+j+k-s-t} (\nu+1)^{i+j+k} \nu^{i+j-s-t} y^{i+k+l+2s+2t}\,{}_i C_s\, {}_j C_t \,\tilde{y}^{-(3i+2j+k)}.
\end{align*}
Consequently 
\begin{align}
&\eta[i, j, k,, l]\nonumber\\
&=\sum_{s=0}^i \sum_{t=0}^j 2^i (-1)^{2i+j+k-s-t} (\nu+1)^{i+j+k} \nu^{-s-t-2i-j-k}\, {}_i C_s\, {}_j C_t \nonumber\\
&\quad \times c(\nu) H(i+k+l+2s+2t, \nu+6i+4j+2k) \label{eta_t},
\end{align}
where 
\begin{align*}
H(a, b)&\triangleq \int_{-\infty}^\infty y^{a} (1+y^2/\nu)^{-(b+1)/2} dy\\
&= \int_{-\infty}^\infty y^a \biggl(1+\frac{y^2(b/v)}{b}\biggr)^{-(b+1)/2} dy\\
&=\int_{-\infty}^\infty \tilde{y}^a \biggl( \frac{\nu}{b} \biggr)^{a/2} \biggl( 1+\frac{\tilde{y}^2}{b}\biggr)^{-(b+1)/2} \biggl(\frac{\nu}{b}\biggr)^{1/2}d\tilde{y},\qquad 
\tilde{y}\triangleq y\sqrt{\frac{b}{v}}\\
&=\biggl(\frac{\nu}{b}\biggr)^{(a+1)/2} \int_{-\infty}^\infty {\tilde{y}}^a \biggl( 1+\frac{\tilde{y}^2}{b}\biggr)^{-(b+1)/2} d\tilde{y}\\
&=\biggl(\frac{\nu}{b}\biggr)^{(a+1)/2} \frac{1}{c(b)} E[t^a],\quad t\sim t(b)\\
&=
\begin{cases}
0, &\text{ if $0< a < b$, and $a$ is an odd number,} \\
\nu^{(a+1)/2}\:\frac{\Gamma((a+1)/2)\Gamma((b-a)/2)}{\Gamma((b+1)/2)}, &\text{ if  $0< a < b$, and  $a$ is a even number.} 
\end{cases}
\end{align*}
If we insert these results, we get $\overset{\alpha}{E\!D}$ as a function of $p, n, \alpha, \nu$ and the joint moments of $x$. We are interested in how $\overset{\alpha}{E\!D}$ is effected by a long tail error distribution compared to the standard normal distribution, hence we set $\nu$ as 3. (Mathematica program for a general $\nu$ is available in Appendix B) We have the following result;
\begin{align}
&\overset{\alpha}{E\!D} \nonumber\\
&=\frac{p+2}{2n}\nonumber\\
&\qquad+\frac{1}{384 n^2}\Bigl(6 \alpha^2 \left(13 + (10 - 3 m_{22} + 3 m_{4}) p + 3 m_{22} p^2\right)\nonumber\\
&\hspace{30mm}-2\alpha\left(-77 + (-72 - 51 m_{22} + 51 m_{4}) p + 3 (-5 + 17 m_{22}) p^2\right)\nonumber\\
&\hspace{30mm}+3\left(287 + (296 + 90 m_{22} - 90 m_{4}) p + (65 - 90 m_{22}) p^2\right)\Bigr)\nonumber\\
&\qquad+o(n^{-2})\label{ED_t_1}\\
&=\frac{p+2}{2 n}\nonumber\\
&\quad+\frac{1}{384 n^2}\Bigl(6p(-45 - 17\alpha + 3 \alpha^2) m_{4}+6p(p-1)(-45 - 17\alpha + 3 \alpha^2)m_{22}\nonumber\\
&\hspace{25mm}+\alpha^2 (78+ 60 p) + \alpha(154 + 144 p + 30 p^2) \nonumber\\
&\hspace{25mm}+861+ 888 p + 195 p^2  \Bigr)\nonumber\\
&\qquad+o(n^{-2})\label{ED_t_2}\\
&=\frac{p+2}{2 n}\nonumber\\
&\quad+\frac{1}{384 n^2}\Bigl( (195 + 30 \alpha - 270 m_{22} - 102 \alpha m_{22} + 18 \alpha^2 m_{22}) p^2\nonumber\\
&\hspace{25mm}+ (888 + 144 \alpha + 60 \alpha^2 + 270 m_{22} + 102 \alpha m_{22} - 18 \alpha^2 m_{22} \nonumber\\
&\hspace{25mm}\qquad -270 m_{4} - 102 \alpha m_{4} + 18 \alpha^2 m_{4}) p\nonumber\\
&\hspace{25mm}+861 + 154 \alpha + 78 \alpha^2\Bigr)\nonumber\\
&\qquad+o(n^{-2})\label{ED_t_3}
\end{align}
The $n^{-2}$ order term has similar properties as in the case of $N(0,1)$.
\begin{enumerate}
\item The dimension of  $p$ is two, hence $\overset{\alpha}{E\!D}$ is asymptotically determined by the $p-n$ ratio.
\item  Other moments than $m_4$ and $m_{22}$ do not appear.
\item The coefficients of $m_4$ and $m_{22}$ are non-positive when
 $3 \alpha^2 - 17\alpha -45 \leq 0$, that is, 
\begin{equation}
\label{ED_t_alpha_int}
-1.97\cdots \leq \alpha \leq 7.63\cdots.
\end{equation}
For $\alpha$ within this interval, the larger $m_{22}$ or $m_4$ gets, the less $\overset{\alpha}{E\!D}$ becomes. The divergences often used in statistical literature are all included in this interval.
\end{enumerate}
We noticed that if the error term distribution is the standard normal distribution (see \eqref{ED_normal_2}) or $t(3)$ distribution (see \eqref{ED_t_2}), only $m_{4}$ and $m_{22}$ among the moments \eqref{4th_moments} and \eqref{3rd_moments} appear in the asymptotic expansion of $\overset{\alpha}{E\!D}$ up to the $n^{-2}$ order. 
On this phenomena, we have the following general result.
\begin{prop} 
\label{disappear_m^6}
If the error term distribution is quadratic, namely $f(\epsilon)=g(\epsilon^2)$ for some function $g(\cdot)$, then the asymptotic expansion of $\overset{\alpha}{E\!D}$ up to the order $n^{-2}$ includes only $m_{4}$ and $m_{22}$ among the third and forth order joint moments of $x$.
\end{prop}
\textit{<Proof>} From \eqref{def_L}, we notice that the third or forth order moments of $x$  in the expansion of $\overset{\alpha}{E\!D}$ up to the order $n^{-2}$  are generated from  the terms $m[i,j,k,l]$ and $m[i,j,k]m[s,t,u]$.

The forth order moments arise from $m[i, j, k, l] (1\leq i, j, k, l \leq p)$ in $L11$ to $L15$. 
Since $m[i,j,k,l]$ is multiplied with $g^{ij}g^{kl}$ as in \eqref{def_L}, and $g^{ij}$ vanishes unless $i=j$, the possible moments coming from  $E[x_i x_j x_k x_l]$ are only $m_{4}$ and $m_{22}$.

On the other hand, $m[i,j,k]m[s,t,u]$ come from either term of $L21,\ldots, L26$.
We notice that  if the third moments are generated from these terms, they are always multiplied 
 with $\eta[0,1,1,0]$ or $\eta[0,0,3,0]$. (See \eqref{L_(ij)k}, \eqref{L_ijk}.)
If $f(\epsilon)=g(\epsilon^2)$, then we have
\begin{eqnarray*}
&\frac{d }{dy}\log f(y)=2y\frac{g'(y^2)}{g(y^2)},\\
&\frac{d^2 }{dy^2}\log f(y)=\frac{2}{g^2(y^2)}\Bigl(g'(y^2)+2y^2g''(y^2)g(y^2)+2y^2(g'(y^2))^2\Bigr).\\
\end{eqnarray*}
Therefore $\eta[0,1,1,0]$ and $\eta[0,0,3,0]$ vanishes.  \hfill \textit{ Q.E.D.}
\\
Typical four cases $\alpha=-1,1,0,3$ are given as follows;
\begin{align}
&\overset{\scalebox{0.6}{\!\!$-1$}}{E\!D}=\frac{p+2}{2n}\nonumber\\
&\qquad\quad+\frac{1}{384 n^2}\Bigl(785 + 6 (134 + 25 m_{22} - 25 m_{4}) p + (165 - 150 m_{22}) p^2\Bigr)+o(n^{-2}) ,\label{ED_-1_t}\\
&\overset{\scalebox{0.6}{\!\!$1$}}{E\!D}=\frac{p+2}{2n}\nonumber\\
&\qquad\quad+\frac{1}{384 n^2}\Bigl(1093 + 6 (182 + 59 m_{22} - 59 m_{4}) p + (225 - 354 m_{22}) p^2\Bigr)+o(n^{-2}) ,\label{ED_1_t}\\
&\overset{\scalebox{0.6}{\!\!$0$}}{E\!D}=\frac{p+2}{2n}\nonumber\\
&\qquad\quad+\frac{1}{384 n^2}\Bigl(287 + (296 + 90 m_{22} - 90 m_{4}) p + (65 - 90 m_{22}) p^2\Bigr)+o(n^{-2}) ,\label{ED_0_t}\\
&\overset{\scalebox{0.6}{\!\!$3$}}{E\!D}=\frac{p+2}{2n}\nonumber\\
&\qquad\quad+\frac{1}{128 n^2}\Bigl(675 + 2 (310 + 69 m_{22} - 69 m_{4}) p + (95 - 138 m_{22}) p^2\Bigr)+o(n^{-2}) .\label{ED_2_t}
\end{align}
For each distribution (moments) of $x$ introduced in Section \ref{general_result}, we have the following results.
\begin{enumerate}
\item  $x\sim N_p(0, I_p)$
\begin{align}
&\overset{\alpha}{E\!D} \nonumber\\
&=\frac{p+2}{2 n}\nonumber\\
&\quad+\frac{1}{384 n^2}\Bigl(6 \alpha^2 (13 + 16 p + 3 p^2) - 2 \alpha (-77 + 30 p + 36 p^2)\nonumber\\
&\qquad\qquad\qquad 861 + 348 p - 75 p^2\Bigr)+o(n^{-2})
\label{tncase_ED}
\end{align}
\item $x \sim t_p(0, I_p, 4.2)$
\begin{align}
&\overset{\alpha}{E\!D} \nonumber\\
&=\frac{p+2}{2n}\nonumber\\
&\quad+\frac{1}{384 n^2}\Bigl(6 \alpha^2 (13 + 76 p + 33 p^2) - 14 \alpha (-11 + 150 p + 78 p^2)\nonumber\\
&\qquad\qquad\qquad 3 (287 - 1684 p - 925 p^2)\Bigr)+o(n^{-2})
\label{ttcase_ED}
\end{align}
\item $x$ is controlled.
\begin{align}
&\overset{\alpha}{E\!D} \nonumber\\
&=\frac{p+2}{2n}\nonumber\\
&\quad+\frac{1}{384 n^2}\Bigl( 6 \alpha^2 (13 + 10 p + 3 p^2)+\alpha (154 + 144 p - 72 p^2) \nonumber\\
&\qquad\qquad\qquad 861 + 888 p - 75 p^2\Bigr)+o(n^{-2})
\label{tccase_ED}
\end{align}
\item $x_i$ is i.i.d. as $P(4.2)$ $(i=1,\ldots,p)$
\begin{align}
&\overset{\alpha}{E\!D} \nonumber\\
&=\frac{p+2}{2n}\nonumber\\
&\quad+\frac{1}{2688 n^2}\Bigl(6 \alpha^2 (91 + 8178 p + 21 p^2) - 2 \alpha (-539 + 137332 p + 252 p^2)\nonumber\\
&\qquad\qquad\qquad 3 (2009 - 241168 p - 175 p^2)\Bigr)+o(n^{-2})
\label{tpcase_ED}
\end{align}
\end{enumerate}

\begin{figure}
\centering
\includegraphics[width=11cm,clip]{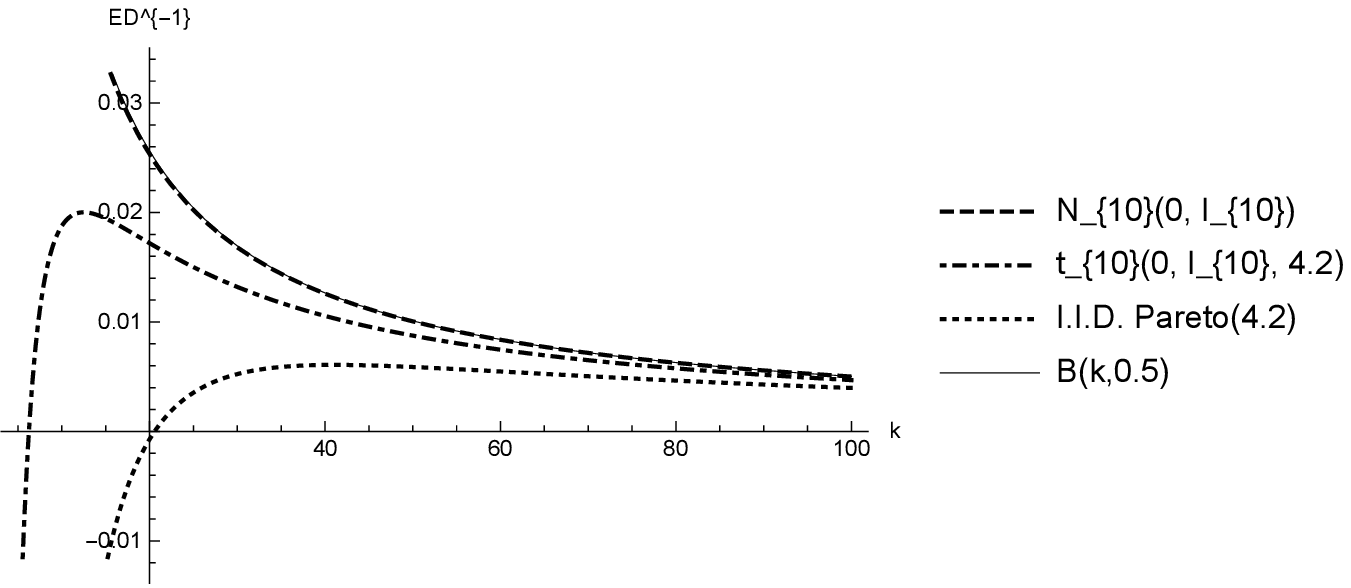}
\caption{$\overset{-1}{E\!D}$ when $\epsilon \sim t(3)$}
\label{fig:t_x_-1_1}
\end{figure}

\begin{figure}
\centering
\includegraphics[width=11cm,clip]{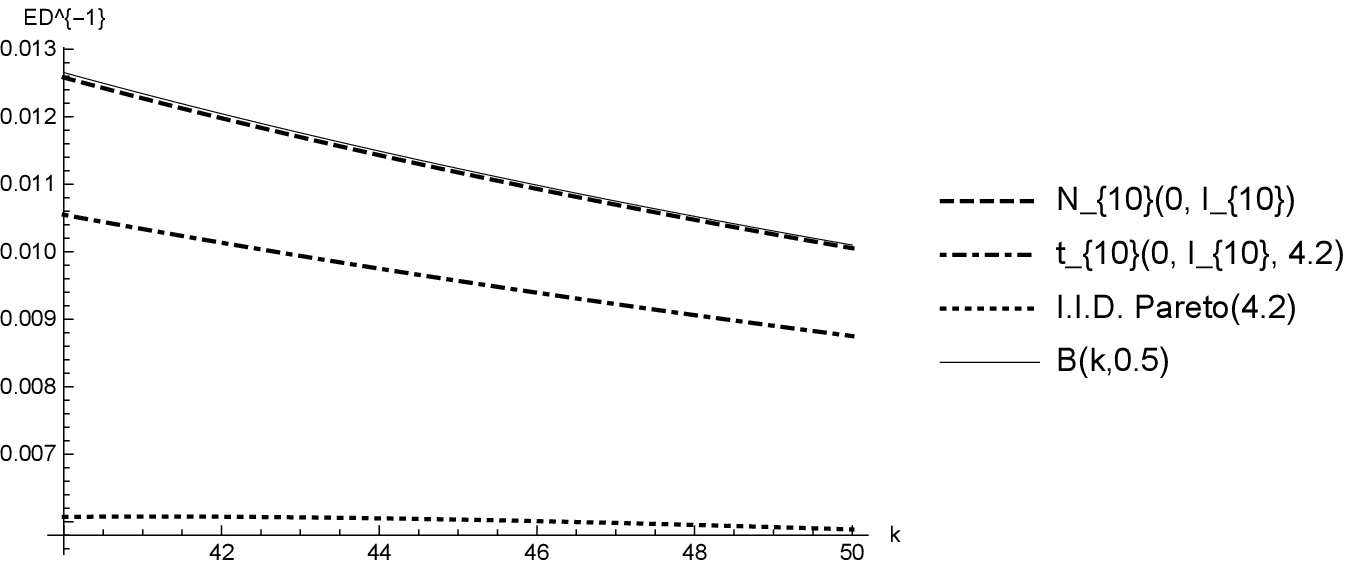}
\caption{$\overset{-1}{E\!D}$ when $\epsilon \sim t(3)$}
\label{fig:t_x_-1_2}
\end{figure}

We made a numerical comparison under the condition $p=10$ and $n=12k$, which means $p-n$ ratio equals $1/k$. Figure \ref{fig:t_x_-1_1} is the graph of the approximated $\overset{\scalebox{0.6}{\!\!$-1$}}{E\!D}$'s corresponding to each distribution above-mentioned  except for the controlled distribution as $k$ varies from 5 to 100. Figure \ref{fig:t_x_-1_2} magnifies the part of $k$ from 40 to 50. We put as the benchmark the approximated $\overset{\scalebox{0.6}{\!\!$-1$}}{E\!D}$ of the binomial model $B(k,0.5)$.

Just like the case of $N(0,1)$,  heavy tail property of Pareto distribution $P(4.2)$ or $t$-distribution $t(4.2)$ eases difficulty in estimating the parameter. On the contrary, if $m_{4}$ and $m_{22}$ are as small as those of $N(0, 1)$ (or the controlled distribution), then the difficulty of estimation is close to the normal coin toss. 

\begin{figure}
\centering
\includegraphics[width=11cm,clip]{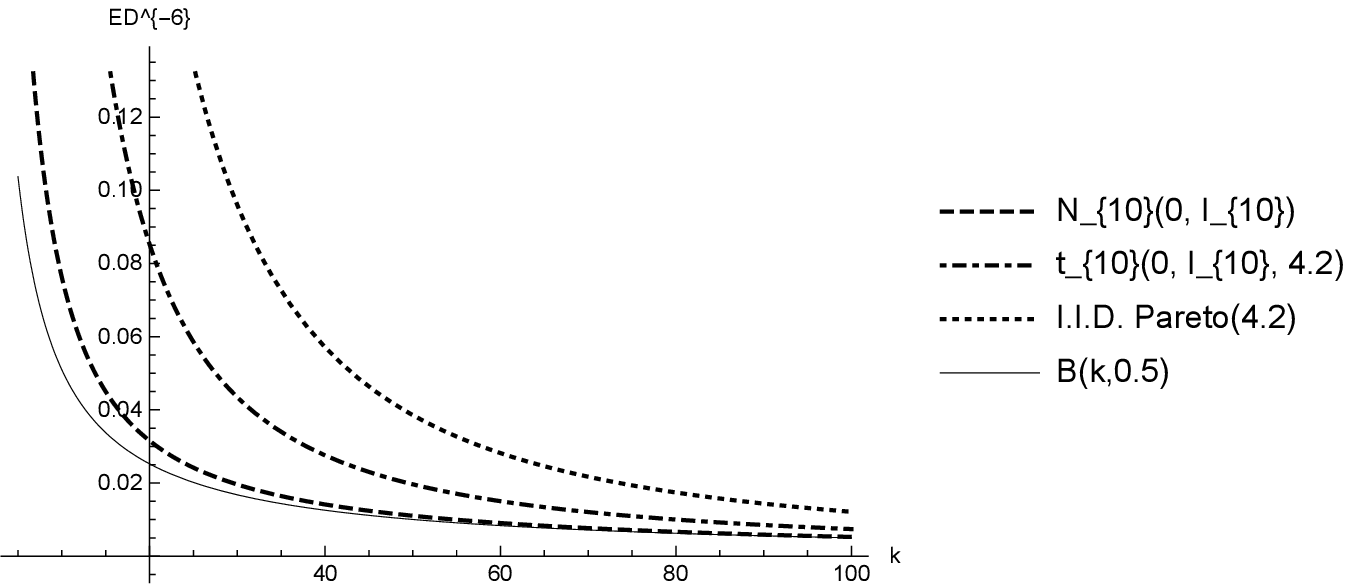}
\caption{$\overset{-6}{E\!D}$ when $\epsilon \sim t(3)$}
\label{fig:t_x_-6}
\end{figure}

It was also observed that the effects of $m_{4}$ and $m_{22}$ depends on $\alpha$. If $\alpha$ is outside the interval \eqref{ED_t_alpha_int}, the large value of $m_{44}$ or $m_{22}$ enhances the difficulty of the estimation. For example, see Figure \ref{fig:t_x_-6} for $\alpha=-6$, where the $\overset{-6}{E\!D}$ for the $t$ or Pareto distribution is larger than that of the normal distribution.

We considered I.D.E. and R.S.S. w.r.t. Kullback-Leibler divergence under the condition $p=10$ for each distribution of $x$. Table \ref{ide_ssd_t} shows the result.The same comments hold as in the case of $N(0,1)$. The large value of $m_{4}$ or $m_{22}$ of the $t$-distribution or Pareto distribution makes the estimation easier compared to the normal distribution or the controlled distribution. Generally speaking, irrespective of the above difference, the estimation of the regression model under $t$-distribution error is not so hard. With 10 times as large sample size as the parameter dimension, we can estimate the parameter without much trouble.
\begin{table}
\centering
\caption{I.D.E. \& R.S.S. for $t(3)$ error distribution}
\label{ide_ssd_t}
\begin{tabular}{|c|c|c|}
\hline 
      & I.D.E. &  R.S.S.\\
\hline
$x\sim N_{10}(0, I_{10})$ & * & 117(10)\\
\hline
$x\sim t_{10}(0, I_{10}, 4.2)$& * & 246(30)\\
\hline
$x$ is controlled & * & 118(10) \\
\hline
$x$ is $i.i.d. P(4.2)$ &  *   & 689(90) \\
\hline
\end{tabular}
\end{table}
\subsection{Skew-Normal Error Term Distribution}
\label{Skew_Normal_Distribution_Error}
In this subsection, we take a skew-normal distribution as an error term distribution so that we investigate how the skewness of error effects $\overset{\alpha}{E\!D}.$
Suppose that $f(y)$ is given by 
$$
f(y)=2\phi(y) \Phi(b y),
$$
where $\phi(y)$ is the p.d.f. of the standard normal distribution, and $\Phi(y)$ is its cumulative distribution function. This is the p.d.f. of the (standard) skew-normal distribution with the shape parameter $b$. If $b$ is positive (negative), the distribution is right (left) skewed. When $b=0$, it is the standard normal distribution.

\begin{figure}
\centering
\includegraphics[width=11cm,clip]{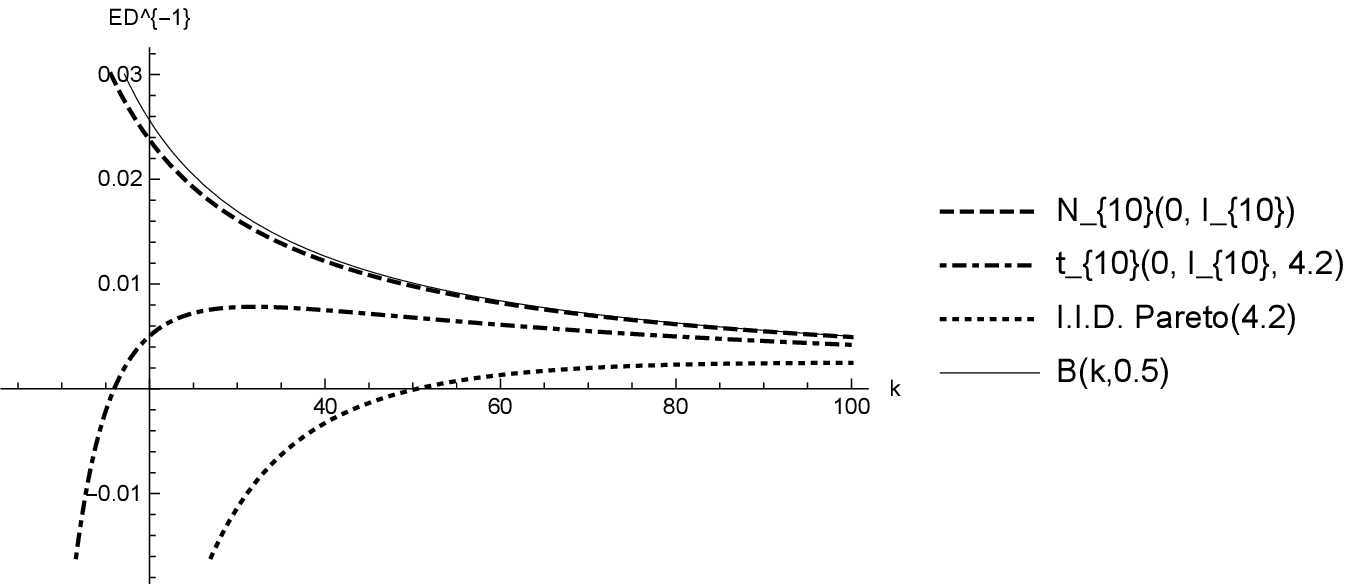}
\caption{$\overset{-1}{E\!D}$ when $\epsilon \sim SN(3)$}
\label{fig:sn_x_-1_1}
\end{figure}

\begin{figure}
\centering
\includegraphics[width=11cm,clip]{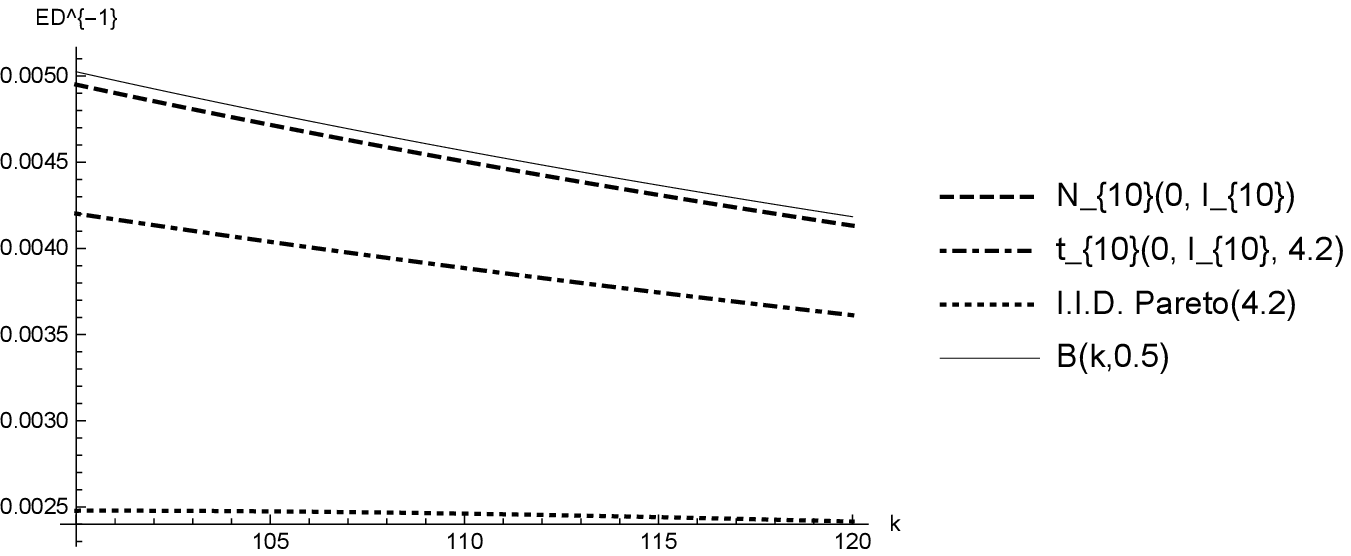}
\caption{$\overset{-1}{E\!D}$ when $\epsilon \sim SN(3)$}
\label{fig:sn_x_-1_2}
\end{figure}

For this distribution, we have
\begin{align*}
\log f(y)&=\log 2+\log\phi(y)+\log\Phi(by), \\
\frac{d }{dy}\log f(y)&=-y+\frac{b \phi(by)}{\Phi(by)},\\
\frac{d^2 }{dy^2}\log f(y)&=-1-b^2\frac{\phi^2(by)}{\Phi^2(by)}+\frac{b^2\phi'(by)}{\Phi(by)}=-1-b^2\frac{\phi^2(by)}{\Phi^2(by)}-b^3y\frac{\phi(by)}{\Phi(by)},\\
\frac{d^3 }{dy^3}\log f(y)&=2b^3\frac{\phi^3(by)}{\Phi^3(by)}-2b^3\frac{\phi(by)\phi'(by)}{\Phi^2(by)}-b^3\frac{\phi(by)\phi'(by)}{\Phi^2(by)}+b^3\frac{\phi''(by)}{\Phi(by)}\\
&=2b^3\frac{\phi^3(by)}{\Phi^3(by)}-3b^3\frac{\phi(by)\phi'(by)}{\Phi^2(by)}+b^3\frac{\phi''(by)}{\Phi(by)}\\
&=b^3\biggl(2\frac{\phi^3(by)}{\Phi^3(by)}+3by\frac{\phi^2(by)}{\Phi^2(by)}+(b^2y^2-1)\frac{\phi(by)}{\Phi(by)}\biggr)
\end{align*}
If we insert these results into the definition \eqref{def_eta}, we get the formal form of  $\eta[i, j, k, l]$. However, since the value of $\eta[i, j, k, l]$ can not be gained theoretically, we have to calculate it numerically after a specific value of $b$ is chosen.  Here we put $b=3$ for a relatively strong right skewness. The asymptotic expansion of $\overset{\alpha}{E\!D}$ is given as follows (the numbers are rounded off to three decimal place).
\begin{align}
&\overset{\alpha}{E\!D}\nonumber\\
&=\frac{p+2}{2n}\nonumber\\
&\quad+\frac{1}{n^2}\Bigl(\alpha^2\bigl(0.175 m_{22}p^2+(0.175m_4-0.175m_{22}+0.689)p+0.988\bigr) \nonumber\\
&\qquad\qquad+\alpha\bigl((0.087 m_{111}^2 + 0.130 m_{21}^2) p^3\nonumber\\
&\qquad\qquad\qquad+(0.385 - 0.260 m_{111}^2 - 0.430 m_{22} + 0.260 m_{21} m_3) p^2
\nonumber\\
&\qquad\qquad\qquad+(2.074 + 0.174 m_{111}^2 - 0.130 m_{21}^2 + 0.430 m_{22} - 
   0.260 m_{21} m_3\nonumber\\
&\hspace{30mm} + 0.217m_3^2 - 0.430 m4) p+2.352\bigr)\nonumber\\
&\qquad\qquad + (0.065 m_{21}^2-0.130 m_{111}^2) p^3\nonumber\\
&\qquad\qquad+(0.583+ 0.392 m_{111}^2 - 0.522 m_{21}^2 - 
 1.503 m_{22}  + 0.130 m_{21} m3) p^2\nonumber\\
&\qquad\qquad+(2.823 - 0.261 m_{111}^2+ 0.457 m_{21}^2 + 1.503 m_{22} - 0.130 m_{21} m_3\nonumber\\
&\qquad\qquad\qquad - 0.065 m_3^2-1.503m_4)p+3.385\Bigr)\nonumber\\
&\quad +o(n^{-2})\label{ED_s_1}\\
&=\frac{p+2}{2n}\nonumber\\
&\quad+\frac{1}{n^2}\Bigl(0.175 p(-8.570- 2.451\alpha + \alpha^2) m_4\nonumber\\
&\qquad\qquad+0.175 p(p-1)(-8.570- 2.451\alpha + \alpha^2)m_{22}\nonumber\\
&\qquad\qquad+0.217p (-0.302+\alpha)m_{3}^2\nonumber\\
&\qquad\qquad+p(0.065 p^2 + 0.130 \alpha p^2 -0.522 p+0.457 - 0.130 \alpha)m_{21}^2\nonumber\\
&\qquad\qquad+0.087p(-1.504 p^2 + \alpha p^2+4.513 p - 3\alpha p-3.001+2\alpha)m_{111}^2\nonumber\\
&\qquad\qquad+0.260p(p-1)(0.500+\alpha) m_3m_{21}\nonumber\\
&\qquad\qquad+(0.988+0.689p)\alpha^2 + (2.352+2.074p+0.385p^2) \alpha \nonumber\\
&\qquad\qquad+3.385+2.823p+0.583p^2\Bigr)\nonumber\\
&\quad+o(n^{-2})\label{ED_s_2}\\
&=\frac{p+2}{2n}\nonumber\\
&\quad+\frac{1}{n^2}\Bigl(
p^3 \bigl(-0.131 m_{111}^2 + 0.065 m_{21}^2  + \alpha (0.087 m_{111}^2 + 0.130 m_{21}^2)\bigr)\nonumber\\
&\qquad\qquad+p^2 \bigl( \alpha^2 0.175 m_{22}+\alpha (0.385 - 0.260 m_{111}^2 - 0.430 m_{22} + 
 0.260 m_{21} m_3)\nonumber\\
&\qquad\qquad\qquad+0.583 + 0.392 m_{111}^2 - 0.522 m_{21}^2 - 1.503 m_{22} + 
 0.130 m_{21} m_3\bigr)\nonumber\\
&\qquad\qquad+p
\bigl(\alpha^2( 0.689- 0.175 m_{22} + 0.175 m_4)\nonumber\\
&\qquad\qquad\qquad+\alpha(2.074 + 0.174 m_{111}^2 - 0.130 m_{21}^2 + 0.430 m_{22}\nonumber\\
&\qquad\qquad\qquad\qquad - 0.260 m_{21} m_3 + 0.217 m_3^2 - 0.430 m_4)\nonumber\\
&\qquad\qquad\qquad+2.823 - 0.261 m_{111}^2 + 0.457 m_{21}^2 + 1.503 m_{22} - 
 0.130 m_{21} m_3 - 0.065m_3^2 \bigr) \nonumber\\
&\qquad\qquad+3.385 + 2.352 \alpha + 0.988 \alpha^2\Bigr)\nonumber\\
&\quad +o(n^{-2})\label{ED_s_3}
\end{align}
Typical four cases $\alpha=-1,1,0,3$ are given as follows;
\begin{align}
&\overset{\scalebox{0.6}{\!\!$-1$}}{E\!D}\nonumber\\
& =\frac{p+2}{2n}+\frac{1}{2n^2}\Bigl(0.435 m_{111}^2 p^3 - 0.130 m_{21}^2 p^3\nonumber\\
&\qquad 0.398 p^2 + 1.304 m_{111}^2 p^2 - 1.043 m_{21}^2 p^2 - 
 1.796 m_{22} p^2 - 0.261 m_{21} m_3 p^2\nonumber\\
&\qquad  + 2.876 p - 0.869 m_{111}^2 p + 1.173 m_{21}^2 p + 
 1.796 m_{22} p + 0.261 m_{21} m_3 p \nonumber\\
&\qquad- 0.565 m_3^2 p - 1.796 m_4 p +4.043\Bigr)+o(n^{-2}),\label{ED_-1_s}\\
&\overset{\scalebox{0.6}{\!\!$1$}}{E\!D}\nonumber\\
& =\frac{p+2}{2n}+\frac{1}{2n^2}\Bigl(-0.088 m_{111}^2 p^3 +0.391 m_{21}^2 p^3\nonumber\\
&\qquad 1.936 p^2 + 0.263 m_{111}^2 p^2 - 1.044 m_{21}^2 p^2 - 
 3.516 m_{22} p^2 +0.781 m_{21} m_3 p^2+\nonumber\\
&\qquad  + 11.171 p - 0.175 m_{111}^2 p + 0.653 m_{21}^2 p + 
 3.516 m_{22} p -0.781 m_{21} m_3 p \nonumber\\
&\qquad+0.303 m_3^2 p - 3.516 m_4 p +13.450\Bigr)+o(n^{-2}),\label{ED_1_s}\\
&\overset{\scalebox{0.6}{\!\!$0$}}{E\!D}\nonumber\\
& =\frac{p+2}{2n}+\frac{1}{2n^2}\Bigl(0.261 m_{111}^2 p^3 + 0.130 m_{21}^2 p^3\nonumber\\
&\qquad 1.167 p^2 + 0.783 m_{111}^2 p^2 - 1.043 m_{21}^2 p^2 - 
 3.007 m_{22} p^2 +0.260 m_{21} m_3 p^2\nonumber\\
&\qquad  + 5.645 p - 0.522 m_{111}^2 p + 0.913 m_{21}^2 p + 
 3.007 m_{22} p - 0.260 m_{21} m_3 p \nonumber\\
&\qquad- 0.131 m_3^2 p - 3.007 m_4 p +6.771\Bigr)+o(n^{-2}),\label{ED_0_s}\\
&\overset{\scalebox{0.6}{\!\!$3$}}{E\!D}\nonumber\\
& =\frac{p+2}{2n}+\frac{1}{2n^2}\Bigl(0.260 m_{111}^2 p^3 +0.911 m_{21}^2 p^3\nonumber\\
&\qquad 3.474 p^2 -0.778 m_{111}^2 p^2 - 1.044 m_{21}^2 p^2 - 
 2.430 m_{22} p^2 +1.823 m_{21} m_3 p^2\nonumber\\
&\qquad  + 30.488 p +0.519 m_{111}^2 p + 0.133 m_{21}^2 p + 
 2.430 m_{22} p -1.823 m_{21} m_3 p \nonumber\\
&\qquad+1.171 m_3^2 p -2.430 m_4 p +38.661\Bigr)+o(n^{-2}),\label{ED_2_s}
\end{align}
We observe the following points for $n^{-2}$ order term.
\begin{enumerate}
\item The dimension of  $p$ is three, hence if $p$ increases with a constant $p-n$ ratio, $n^{-2}$ order term could diverge for some given $\alpha$ and the moments of $x$. Then it is not enough to increase the sample size proportionally to the number of the explanatory variables in order to keep $\overset{\alpha}{E\!D}$ at a certain level.
\item  $m_{3}$, $m_{21}$ and $m_{111}$ appear in the expansion that do not appear in the case of $N_p(0, I_p)$ or $t_p(0, I_p, \nu).$  The effect of these moments are rather complicated and depends on $\alpha$ and $p$, For example, when $p$ is large enough, the larger absolute value of $m_{21}$ decreases the approximated $\overset{\alpha}{E\!D}$ for $\alpha=-1$, but vice versa for $\alpha=1, 0, 3$.
\item The larger $m_4$ and $m_{22}$ decreases $\overset{\alpha}{E\!D}$ if $\alpha^2- 2.451 \alpha-8.570 <0$, namely
$$
-1.95\cdots < \alpha < 4.40 \cdots.
$$
Again $\alpha$'s such as $-1, 0, 1, 3$  are all included  in this interval.
\end{enumerate}
For each distribution of $x$ introduced in Section \ref{general_result}, we have the following results.
\begin{enumerate}
\item  $x\sim N_p(0, I_p)$ \\
Since $m_4=3, m_{22}=1, m_3=0,\ m_{21}=0,\ m_{111}=0$, we have 
\begin{align}
&\overset{\alpha}{E\!D} \nonumber\\
&=\frac{p+2}{2 n}\nonumber\\
&\quad+\frac{1}{n^2}\Bigl( \alpha^2 (0.988 + 1.040 p + 0.175 p^2) + 
    \alpha (2.352 + 1.214 p - 0.046 p^2) \nonumber\\
&\qquad\qquad\qquad 3.385 - 0.184 p - 0.920 p^2\Bigr)+o(n^{-2})
\label{sncase_ED}
\end{align}
\item $x \sim t_p(0, I_p, 4.2)$\\
Since  $m_4=33,\ m_{22}=11,\ m_3=0,\ m_{21}=0,\ m_{111}=0$, we have
\begin{align}
&\overset{\alpha}{E\!D} \nonumber\\
&=\frac{p+2}{2n}\nonumber\\
&\quad+\frac{1}{n^2}\Bigl(\alpha^2 (0.988 + 4.548 p + 1.930 p^2)+\alpha (2.352 - 7.387 p - 4.346 p^2)\nonumber\\
&\qquad\qquad\qquad 3.385 - 30.253 p - 15.955 p^2\Bigr)+o(n^{-2})
\label{stcase_ED}
\end{align}
\item $x$ is controlled. \\
Since $m_4=1, m_{22}=1, m_3=0, m_{211}=0, m_{111}=0$, we have
\begin{align}
&\overset{\alpha}{E\!D} \nonumber\\
&=\frac{p+2}{2n}\nonumber\\
&\quad+\frac{1}{n^2}\Bigl(\alpha^2(0.988 + 0.689 p + 0.175 p^2)+\alpha(2.352 + 2.074 p - 0.046 p^2)\nonumber\\
&\qquad\qquad\qquad 3.385 + 2.823 p - 0.920 p^2\Bigr)+o(n^{-2})
\label{sccase_ED}
\end{align}
\item $x_i$ is i.i.d. as $P(4.2)$ $(i=1,\ldots,p)$\\
Since  $m_4=387.095, m_{22}=1, m_3=6.272, m_{21}=0, m_{111}=0$, we have
\begin{align}
&\overset{\alpha}{E\!D} \nonumber\\
&=\frac{p+2}{2n}\nonumber\\
&\quad+\frac{1}{n^2}\Bigl( \alpha^2 (0.988+ 68.420 p + 0.175 p^2)+\alpha(2.352 - 155.433 p - 0.046 p^2)\nonumber\\
&\qquad\qquad\qquad 3.385 - 580.229 p - 0.920 p^2\Bigr)+o(n^{-2})
\label{spcase_ED}
\end{align}
\end{enumerate}
We made a numerical comparison under the condition $p=10$ and $n=12k$, which means $p-n$ ratio equals $1/k$. Figure \ref{fig:sn_x_-1_1} is the graph of the approximated $\overset{\scalebox{0.6}{\!\!$-1$}}{E\!D}$'s corresponding to each distribution above-mentioned except for the controlled distribution as $k$ varies from 5 to 100.  Figure \ref{fig:sn_x_-1_2} magnifies the part of $k$ from 100 to 120. We put as the benchmark the approximated $\overset{\scalebox{0.6}{\!\!$-1$}}{E\!D}$ of the binomial model $B(k,0.5)$. 
\begin{figure}
\centering
\includegraphics[width=11cm,clip]{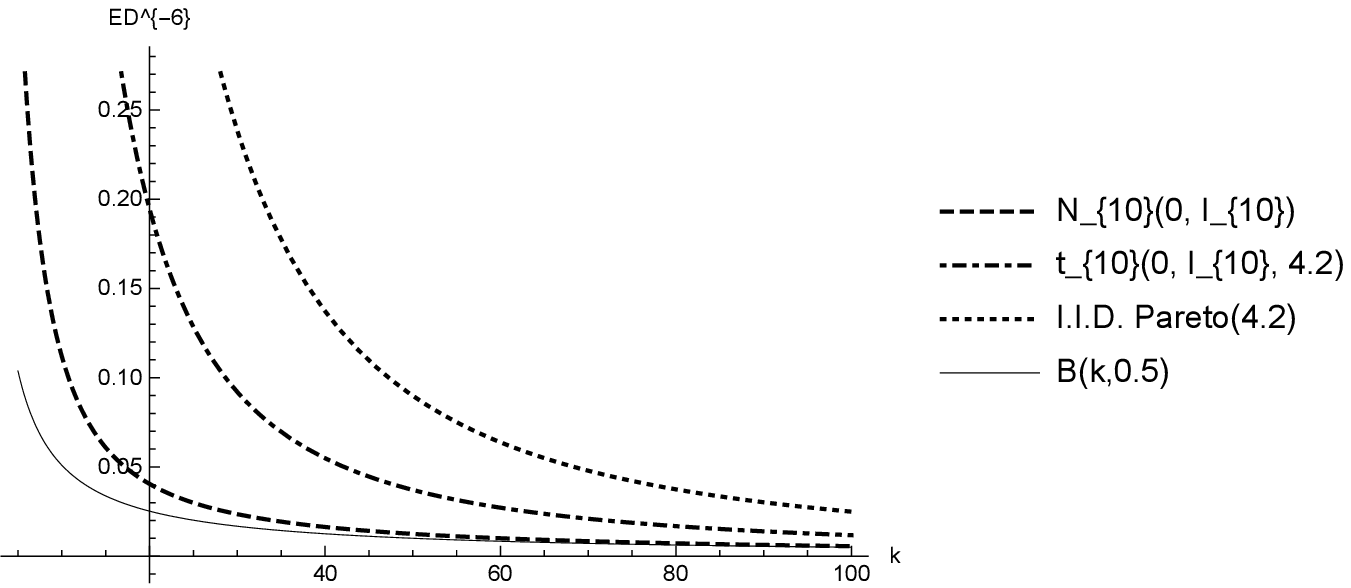}
\caption{$\overset{-6}{E\!D}$ when $\epsilon \sim SN(3)$}
\label{fig:sn_x_-6}
\end{figure}
\begin{table}
\centering
\caption{I.D.E. \& R.S.S. for $SN(3)$ error distribution}
\label{ide_ssd_sn}
\begin{tabular}{|c|c|c|}
\hline 
      & I.D.E. &  R.S.S.\\
\hline
$x\sim N_{10}(0, I_{10})$ & * & 101(10)\\
\hline
$x\sim t_{10}(0, I_{10}, 4.2)$& * & 536(70)\\
\hline
$x$ is controlled & * & 105(10) \\
\hline
$x$ is $i.i.d. P(4.2)$ &  *  & 1499(210) \\
\hline
\end{tabular}
\end{table}
The graph of the approximated $\overset{-1}{E\!D}$ for the case where $x$ has Pareto distribution is still decreasing when $k$ is around 100, hence the approximation is only feasible when $k>100$. We observe again that large values of $m_{4}$ and  $m_{22}$ of Pareto distribution $P(4.2)$ or $t$-distribution $t(4.2)$ lead to easier estimation than the case of normal distribution $N(0,1)$ when $\alpha=-1$.  However, just like the case when the error term has a normal distribution or t-distribution, the order of difficulty in the estimation is completely reversed with another $\alpha$. For example, see Figure \ref{fig:sn_x_-6} for the case when $\alpha=-6$.

We considered I.D.E. and R.S.S. w.r.t. Kullback-Leibler divergence under the condition $p=10$ for each distribution of $x$. Table \ref{ide_ssd_sn} shows the result. I.D.E. tells us that with any case of the moments of $x$, the regression model is easier to be estimated than the binomial model with the same $p-n$ ratio. If we divide R.S.S. with the number in the parenthesis, it is always less than 12. This means the $p-n$ ratio is always larger than that of the binomial model which has the same level of estimation difficulty as the regression model.  Especially when the distribution (moments) of $x$ is given as $t$-distribution or Pareto distribution, it makes the estimation easier.
\subsection{Comparison between different error distributions}
\label{comp_error_dists}
In this subsection, under a fixed distribution (moments) of $x$, we compared the approximated $\overset{\alpha}{E\!D}$'s for the three error distributions: the standard normal distribution (say $\overset{\alpha}{E\!D}_n$), the t-distribution with the d.f. of  3  (say $\overset{\alpha}{E\!D}_t$) and the skew-normal distribution with the shape parameter of 3 (say $\overset{\alpha}{E\!D}_s$). All comparisons are made  under the condition $p=10, n=12k$.

The order of the approximated $\overset{\alpha}{E\!D}$ among the three error distributions depends on $\alpha$. We pick up two values of $\alpha$, $\alpha=-1$ and $\alpha=-6$ as contrasting cases and summarized the results in Table \ref{table:comp_between_erros}. We notice that the order is completely reversed between $\alpha=-1$ and $\alpha=-6$. Under the fixed $\alpha$, $\overset{\alpha}{E\!D}_n$, $\overset{\alpha}{E\!D}_t$, $\overset{\alpha}{E\!D}_s$ keep the same order irrespective of the distribution of $x$.

We also present the graphs of  $\overset{\alpha}{E\!D}_n$, $\overset{\alpha}{E\!D}_t$, $\overset{\alpha}{E\!D}_s$ for each fixed distribution (moments) of $x$ with the reference to that of the normal coin toss model $B(k, 0.5)$. We notice that there is only little difference among the three error distributions and the normal coin toss model for the Kullback-Leibler divergence, especially when the distribution of $x$ is $N_{10}(0, I_{10})$ or controlled.

As for I.D.E. and R.S.S., we can make the comparison between different error term distributions if we look through Tables \ref{ide_ssd_normal}, \ref{ide_ssd_t} and \ref{ide_ssd_sn} with a fixed distribution of $x$. I.D.E. again indicates that the regression model can be more easily estimated than the coin toss model with any of the three error distributions.Though R.S.S. shows the sample size required  do not differ so much among the error term distributions, if pressed, t(3) requires a bit larger size of samples. If we divide R.S.S with the number in the parenthesis, we notice that $t(3)$ is always larger than the other distributions.
\begin{table}
\caption{Comparison between different error distributions}
\label{table:comp_between_erros}
\begin{center} 
\begin{tabular}{|l|c|c|}
\hline
  & $\alpha=-1$ & $\alpha=-6$ \\
\hline
$x\sim N_{10}(0, I_{10})$ & $\overset{-1}{E\!D}_t > \overset{-1}{E\!D}_n > \overset{-1}{E\!D}_s$  & 
 $\overset{-6}{E\!D}_s > \overset{-6}{E\!D}_n > \overset{-6}{E\!D}_t$\\
\hline   
$x\sim t_{10}(0, I_{10}, 4.2)$ & $\overset{-1}{E\!D}_t> \overset{-1}{E\!D}_n >\overset{-1}{E\!D}_s$ 
& $\overset{-6}{E\!D}_s >\overset{-6}{E\!D}_n > \overset{-6}{E\!D}_t$ \\
\hline
$x$ is controlled &  $\overset{-1}{E\!D}_t > \overset{-1}{E\!D}_n > \overset{-1}{E\!D}_s$ 
&  $\overset{-6}{E\!D}_s > \overset{-6}{E\!D}_n > \overset{-6}{E\!D}_t $\\
\hline 
$x$ is $i.i.d. P(4.2)$ &  $\overset{-1}{E\!D}_t> \overset{-1}{E\!D}_n> \overset{-1}{E\!D}_s$  
&  $\overset{-6}{E\!D}_s > \overset{-6}{E\!D}_n> \overset{-6}{E\!D}_t$ \\
\hline
\end{tabular}
\end{center}
\end{table}
\begin{figure}
\centering
\includegraphics[height=4.5cm,clip]{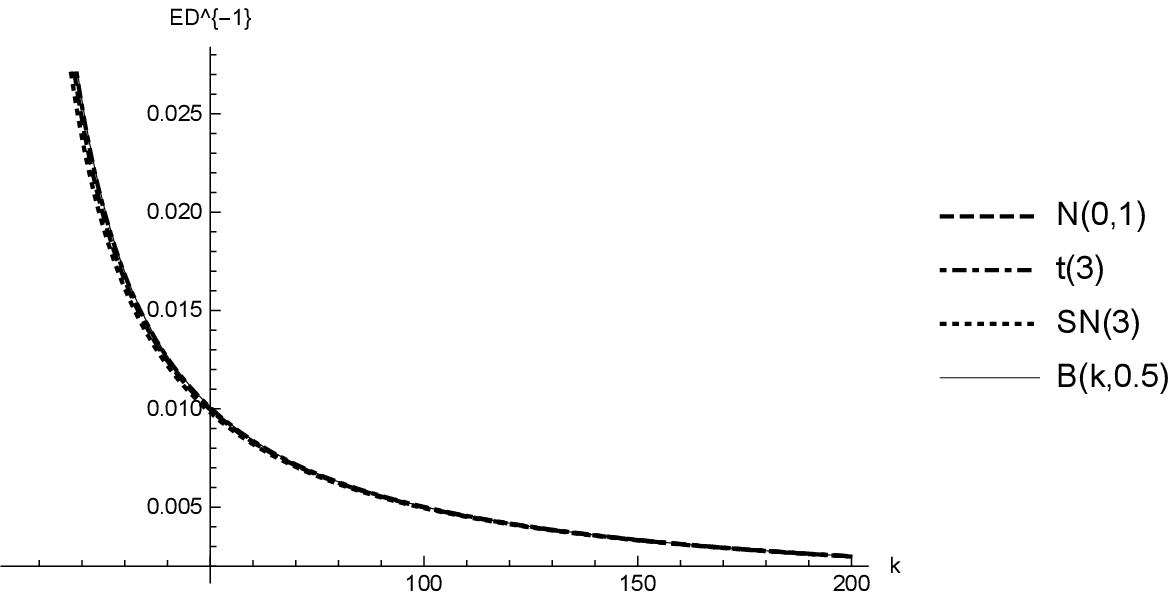}
\caption{$\overset{-1}{E\!D}$ when $x \sim N_{10}(0, I_{10})$}
\label{fig:x_N_p_-1}
\includegraphics[height=4.5cm,clip]{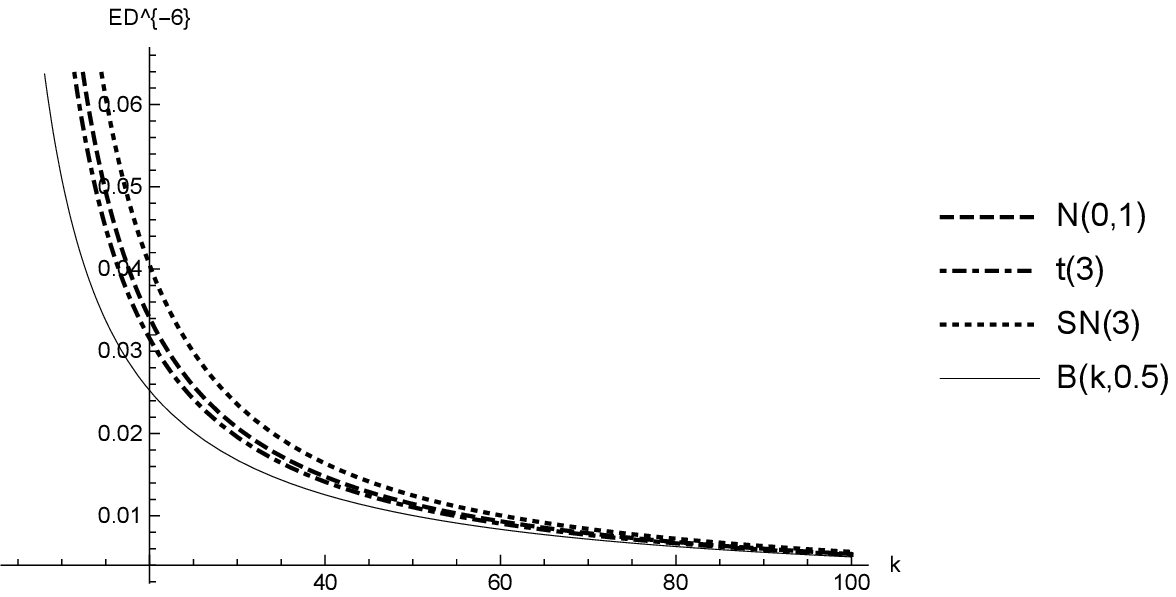}
\caption{$\overset{-6}{E\!D}$ when $x \sim N_{10}(0, I_{10})$}
\label{fig:x_N_p_-6}
\includegraphics[height=4.5cm,clip]{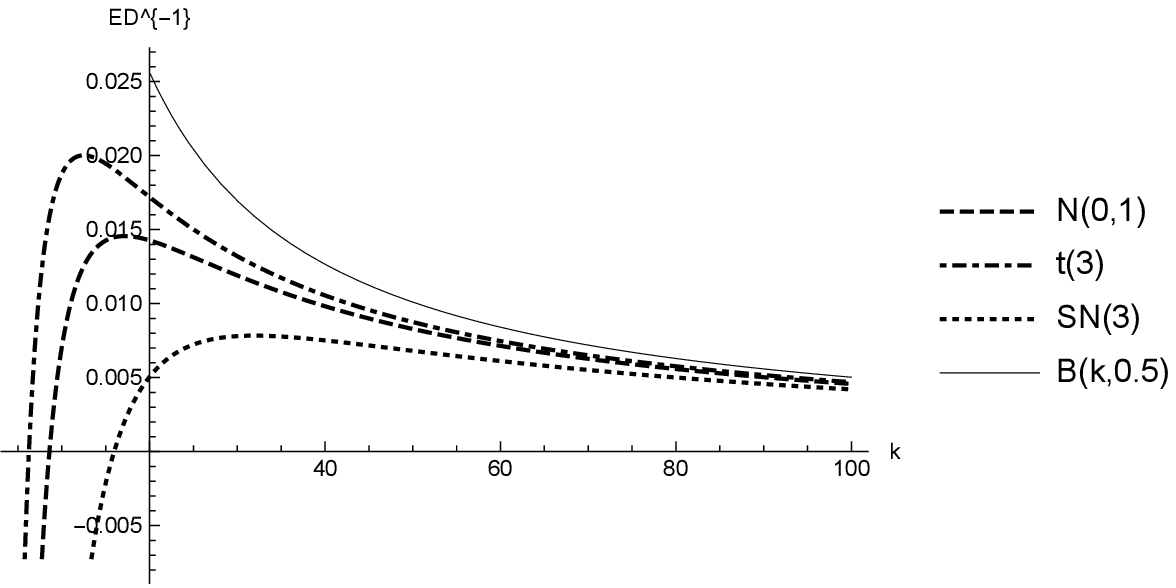}
\caption{$\overset{-1}{E\!D}$ when $x \sim t_{10}(0, I_{10}, 4.2)$}
\label{fig:x_t_p_-1}
\includegraphics[height=4.5cm,clip]{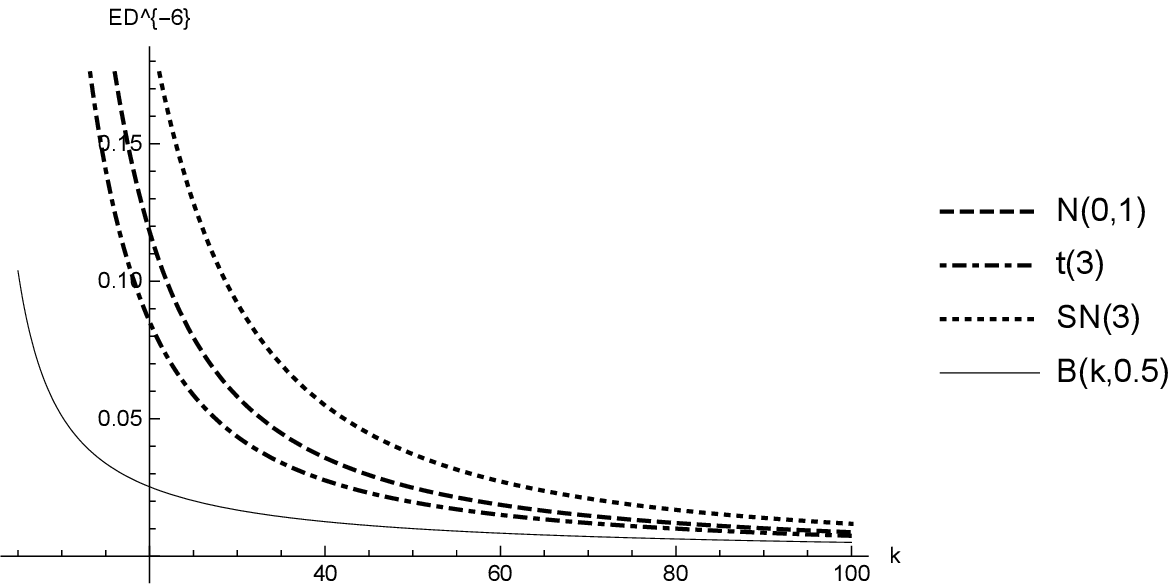}
\caption{$\overset{-6}{E\!D}$ when $x \sim t_{10}(0, I_{10}, 4.2)$}
\label{fig:x_t_p_-6}
\end{figure}
\begin{figure}
\centering
\includegraphics[height=4.5cm,clip]{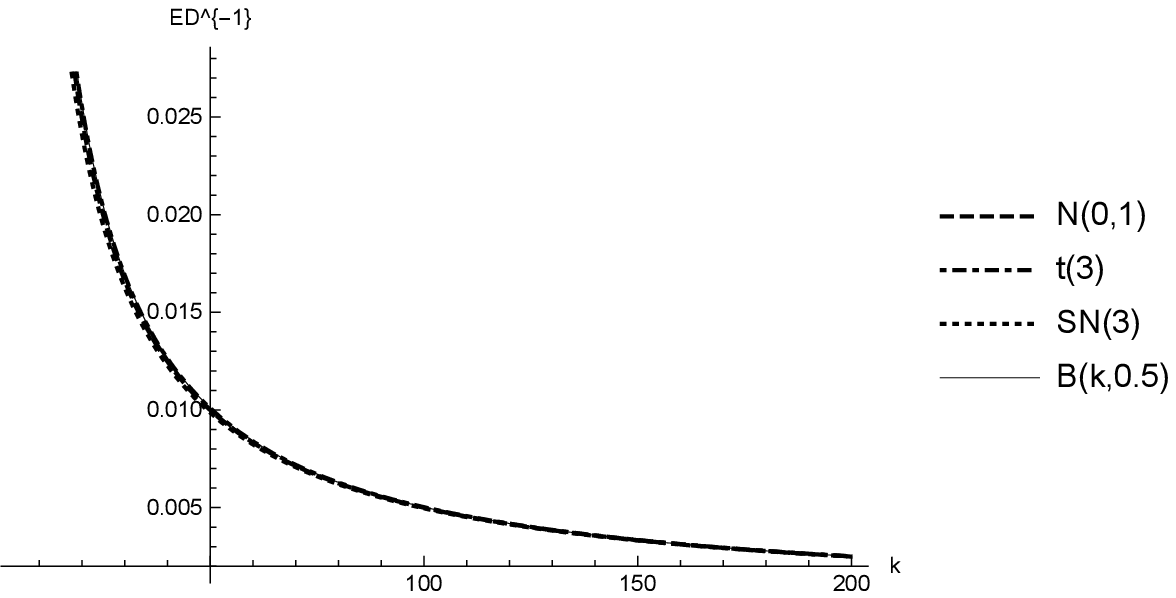}
\caption{$\overset{-1}{E\!D}$ when $x$ is controlled}
\label{fig:x_c_-1}
\includegraphics[height=4.5cm,clip]{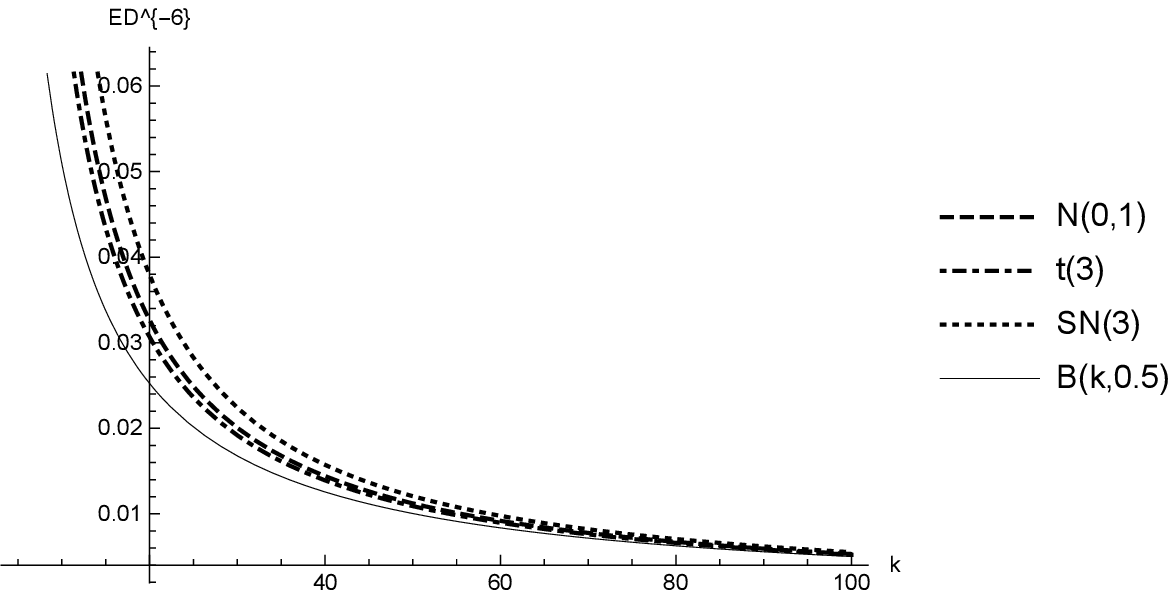}
\caption{$\overset{-6}{E\!D}$ when $x$ is controlled}
\label{fig:x_c_-6}
\includegraphics[height=4.5cm,clip]{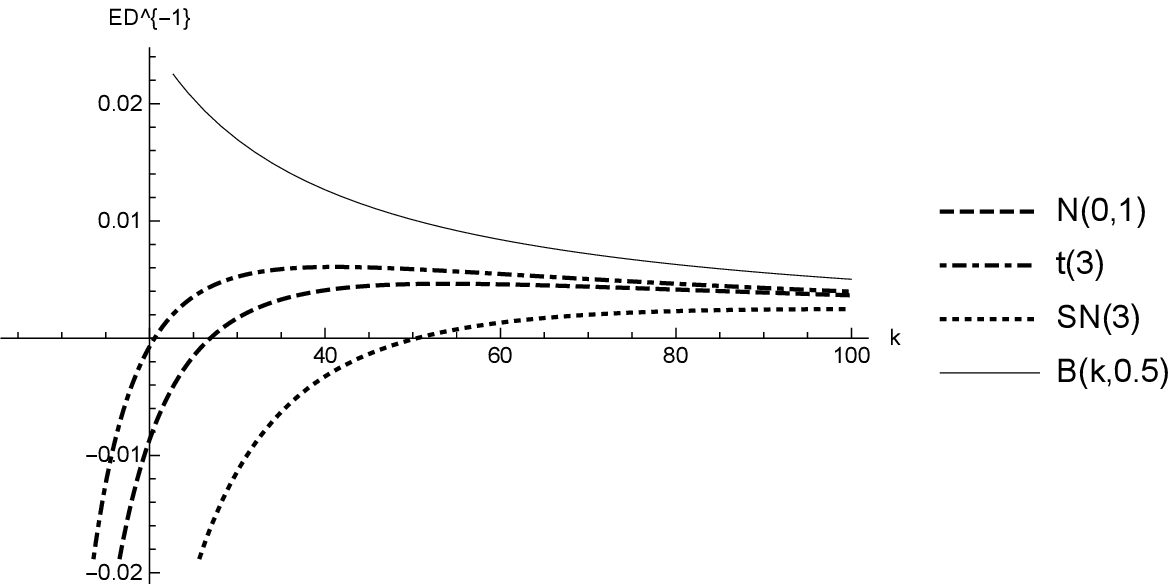}
\caption{$\overset{-1}{E\!D}$ when $x$ is i.i.d. as Pareto(4.2)}
\label{fig:x_P_-1}
\includegraphics[height=4.5cm,clip]{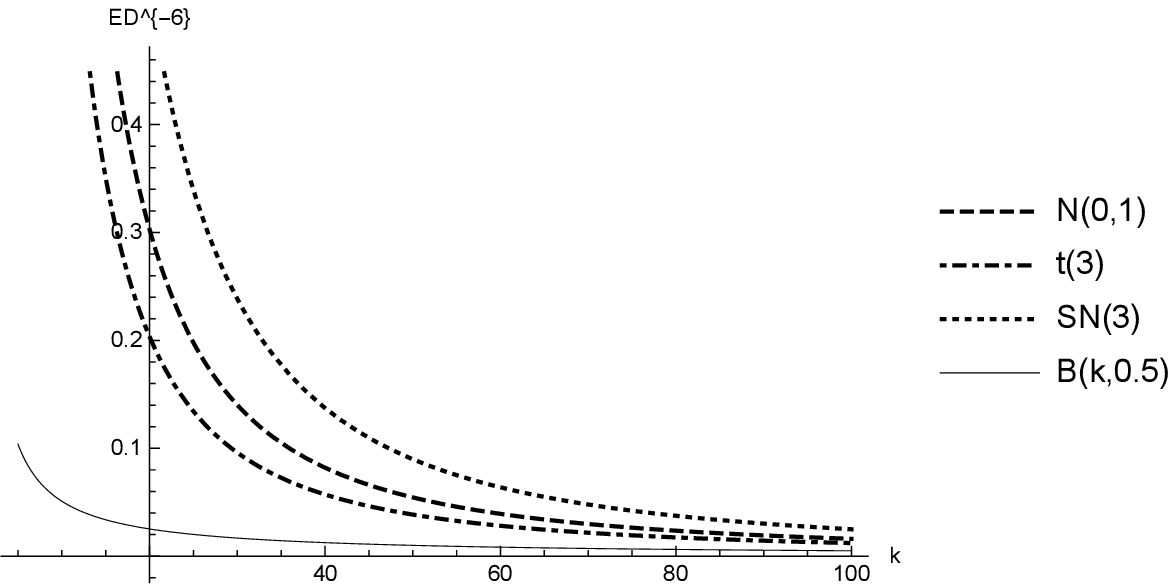}
\caption{$\overset{-6}{E\!D}$ when $x$ is i.i.d. as Pareto(4.2)}
\label{fig:x_P_-6}
\end{figure}
\section{Real Data --non-homogeneous explanatory variables--}
\label{real_data}
In this section we deal with two real datasets. As well as examples of non-homogeneous explanatory variables, these datasets serves as concrete cases to which the general results in the previous sections can be applied.  We calculate the sample moments of  the explanatory variables of those datasets and use them as examples of the following moments of $x$ (These datasets also include the dependent variables, but we do not use them here.)
\begin{equation}
\label{joint_moments_x_def}
m[i, j, k]= E[x_i x_j x_k] \qquad m[i, j, k, l]= E[x_i x_j x_k x_l] \qquad   1 \leq i, j, k, l  \leq p.
\end{equation}
First in order to standardize $x$ as in \eqref{stand_x}, we transform $x$ into its principal component scores. Then we calculate the  moments of the transformed $x$,
$$
n^{-1}\sum_{t=1}^n x_{ti} x_{tj} x_{tk} \qquad n^{-1}\sum_{t=1}^n x_{ti} x_{tj} x_{tk} x_{tl} \qquad   1 \leq i, j, k, l  \leq 11,
$$
and use them instead of \eqref{joint_moments_x_def} for the calculation of the aggregated sample moments 
\begin{align}
M_{2a}&\triangleq \sum_{i,j,k \in \mathcal{I}}m^2[i,j,k]\label{def_M2a}\\
M_{2b}&\triangleq \sum_{i,j,k \in \mathcal{I}}m[i,i,k] m[j,j,k] \label{def_M2b}\\
M_1&\triangleq \sum_{i,k \in \mathcal{I}}m[i,i,k,k].\label{def_M1}
\end{align}
Actually  $\overset{\alpha}{E\!D}$ is affected by the moments of $x$ only through these aggregated moments. (See the last part of Appendix A) 

\begin{figure}
\centering
\includegraphics[height=4.5cm,clip]{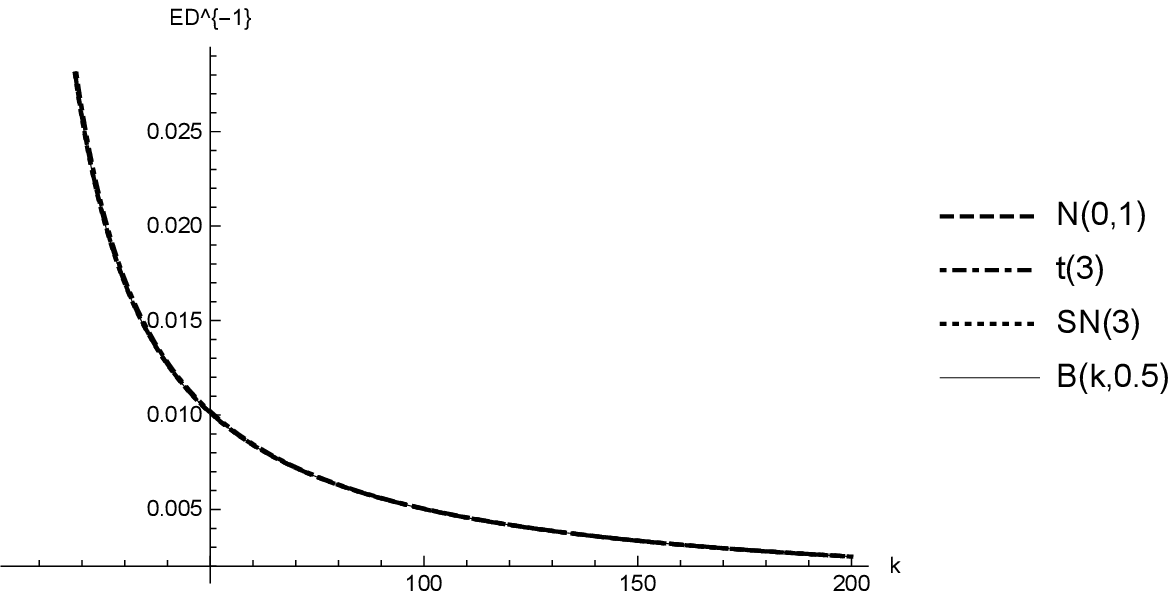}
\caption{$\overset{-1}{E\!D}$ for the wine data}
\label{fig:x_wine_-1_1}
\includegraphics[height=4.5cm,clip]{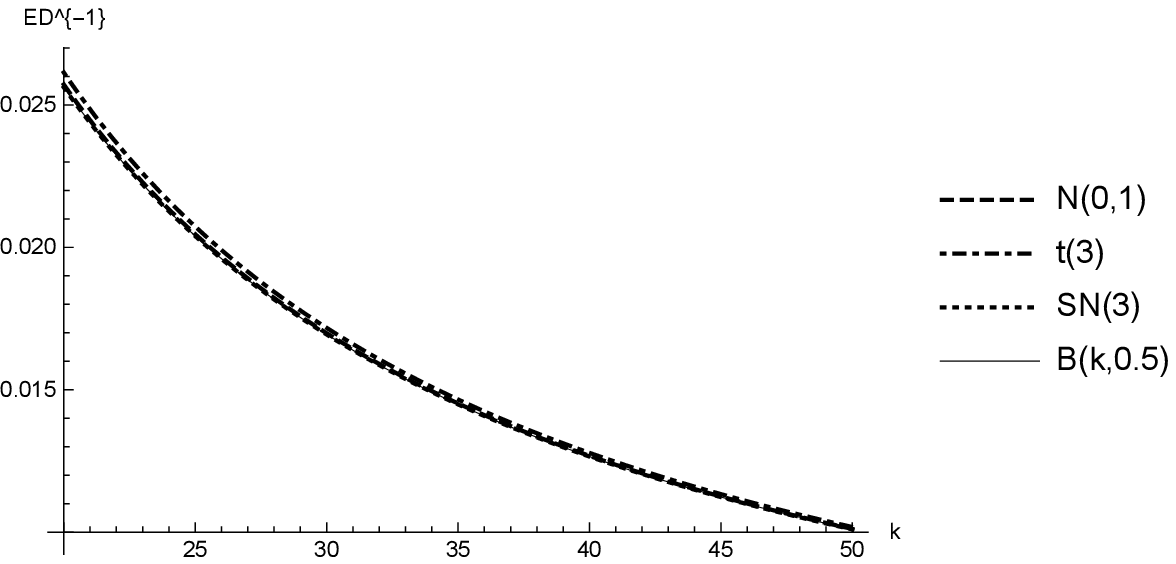}
\caption{$\overset{-1}{E\!D}$ for the wine data}
\label{fig:x_wine_-1_2}
\end{figure}

Since the results for those datasets are quite similar among different $\alpha$'s ($\alpha=-1, 0, 1,-6, 6$), we focus  ourselves on the case $\alpha=-1$.
\bigskip
\\
-- Example 1: Wine Quality --\\
This is the famous dataset on wine quality used in Cortez et.al. \cite{Cortez_et_al}. The data file is available at U.C.I. Machine Learning Repository  (https://archive.ics.uci.edu/ml/datasets\\/Wine+Quality). 
We used the white wine dataset. The dataset is as follows;
\\
$y \text{ (dependent variable)} = (y_t)_{1\leq t \leq n}$: the quality score of the wine form 0 to 10.
\\
$x \text{ (explanatory variables)} = (x_{ti})_{1\leq t \leq n, 1\leq i \leq 11}$: $n \times 11$ real value data on the quantity of the chemical substances in the wines . Each column is the data for the corresponding explanatory variable. $x_1$: ``fixed acidity'',  $x_2$: ``volatile acidity'', ... , $x_{11}$: ``alcohol''. 
\\
$n$ (sample size): 4898
\bigskip
\\
The values of \eqref{def_M2a} to \eqref{def_M1} for this dataset are as follows;
\begin{equation}
\label{aggre_moment_wine}
M_{2a}=0.000326899,\qquad M_{2b}=0.000230836, \qquad M_1=0.116967. 
\end{equation}
\bigskip
$M_{2a}$ or $M_{2b}$ is the summation over $11^3$ pieces of the squared 3-dimensional joint moments of  $x$. 
Since their averages $M_{2a}/11^3$ and $M_{2b}/11^3$ are quite small compared to the unit variance of $x_i$, this indicates $x$ are quite symmetric around the origin. $M_1/11^2$ is also much  smaller than 1, hence the distribution of $x$ has shorter tail  than the normal distribution.  $\overset{\alpha}{E\!D}$ is given as 
\begin{align}
\overset{\alpha}{E\!D}=
\begin{cases}
\displaystyle \frac{6.5}{n}+\frac{6.386 \alpha^2+48.804\alpha+91.026}{n^2}&\text{ if $\epsilon \sim N(0,1)$,}\\
\displaystyle \frac{6.5}{n}+\frac{1.927 \alpha^2+13.948\alpha+89.043}{n^2}&\text{ if $\epsilon \sim t(3)$,}\\
\displaystyle \frac{6.5}{n}+\frac{8.586 \alpha^2+71.639\alpha+104.856}{n^2}&\text{ if $\epsilon \sim SN(3)$.}
\end{cases}
\end{align}
Figure \ref{fig:x_wine_-1_1} ($k$ varies from 5 to 200) and \ref {fig:x_wine_-1_2} ($k$ varies from 20 to 50) show the graphs of $\overset{\scalebox{0.6}{\!\!$-1$}}{E\!D}$ for three error distributions under this moments of $x$ and $n=13k$.  We also put the graph of $\overset{\scalebox{0.6}{\!\!$-1$}}{E\!D}$ of $B(0.5,k)$ as a reference. We see that $\overset{\scalebox{0.6}{\!\!$-1$}}{E\!D}$'s of the four cases are quite close to each other. There is almost no difference among the error distributions. Besides, the estimation difficulty of the regression model is similar to that of the normal coin toss with the same $p-n$ ratio irrespective of the error distributions.
\begin{table}
\centering
\caption{I.D.E. \& R.S.S. for the wine data}
\label{ide_rss_wine}
\begin{tabular}{|c|c|c|}
\hline 
      & I.D.E. &  R.S.S.\\
\hline
$N(0, 1)$ & 0.66 & 130(10)\\
\hline
$t(3)$ &  0.81 & 135(10)\\
\hline
$SN(3)$  & * & 130(10)\\
\hline
\end{tabular}
\end{table}

I.D.E. and R.S.S. is stated in Table \ref{ide_rss_wine}. I.D.E. shows that when the error distribution is $SN(3)$, the estimation is the easiest and that when $\epsilon \sim t(3)$, the estimation gets slightly more difficult. As for R.S.S., we notice that there is very little difference among the three error distributions and that the estimation is relatively easy. Around 130 samples guarantee as easy estimation as the 10-times normal coin toss problem. 

We can evaluate the actual sample size 4898 of this dataset by answering the following question; how large sample size $n$ for the normal coin toss model $B(0.5, n)$ is required in order to attain the same level of easiness in the estimation as the regression model with the moments of $x$ as in \eqref{aggre_moment_wine} and the sample size 4898 ? For example, if the error distribution is $N(0,1)$, then the answer is given as the solution of  the equation
\begin{equation}
\label{sample_size_eval_wine}
\frac{1}{2n}+\frac{1}{8n^{2}}((\alpha')^2-5\alpha'+6)=
\frac{6.5}{4898}+\frac{6.386 \alpha^2+48.804\alpha+91.026}{4898^2},
\end{equation}
where the left-hand side is \eqref{expan_binom} with $M=4$.
\\
The rounded solution when $\alpha=-1$ equals 376 or 377 for the three error distributions, which means the sample size 4898 for the regression is equivalent to the 376 (377) times normal coin toss in view of the estimation difficulty. We see that the estimation is fairly easy with this sample size.
\bigskip
\\
-- Example 2: Communities and Crime --\\
This data combines socio-economic data for each community within USA from the 1990 US Census, law enforcement data from the 1990 US LEMAS survey, and crime data from the 1995 FBI U.C.R.. You can download the data file from U.C.I. Machine Learning Repository (https://archive.ics.uci.edu/ml/datasets/Communities+and+Crime).

\begin{figure}
\centering
\includegraphics[height=4.5cm,clip]{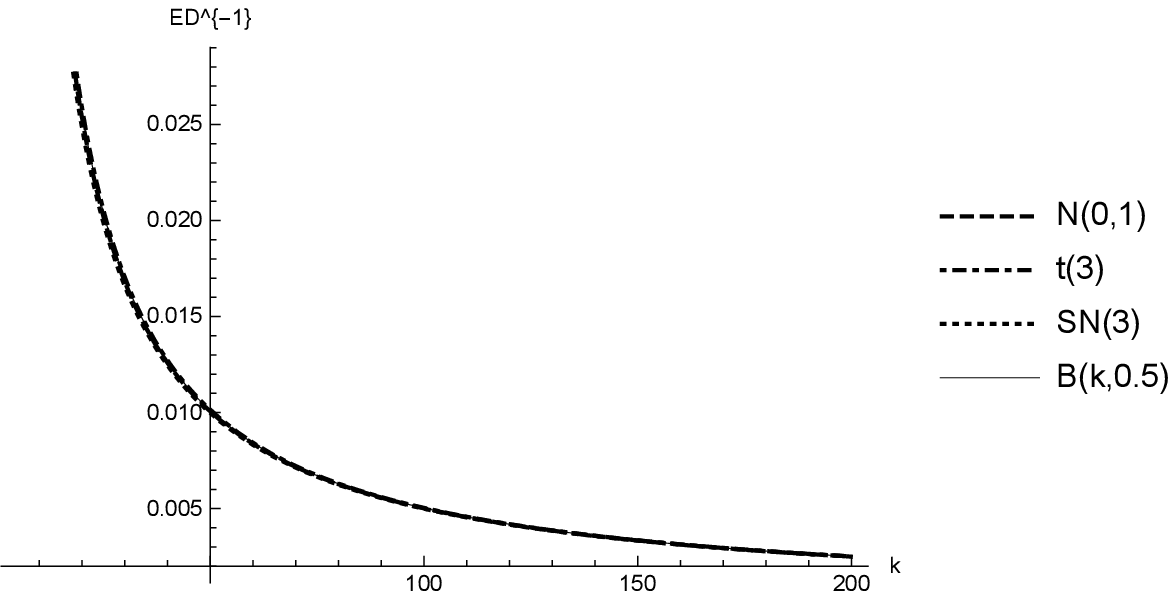}
\caption{$\overset{-1}{E\!D}$ for the crime data}
\label{fig:x_crime_-1_1}
\includegraphics[height=4.5cm,clip]{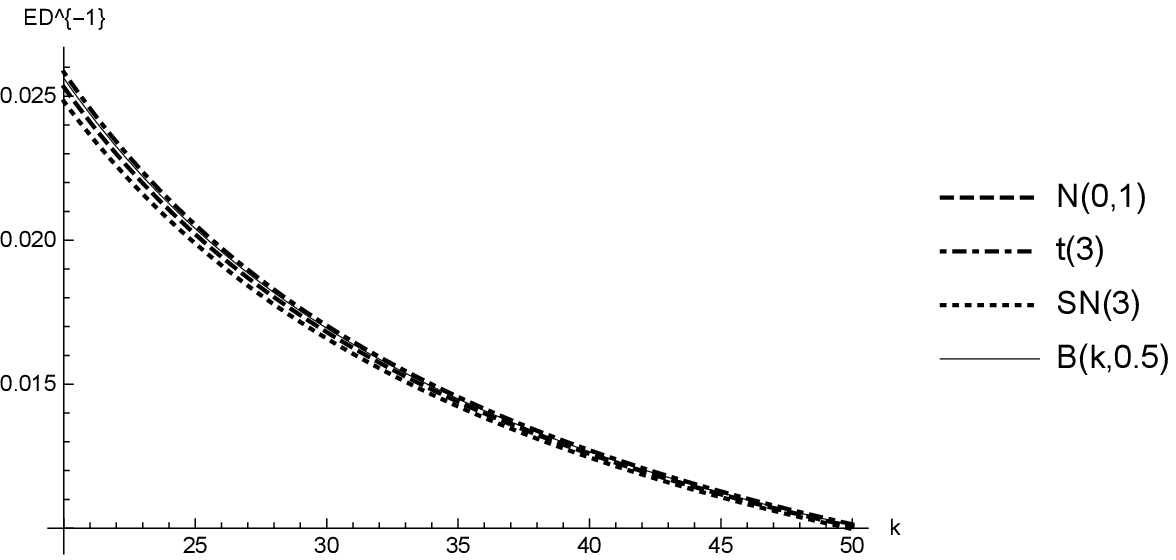}
\caption{$\overset{-1}{E\!D}$ for the crime data}
\label{fig:x_crime_-1_2}
\end{figure}

The original data contains 124 explanatory variables from ``population'' to ``PolicBudgPerPop''. We excluded the explanatory variables that contains missing data (denoted by "?" in the original dataset) . Besides we excluded the variable "numbUrban","PctRecImmig8" and "OwnOccMedVal" because the following correlations exceed 0.99: Corr(``population'', ``numUrban''), Corr(``PctRecImmig5'',''PctRecImmig8''), Corr(``PctRecImmig8'',''PctRecImmig10''),  Corr(``OwnOccLowQuart'',''OwnOccMedVal''). After this process, the dataset is as follows;
\\
$y \text{ (dependent variable)} = (y_t)_{1\leq t \leq n}$:  The candidates of $y$ are 18 attributes from ``murders'' (the number of the murders committed in the community) to  ``nonViolPerPop''(the number per capita of non-violent crimes committed in the community). They are the numbers of the committed crimes categorized in various ways.
\\
$x \text{ (explanatory variables)} = (x_{ti})_{1\leq t \leq n, 1\leq i \leq 99}$: $n \times 99$ real value data on the socio-economic character of the community.  $x_1$: ``population'', $x_2$: ``household''(mean people per household) ,..., $x_{99}$: ``LemasPctOfficDrugUn''(the percent of officers assigned to drug units ).
\\
$n$ (sample size): 2215
\\
We used principle component sores as the standardized $x$. The aggregated sample moments are given by
$$
M_{2a}=1708.97,\qquad M_{2b}=1749.28, \qquad M_1=2604.5.
$$
$M_{2a}/99^3$, $M_{2b}/99^3$ and $M_1/99^2$ are much smaller than unit. Like the wine data, the distribution of $x$ is symmetric and short-tailed.  Using these values we calculated $\overset{\alpha}{E\!D}$, which is given by
\begin{align}
\overset{\alpha}{E\!D}=
\begin{cases}
\displaystyle \frac{50.5}{n}+\frac{294.547 \alpha^2+1949.71 \alpha+2953.58}{n^2}&\text{ if $\epsilon \sim N(0,1)$,}\\
\displaystyle \frac{50.5}{n}+\frac{137.758 \alpha^2+111.409 \alpha+3376.96 }{n^2}&\text{ if $\epsilon \sim t(3)$,}\\
\displaystyle \frac{50.5}{n}+\frac{526.088 \alpha^2+3232.22 \alpha+1976.26 }{n^2}&\text{ if $\epsilon \sim SN(3)$.}
\end{cases}
\end{align}
Figure \ref{fig:x_crime_-1_1} ($k$ varies from 5 to 200) and \ref {fig:x_crime_-1_2} ($k$ varies from 20 to 50) show shows the graphs of $\overset{\scalebox{0.6}{\!\!$-1$}}{E\!D}$ for the three error distributions under these moments of $x$ and $n=101k$.  We also put the graph of $\overset{\scalebox{0.6}{\!\!$-1$}}{E\!D}$ of $B(0.5,k)$ as a reference. The comment for Example 1 holds for this data. We see that $\overset{\scalebox{0.6}{\!\!$-1$}}{E\!D}$'s for the three error distributions are almost same. Compared to the normal coin toss  with the same $p-n$ ratio, the regression model is on the same level for the estimation difficulty.

You can see I.D.E. and R.S.S. in Table \ref{ide_rss_crime}. We notice that it is slightly harder to estimate the parameters when $\epsilon \sim t(3)$, but, generally speaking, for the regression model with these moments of $x$, estimating the parameters is not a hard task if we have around 1000 samples. 
\begin{table}
\centering
\caption{I.D.E. \& R.S.S. for the crime data}
\label{ide_rss_crime}
\begin{tabular}{|c|c|c|}
\hline 
      & I.D.E. &  R.S.S.\\
\hline
$N(0, 1)$ & * & 987(10)\\
\hline
$t(3)$ &  0.72 & 1025(10)\\
\hline
$SN(3)$  & * & 947(10)\\
\hline
\end{tabular}
\end{table}
We evaluate the sample size 2215 in a similar way to \eqref{sample_size_eval_wine}. If the error distribution is $N(0,10)$, then solving the equation
\begin{equation}
\label{sample_size_eval_crime}
\frac{1}{2n}+\frac{1}{8n^{2}}((\alpha')^2-5\alpha'+6)=
\frac{50.5}{2215}+\frac{294.547 \alpha^2+1949.71 \alpha+2953.58}{2215^2}
\end{equation}
gives us an evaluation of the actual sample size. When $\alpha=-1$, the rounded solution is 22 or 23 for the three error distributions. Though this number is much smaller than 376(377) in Example 1, the estimation is still not a hard task since 22-times normal coin toss gives us plenty of information.
\section{Summary and Discussion}
\label{discussion}
\begin{itemize}
\item $\overset{\alpha}{E\!D}$ is constant for the parameter $\beta, \sigma$.
\item The main term ($n^{-1}$ term) of the asymptotic expansion of $\overset{\alpha}{E\!D}$ is $(p+2)/n$, that is, $p-n$ ratio. 
\item For the second term ($n^{-2}$ term) of the expansion, we observe the following points.
\begin{enumerate}
\item The maximum dimension of $p$ depends on the error term distribution. It can be more than two as in the case $\epsilon \sim SN(3)$, where it is not enough to increase the sample size proportionally to $p$ for reliable estimation (so called "the curse of dimension").
\item The joint moments that appear in the term is maximally of the forth order. What moments appear is different among the error term distributions. If it is a quadratic distribution (e.g. $N(0,1)$, $t(\nu)$ ), then the moments $m_4$ and $m_{22}$ only appear.
\item The effect of  $m_4$ and $m_{22}$ depends on $\alpha.$ When $\alpha=-1, 0, 1, 3$, the larger $m_4$ and $m_{22}$ decreases the difficulty of the estimation. In a geometrical view, there is no preference among $\alpha$'s. Each $\alpha$ gives its own geometrical structure to Riemannian manifold formed by the parametric distribution model (see e.g.  Amari  and Nagaoka \cite{Amari&Nagaoka}). However there might be values for $\alpha$ that is ``natural'' in a statistical sense or ``appropriate'' for a purpose of the statistical analysis.
\item The effect of the error term distributions also depends on $\alpha$. For example,  the order of the estimation difficulty among the three error distributions is quite different between $\alpha=-1$ and $\alpha=-6$.
\item The difference between the three error term distributions we investigated is relatively small if we use Kullback-Leibler divergence.This might be  due to the assumption that we know the error term distributions, hence are able to use m.l.e.  In most applications, the actual error term distribution is unknown, and m.l.e. is not applicable. It is of much interest what would happen to the risk of the predictive distribution, if we use another estimator such as the least squares estimator. 
\end{enumerate}
\item 
We proposed  measuring the (asymptotic) difficulty of estimation by the approximated $\overset{\alpha}{E\!D}$ and tried to give a suggestion on the sample size. It is a method comparing the approximated $\overset{\alpha}{E\!D}$ of the regression model to that of a binomial model $B(n,m)$.  I.D.E. tells the difficulty of estimation by the value of $m$ of $B(k, m)$, which has the same $p-n$ ratio as the regression model \eqref{regression_model} of the sample size $(p+2)k$. R.S.S. gives the sample size $n$ for the regression model which leads to the same difficulty of estimation as $B(10, 0.5)$ (If it is needed, a more large value than 10 will be used for the binomial model). 
\begin{enumerate}
\item Though there exist small difference between the error term distributions and the moments of $x$, in most cases we investigated, the regression model is easier to be estimated than the normal coin toss $B(k, 0.5)$ under the same $p-n$ ratio $1/k$. 
\item The sample size $n=10(p+2)$ guarantees the good performance of the estimation at the same level as the 10-times normal coin toss irrespective of the error term distributions and the moments of $x$ which we investigated in this paper.
\end{enumerate}
\end{itemize}
\appendix
\section{Calculation of \eqref{g_ij} -- \eqref{eta_ssss} and \eqref{def_L}}
\label{detailed_cal}
First we give the detailed calculation process of the results from \eqref{g_ij} to \eqref{eta_ssss}.
Since 
$$
\log f(y,x \,|\, \beta, \sigma)=-\log \sigma+\log f(y^*)+\log h(x),\qquad y^*\triangleq \frac{y-\beta' \tilde{x}}{\sigma},
$$
and 
$$
\frac{\partial y^*}{\partial \sigma}=-\sigma^{-1} y^*,\qquad \frac{\partial y^*}{\partial \beta_i}=-\sigma^{-1}x_i\text{ for $0 \leq i \leq p$},
$$
we have the following results;
\begin{align}
l_i&=-\sigma^{-1}x_i \log'f(y^*)\text{ for $0\leq i \leq p$},\label{l_i}\\
l_\sigma&=-\sigma^{-1} (1+\log'f(y^*)y^*),\label{l_s}\\
l_{ij} &=\sigma^{-2}x_i x_j \log''f(y^*)\text{ for $0 \leq i, j \leq p$},\label{l_ij}\\
l_{i\sigma}&=\sigma^{-2} x_i (y^* \log''f(y^*)+\log' f(y^*))\text{ for $0\leq i \leq p$},\label{l_is}\\
l_{\sigma\sigma}&=\sigma^{-2}(1+\log'f(y^*)y^*+\log''f(y^*)(y^*)^2+\log'f(y^*)y^*)\nonumber\\
&=\sigma^{-2}(1+2\log'f(y^*)y^*+\log''f(y^*)(y^*)^2),\label{l_ss}\\
l_{ijk}&=-\sigma^{-3}x_i x_j x_k \log'''f(y^*)\text{ for $0 \leq i, j, k\leq p$},\label{l_ijk}\\
l_{ij\sigma}&=-\sigma^{-3}x_i x_j(2\log''f(y^*)+\log'''f(y^*)y^*)\text{ for $0 \leq i, j \leq p$},\label{l_ijs}\\
l_{i\sigma\sigma}&=-\sigma^{-3}x_i(2y^*\log''(y^*)+2\log'f(y^*)+y^*\log''f(y^*)\nonumber\\
&\qquad+(y^*)^2\log'''f(y^*)+y^*\log''f(y^*))\nonumber\\
&=-\sigma^{-3}x_i(4y^*\log''f(y^*)+2\log'f(y^*)+(y^*)^2\log'''f(y^*))\text{ for $0 \leq i \leq p$},\label{l_iss}\\
l_{\sigma\sigma\sigma}&=-\sigma^{-3}(2(y^*)^2\log''f(y^*)+2y^*\log'f(y^*)+\log'''f(y^*)(y^*)^3\nonumber\\
&\qquad+2\log''f(y^*)(y^*)^2+2+4\log'f(y^*)y^*+2(y^*)^2\log''f(y^*))\nonumber\\
&=-\sigma^{-3}(2+6(y^*)^2\log''f(y^*)+6y^*\log'f(y^*)+(y^*)^3\log'''f(y^*)).
\label{l_sss}
\end{align}

Now we state the proof of each result from \eqref{g_ij} to \eqref{eta_ssss}.
The next formula will be  often used in those proofs. 
Suppose that an integrable function $h(y)$ on $R$  allows the following exchangeability between the integral and the differentiation, 
$$
\int_R \frac{d }{d\epsilon}h(y+\epsilon)\biggl |_{\epsilon=0} dy =\frac{d }{d\epsilon}\int_R h(y+\epsilon) dy\biggl |_{\epsilon=0},
$$
then 
\begin{align}
\int_R h'(y) dy &=\int_R \frac{d }{d\epsilon}h(y+\epsilon)\biggl|_{\epsilon=0} dy \nonumber\\
&=\frac{d }{d\epsilon} \int_R h(y+\epsilon) dy \biggl|_{\epsilon=0}\nonumber\\
&=\frac{d }{d\epsilon} \int_R h(y) dy \biggl|_{\epsilon=0}=0. 
\end{align}
In the following proofs, all the functions derived from $f(y)$ are supposed to satisfy this exchangeability. We also use the notation 
$$
E[h(y)]\triangleq \int_{-\infty}^\infty h(y) f(y) dy
$$
-- Proof of \eqref{g_ij} --
\begin{align*}
g_{ij}\\
&=L_{ij}\\
&=\sigma^{-3}\int_{R^p} x_i x_j h(x) \int_R \Bigl(\frac{f'(y^*)}{f(y^*)}\Bigr)^2
f(y^*) dy dx \\
&=\sigma^{-2}\int_{R^p} x_i x_j h(x) \int_R \Bigl(\frac{f'(y)}{f(y)}\Bigr)^2 f(y) dy dx\\
&=\sigma^{-2}\int_{R^p} x_i x_j h(x) dx \:E\bigl[\bigl(\log'f(y)\bigr)^2\bigr]\\
&=\sigma^{-2}m[i,j]\eta[0,0,2,0]
\end{align*}
From the fact $x_0 \equiv 1$ and \eqref{stand_x}, $m[i,j]=\delta_{ij}$.
The following equation also holds;
\begin{align*}
\eta[0,0,2,0]+\eta[0,1,0,0]&=\int_R \Bigl(\bigl(\log' f(y)\bigr)^2+\log''f(y) \Bigr) f(y) dy \\
&=\int_ R f''(y) dy\\
&=0 .
\end{align*}
-- Proof of \eqref{g_is} --
\begin{align*}
&g_{i\sigma}\\
&=L_{i\sigma}\\
&=\sigma^{-3}\int_{R^p} x_i  h(x) \int_R \bigl(y^*\bigl(\log'f(y^*)\bigr)^2+\log'f(y^*)\bigr)
f(y^*) dy dx \\
&=\sigma^{-2}\int_{R^p} x_i h(x) dx \int_R \bigl(y\bigl(\log'f(y)\bigr)^2+\log'f(y)\bigr)
f(y) dy \\
&=\sigma^{-2}\int_{R^p} x_i h(x) dx \int_R y\bigl(\log'f(y)\bigr)^2
f(y) dy \\
&=\sigma^{-2}m[i]\eta[0,0,2,1]\\
&=
\begin{cases}
\sigma^{-2}\eta[0,0,2,1] &\text{ if $i=0$,}\\
0 &\text{ if $1\leq i \leq p$.}
\end{cases}
\end{align*}
The fourth equation is due to the following relation;
\begin{align*}
&\int_R \log'f(y) f(y) dy\\
&=\int_R f'(y) dy \\
&=0.
\end{align*}
We also have the equation,
\begin{align*}
\eta[0,0,2,1]+\eta[0,1,0,1]&=\int_R \Bigl(\bigl(\log' f(y)\bigr)^2+\log''f(y) \Bigr) y f(y) dy \\
&=\int_ R yf''(y) dy\\
&=\int_R yf''(y)+f'(y) dy\\
&=\int_R (y f'(y))'dy\\
&=0.
\end{align*}
\\
-- Proof of \eqref{g_ss} --
\begin{align*}
&g_{\sigma\sigma}\\
&=L_{\sigma\sigma}\\
&=\sigma^{-3} \int_R \bigl(1+\log'f(y^*)y^*\bigr)^2 f(y^*) dy\\
&=\sigma^{-2} \int_R \bigl(1+\log'f(y)y\bigr)^2 f(y) dy \\
&=\sigma^{-2} \int_R \bigl(1+2\log'f(y)y+y^2\bigl(\log'f(y)\bigr)^2\bigr)f(y) dy \\
&=\sigma^{-2}(1+2\eta[0,0,1,1]+\eta[0,0,2,2]).
\end{align*}
The following equations hold;
\begin{align*}
&\eta[0,1,0,2]+2\eta[0,0,1,1]+\eta[0,0,2,2]\\
&=\int_R \bigl(y^2 \log''f(y) +2 y\log'f(y)+\bigl(\log' f(y)\bigr)^2y^2\bigr) f(y) dy
&=\int_R \big(\log'f(y) y^2 f(y)\bigr)' dy
&=0,\\
&1+\eta[0,0,1,1]\\
&=\int_R \bigl( y f(y) \bigr)' dy\\
&=0.
\end{align*}
Therefore 
$$
1+2\eta[0,0,1,1]+\eta[0,0,2,2]=-(1+2\eta[0,0,1,1]+\eta[0,1,0,2]).
$$
\\
-- Proof of \eqref{g^ij}, \eqref{g^0i}, \eqref{g^00}, \eqref{g^0s} and \eqref{g^ss}--\\
From \eqref{g_ij}, \eqref{g_is} and \eqref{g_ss}, we notice that
\begin{align*}
\begin{pmatrix}
g^{00} & g^{0\sigma} \\
g^{0\sigma} & g^{\sigma\sigma}
\end{pmatrix}
&=
\begin{pmatrix}
g_{00} & g_{0\sigma} \\
g_{0\sigma} & g_{\sigma\sigma}
\end{pmatrix}^{-1}\\
&=\sigma^2
\begin{pmatrix}
\eta[0,0,2,0] & -\eta[0,1,0,1] \\
-\eta[0,1,0,1] & 1+2\eta[0,0,1,1]+\eta[0,0,2,2]
\end{pmatrix}^{-1}\\
&=\Delta^{-1}\sigma^2
\begin{pmatrix}
 1+2\eta[0,0,1,1]+\eta[0,0,2,2] & \eta[0,1,0,1] \\
\eta[0,1,0,1] & \eta[0,0,2,0]
\end{pmatrix}
\end{align*}
\eqref{g^ij} and \eqref{g^0i} are obvious.\\
\\
-- Proof of \eqref{eta_(ij)k}, \eqref{eta_(is)k},\eqref{eta_(ij)s}, \eqref{eta_(is)s}, \eqref{eta_(ss)i} and \eqref{eta_(ss)s}--\\
The equations from \eqref{eta_(ij)k} to \eqref{eta_(ss)i} are straightforward  from the results \eqref{l_i} to \eqref{l_ss}. \eqref{eta_(ss)s} is gained as follows;
\begin{align*}
&L_{(\sigma\sigma)\sigma}\\
&=-\sigma^{-4}\int_{R^p} h(x) \int_R \bigl(1+\log''f(y^*) (y^*)^2+ 2\log'f(y^*)y^*\bigr)\\
&\hspace{40mm}\times\bigl(1+\log'f(y^*)y^*\bigr)f(y^*) dy dx\\
&=-\sigma^{-3} \int_R \bigl(1+\log''f(y) y^2+ 2\log'f(y)y\bigr)\bigl(1+\log'f(y)y\bigr)f(y) dy\\
&=-\sigma^{-3} \bigl(1+\eta[0,1,0,2]+2\eta[0,0,1,1]+\eta[0,0,1,1]+\eta[0,1,1,3]+2\eta[0,0,2,2]\bigr)\\
&=-\sigma^{-3}\bigl(1+3\eta[0,0,1,1]+\eta[0,1,0,2]+2\eta[0,0,2,2]+\eta[0,1,1,3]\bigr)
\end{align*}
-- Proof of \eqref{eta_ijk}, \eqref{eta_ijs},\eqref{eta_iss} and \eqref{eta_sss} --\\
These equations are almost obvious from \eqref{l_i} and \eqref{l_s}.
\\
\\
-- Proof of \eqref{eta_(ij)(kl)}, \eqref{eta_(is)(kl)}, \eqref{eta_(is)(js)}, \eqref{eta_(ij)(ss)}, \eqref{eta_(is)(ss)} and \eqref{eta_(ss)(ss)} --\\
\eqref{eta_(ij)(kl)}, \eqref{eta_(is)(kl)}, \eqref{eta_(is)(js)}, \eqref{eta_(ij)(ss)} are obvious from \eqref{l_ij} to \eqref{l_ss}.
\eqref{eta_(is)(ss)} and \eqref{eta_(ss)(ss)} are gained as follows;
\begin{align*}
&L_{(i\sigma)(\sigma\sigma)}\\
&=\sigma^{-4}\int_{R^p} x_i dx \int_R \bigl(y\log''f(y)+\log'f(y)\bigr)\\
&\hspace{40mm}\times\bigl(1+2\log'f(y)y+\log''f(y)y^2\bigr) dy\\
&=\sigma^{-4}m[i] \bigl(\eta[0,1,0,1]+\eta[0,2,0,3]+2\eta[0,1,1,2]+\eta[0,1,1,2]+2\eta[0,0,2,1]\bigr)\\
&=\sigma^{-4}m[i] \bigl(\eta[0,1,0,1]+\eta[0,2,0,3]+3\eta[0,1,1,2]+2\eta[0,0,2,1]\bigr),\\
&L_{(\sigma\sigma)(\sigma\sigma)}\\
&=\sigma^{-4}\int_R \bigl(1+2\log'f(y) y+\log''f(y) y^2\bigr)^2 f(y) dy\\
&=\sigma^{-4}\bigl(1+\eta[0,2,0,4]+4\eta[0,0,2,2]+2\eta[0,1,0,2]+4\eta[0,0,1,1]+4\eta[0,1,1,3]\bigr).
\end{align*}
\\
\\
-- Proof of \eqref{eta_(ijk)l} to \eqref{eta_(sss)s}--\\
We only describe the proof of \eqref{eta_(iss)s}. Other equations are instantly gained from \eqref{l_i}, \eqref{l_s},  \eqref{l_ijk} -- \eqref{l_sss}. 
\begin{align*}
&L_{(i\sigma\sigma)\sigma}\\
&=\sigma^{-4}\int_{R^p} x_i h(x) dx \int_R \bigl(1+\log'f(y) y)\\
&\hspace{40mm}\times \bigl(4y\log''f(y)+2\log'f(y)+y^2\log'''f(y)\bigr)dy\\
&=m[i]\bigl( 4\eta[0,1,0,1]+\eta[1,0,0,2]+4\eta[0,1,1,2]+2\eta[0,0,2,1]+\eta[1,0,1,3]\bigr).
\end{align*}
\\
\\
-- Proof of \eqref{eta_(ij)kl} to \eqref{eta_(ss)ss}--\\
All the equations are instantly gained from \eqref{l_i} -- \eqref{l_ss}.
\\
\\
-- Proof of \eqref{eta_ijkl} to \eqref{eta_ssss}--\\
All the equations are instantly gained from \eqref{l_i} and \eqref{l_s}.
\\
\\
-- Calculation of $\eqref{def_L}$ --\\
We state here the calculation process of $\eqref{def_L}$, which is actually used in Mathematica programming in Appendix \ref{Prog_Mathematica}.
Let 
$$
\tilde{g}^{ij}=\sigma^{-2}g^{ij}, \quad i, j \in \{0, 1, \ldots, p, \sigma\},
$$
and
\begin{align}
M_{2a}&\triangleq \sum_{i,j,k \in \mathcal{I}}m^2[i,j,k] \\
M_{2b}&\triangleq \sum_{i,j,k \in \mathcal{I}}m[i,i,k] m[j,j,k]\\
M_1&\triangleq \sum_{i,k \in \mathcal{I}}m[i,i,k,k]
\end{align}
Note that if $x$ is homogeneous, $M_{2a}, M_{2b}, M_1$ are respectively given by
\begin{align}
M_{2a}&=\sum_{ i \in \mathcal{I}}m^2[i,i,i]+3\sum_{i, j \in \mathcal{I}, i \ne j}m^2[i,i,j]+\sum_{i, j, k \in \mathcal{I}, i \ne j, i\ne k, j\ne k}m^2[i,j,k]\nonumber\\
&=pm_3^2+3p(p-1)m_{21}^2+p(p-1)(p-2)m_{111}^2\\
M_{2b}&=\sum_{ i \in \mathcal{I}}m^2[i,i,i]+\sum_{i, j \in \mathcal{I}, i \ne j}m^2[i,i,j]+\sum_{i, j \in \mathcal{I}, i \ne j}m[i,i,i]m[j,j,i]\nonumber\\
&\hspace{35mm}+\sum_{i, j \in \mathcal{I}, i \ne j}m[i,i,j]m[j,j,j]+\sum_{i, j, k \in \mathcal{I}, i \ne j, i\ne k, j\ne k}m[i,i,k] m[j,j,k]\nonumber\\
&=pm_3^2+p(p-1)m_{21}^2+2p(p-1)m_3m_{21}+p(p-1)(p-2)m_{21}^2\nonumber\\
&=pm_3^2+p(p-1)^2m_{21}^2+2p(p-1)m_3m_{21}\\
M_1&=\sum_{i\in \mathcal{I}} m[i, i, i, i]+\sum_{i, k \in \mathcal{I}, i \ne k} m[i, i, k, k]\nonumber\\
&=pm_4+p(p-1)m_{22}
\end{align}
Preliminarily we have the following results.
\begin{align}
&\sum_{i,j,k \in \mathcal{I}}\tilde{g}^{ii}\tilde{g}^{jj}\tilde{g}^{kk}m^2[i, j,k]\nonumber\\
&=\eta^{-3}[0,0,2,0]M_{2a}\nonumber\\
&\sum_{i,j \in \mathcal{I}}\tilde{g}^{ii}\tilde{g}^{jj}m^2[i, j]\nonumber\\
&=\eta^{-2}[0,0,2,0]\sum_{i,j \in \mathcal{I}}m^2[i, j]\nonumber\\
&=\eta^{-2}[0,0,2,0]\Bigl(pm_2^2+p(p-1)m_{11}^2\Bigr)\nonumber\\
&=\eta^{-2}[0,0,2,0] p, \\
&\sum_{i \in \mathcal{I}}\tilde{g}^{ii}m^2[i]=0.
\end{align}
\begin{align}
&\sum_{i,j,k \in \mathcal{I}}\tilde{g}^{ii}\tilde{g}^{jj}\tilde{g}^{kk}m[i, i,k]m[j, j,k]\nonumber\\
&=\eta^{-3}[0,0,2,0]M_{2b}\\
&\sum_{i,j \in \mathcal{I}}\tilde{g}^{ii}\tilde{g}^{jj}m[i, i]m[j,j]\nonumber\\
&=\eta^{-2}[0,0,2,0]\Bigl(\sum_{i \in \mathcal{I}}m^2[i,i]+\sum_{i,j \in \mathcal{I},i \ne j}m[i, i]m[j,j]\Bigr)\nonumber\\
&=\eta^{-2}[0,0,2,0]\Bigl(pm_2^2+p(p-1)m_{2}^2\Bigr)\nonumber\\
&=\eta^{-2}[0,0,2,0] p^2, \\
&\sum_{i,j \in \mathcal{I},i \ne j}\tilde{g}^{ii}\tilde{g}^{jj}m[i,i,j]m[j]=0\\
&\sum_{i\in \mathcal{I}}\tilde{g}^{ii}m^2[i]=0\\
&\sum_{i\in \mathcal{I}}\tilde{g}^{ii}m[i,i]=\eta^{-1}[0,0,2,0]p\\
&\sum_{i, j, k,l\in \mathcal{I}}\tilde{g}^{ij}\tilde{g}^{kl}m[i, j, k, l]\nonumber\\
&=\eta^{-2}[0,0,2,0]M_1\\
&\sum_{i, j \in \mathcal{I}} \tilde{g}^{ij} m[i, j]=\eta^{-1}[0,0,2,0] p
\end{align}

Let's consider $L21$ first.
\begin{align*}
L21&\triangleq g^{ij}g^{kl}g^{su}L_{(ik)s}L_{jlu}\\
&=\sum_{i,j,k,l,s,u \in \mathcal{I}}  \tilde{g}^{ij}\tilde{g}^{kl}\tilde{g}^{su} m[i,k,s]m[j,l,u]\eta_{(ik)s}\eta_{jlu}\\
&\quad+\sum_{i,j,k,l \in \mathcal{I}, s, u\in \mathcal{S}} \tilde{g}^{ij}\tilde{g}^{kl}\tilde{g}^{su} m[i,k,s]m[j,l,u]\eta_{(ik)s}\eta_{jlu}\\
&\quad+\sum_{i,j,s,u \in \mathcal{I}, k, l\in \mathcal{S}} \tilde{g}^{ij}\tilde{g}^{kl}\tilde{g}^{su} m[i,k,s]m[j,l,u]\eta_{(ik)s}\eta_{jlu}\\
&\quad+\sum_{k,l,s,u  \in \mathcal{I}, i,j\in \mathcal{S}} \tilde{g}^{ij}\tilde{g}^{kl}\tilde{g}^{su} m[i,k,s]m[j,l,u]\eta_{(ik)s}\eta_{jlu}\\
&\quad+\sum_{i,j \in \mathcal{I}, k, l, s, u\in \mathcal{S}} \tilde{g}^{ij}\tilde{g}^{kl}\tilde{g}^{su} m[i,k,s]m[j,l,u]\eta_{(ik)s}\eta_{jlu}\\
&\quad+\sum_{k,l  \in \mathcal{I}, i,j,s,u\in \mathcal{S}} \tilde{g}^{ij}\tilde{g}^{kl}\tilde{g}^{su} m[i,k,s]m[j,l,u]\eta_{(ik)s}\eta_{jlu}\\
&\quad+\sum_{s,u \in \mathcal{I}, i,j,k,l\in \mathcal{S}} \tilde{g}^{ij}\tilde{g}^{kl}\tilde{g}^{su} m[i,k,s]m[j,l,u]\eta_{(ik)s}\eta_{jlu}\\
&\quad+\sum_{i,j,k,l,s,u \in \mathcal{S}} \tilde{g}^{ij}\tilde{g}^{kl}\tilde{g}^{su} m[i,k,s]m[j,l,u]\eta_{(ik)s}\eta_{jlu}\\
&=\eta_{(00)0}\eta_{000}\sum_{i,j,k \in \mathcal{I}} \tilde{g}^{ii}\tilde{g}^{jj}\tilde{g}^{kk}m^2[i,j,k]
\\
&\quad+\sum_{i,j \in \mathcal{I}}\tilde{g}^{ii}\tilde{g}^{jj}m^2[i,j]\sum_{s, u\in \mathcal{S}}\tilde{g}^{su}\eta_{(00)s}\eta_{00u}\\
&\quad+\sum_{i,j\in \mathcal{I}} \tilde{g}^{ii}\tilde{g}^{jj}m^2[i,j]  \sum_{k, l\in \mathcal{S}} \tilde{g}^{kl}\eta_{(0k)0}\eta_{0l0}\\
&\quad+\sum_{k,l \in \mathcal{I}}\tilde{g}^{kk}\tilde{g}^{ll}m^2[k,l] \sum_{i, j\in \mathcal{S}}g^{ij}\eta_{(i0)0}\eta_{j00}\\
&\quad+\sum_{i\in \mathcal{I}}\tilde{g}^{ii}  m^2[i]\sum_{k, l, s, u\in \mathcal{S}}\tilde{g}^{kl}\tilde{g}^{su}\eta_{(0k)s}\eta_{0lu}\\
&\quad+\sum_{k\in \mathcal{I}}\tilde{g}^{kk}m^2[k]\sum_{i,j,s,u\in \mathcal{S}} \tilde{g}^{ij}\tilde{g}^{su} \eta_{(i0)s}\eta_{j0u}\\
&\quad+\sum_{s \in \mathcal{I}}\tilde{g}^{ss} m^2[s]\sum_{i,j,k,l\in \mathcal{S}} \tilde{g}^{ij}\tilde{g}^{kl} \eta_{(ik)0}\eta_{jl0}\\
&\quad+\sum_{i,j,k,l,s,u \in \mathcal{S}} \tilde{g}^{ij}\tilde{g}^{kl}\tilde{g}^{su} m[i,k,s]m[j,l,u]\eta_{(ik)s}\eta_{jlu}\\
&=\eta^{-3}[0,0,2,0]M_{2a} \eta_{(00)0}\eta_{000}\\
&\quad+\eta^{-2}[0,0,2,0]p\sum_{s, u\in \mathcal{S}}\tilde{g}^{su}\eta_{(00)s}\eta_{00u}\\
&\quad+\eta^{-2}[0,0,2,0]p \sum_{k, l\in \mathcal{S}} \tilde{g}^{kl}\eta_{(0k)0}\eta_{0l0}\\
&\quad+\eta^{-2}[0,0,2,0]p\sum_{i, j\in \mathcal{S}}\tilde{g}^{ij}\eta_{(i0)0}\eta_{j00}\\
&\quad+\sum_{i,j,k,l,s,u \in \mathcal{S}} \tilde{g}^{ij}\tilde{g}^{kl}\tilde{g}^{su}\eta_{(ik)s}\eta_{jlu}.
\end{align*}
Similarly we have
\begin{align*}
L22&=\eta^{-3}[0,0,2,0]M_{2b} \eta_{(00)0}\eta_{000}\\
&\quad+\eta^{-2}[0,0,2,0]p^2 \sum_{k, l\in \mathcal{S}} \tilde{g}^{kl}\eta_{(00)k}\eta_{l00}\\
&\quad+\eta^{-1}[0,0,2,0]p \sum_{k, l, s, u\in \mathcal{S}} \tilde{g}^{kl}\tilde{g}^{su}\eta_{(00)k}\eta_{lsu}\\
&\quad+\eta^{-1}[0,0,2,0]p \sum_{i, j, k, l\in \mathcal{S}} \tilde{g}^{ij}\tilde{g}^{kl}\eta_{(ij)k}\eta_{l00}\\
&\quad+\sum_{i,j,k,l,s,u \in \mathcal{S}} \tilde{g}^{ij}\tilde{g}^{kl}\tilde{g}^{su}\eta_{(ij)k}\eta_{lsu}.\\
L23&=\eta^{-3}[0,0,2,0]M_{2a} \eta_{000}\eta_{000}\\
&\quad+\eta^{-2}[0,0,2,0]p \sum_{s, u\in \mathcal{S}} \tilde{g}^{su}\eta_{00s}\eta_{00u}\\
&\quad+\eta^{-2}[0,0,2,0]p \sum_{k, l\in \mathcal{S}} \tilde{g}^{kl}\eta_{0k0}\eta_{0l0}\\
&\quad+\eta^{-2}[0,0,2,0]p \sum_{i, j\in \mathcal{S}} \tilde{g}^{ij}\eta_{i00}\eta_{j00}\\
&\quad+\sum_{i,j,k,l,s,u \in \mathcal{S}} \tilde{g}^{ij}\tilde{g}^{kl}\tilde{g}^{su}\eta_{iks}\eta_{jlu}.\\
L24&=\eta^{-3}[0,0,2,0]M_{2b} \eta_{000}\eta_{000}\\
&\quad+\eta^{-2}[0,0,2,0]p^2 \sum_{k, l\in \mathcal{S}} \tilde{g}^{kl}\eta_{00k}\eta_{l00}\\
&\quad+\eta^{-1}[0,0,2,0]p \sum_{k, l, s, u\in \mathcal{S}} \tilde{g}^{kl}\tilde{g}^{su}\eta_{00k}\eta_{lsu}\\
&\quad+\eta^{-1}[0,0,2,0]p \sum_{i, j, k, l\in \mathcal{S}} \tilde{g}^{ij}\tilde{g}^{kl}\eta_{ijk}\eta_{l00}\\
&\quad+\sum_{i,j,k,l,s,u \in \mathcal{S}} \tilde{g}^{ij}\tilde{g}^{kl}\tilde{g}^{su}\eta_{ijk}\eta_{lsu}.\\
L25&=\eta^{-3}[0,0,2,0]M_{2a} \eta_{(00)0}\eta_{(00)0}\\
&\quad+\eta^{-2}[0,0,2,0]p \sum_{s, u\in \mathcal{S}} \tilde{g}^{su}\eta_{(00)s}\eta_{(00)u}\\
&\quad+\eta^{-2}[0,0,2,0]p \sum_{k, l\in \mathcal{S}} \tilde{g}^{kl}\eta_{(0k)0}\eta_{(0l)0}\\
&\quad+\eta^{-2}[0,0,2,0]p \sum_{i, j\in \mathcal{S}} \tilde{g}^{ij}\eta_{(i0)0}\eta_{(j0)0}\\
&\quad+\sum_{i,j,k,l,s,u \in \mathcal{S}} \tilde{g}^{ij}\tilde{g}^{kl}\tilde{g}^{su}\eta_{(ik)s}\eta_{(jl)u}.\\
L26&=\eta^{-3}[0,0,2,0]M_{2b} \eta_{(00)0}\eta_{(00)0}\\
&\quad+\eta^{-2}[0,0,2,0]p^2 \sum_{k, l\in \mathcal{S}} \tilde{g}^{kl}\eta_{(00)k}\eta_{(00)l}\\
&\quad+\eta^{-1}[0,0,2,0]p \sum_{k, l, s, u\in \mathcal{S}} \tilde{g}^{kl}\tilde{g}^{su}\eta_{(00)k}\eta_{(su)l}\\
&\quad+\eta^{-1}[0,0,2,0]p \sum_{i, j, k, l\in \mathcal{S}} \tilde{g}^{ij}\tilde{g}^{kl}\eta_{(ij)k}\eta_{(00)l}\\
&\quad+\sum_{i,j,k,l,s,u \in \mathcal{S}} \tilde{g}^{ij}\tilde{g}^{kl}\tilde{g}^{su}\eta_{(ij)k}\eta_{(su)l}.\\
L11&=\sum_{i,j,k,l \in \mathcal{I}} \tilde{g}^{ij}\tilde{g}^{kl}m[i,j,k,l]\eta_{(il)jk}\\
&\quad+\sum_{i,j\in \mathcal{I}, k,l\in \mathcal{S}} \tilde{g}^{ij}\tilde{g}^{kl}m[i,j,k,l]\eta_{(il)jk}\\
&\quad+\sum_{i,j\in \mathcal{S}, k,l\in \mathcal{I}} \tilde{g}^{ij}\tilde{g}^{kl}m[i,j,k,l]\eta_{(il)jk}\\
&\quad+\sum_{i,j,k,l\in \mathcal{S}} \tilde{g}^{ij}\tilde{g}^{kl}m[i,j,k,l]\eta_{(il)jk}\\
&=\eta^{-2}[0,0,2,0]M_1\eta_{(00)00}\\
&\quad+\eta^{-1}[0,0,2,0]p\sum_{k, l \in \mathcal{S}}\tilde{g}^{kl}\eta_{(0l)0k}\\
&\quad+\eta^{-1}[0,0,2,0]p\sum_{i, j \in \mathcal{S}}\tilde{g}^{ij}\eta_{(i0)j0}\\
&\quad+\sum_{i,j,k,l\in \mathcal{S}} \tilde{g}^{ij}\tilde{g}^{kl}\eta_{(il)jk}.\\
L12&=\eta^{-2}[0,0,2,0]M_1\eta_{(00)00}\\
&\quad+\eta^{-1}[0,0,2,0]p\sum_{k, l \in \mathcal{S}}\tilde{g}^{kl}\eta_{(00)kl}\\
&\quad+\eta^{-1}[0,0,2,0]p\sum_{i, j \in \mathcal{S}}\tilde{g}^{ij}\eta_{(ij)00}\\
&\quad+\sum_{i,j,k,l\in \mathcal{S}} \tilde{g}^{ij}\tilde{g}^{kl}\eta_{(ij)kl}.\\
L13&=\eta^{-2}[0,0,2,0]M_1\eta_{0000}\\
&\quad+\eta^{-1}[0,0,2,0]p\sum_{k, l \in \mathcal{S}}\tilde{g}^{kl}\eta_{00kl}\\
&\quad+\eta^{-1}[0,0,2,0]p\sum_{i, j \in \mathcal{S}}\tilde{g}^{ij}\eta_{ij00}\\
&\quad+\sum_{i,j,k,l\in \mathcal{S}} \tilde{g}^{ij}\tilde{g}^{kl}\eta_{ijkl}.\\
L14&=\eta^{-2}[0,0,2,0]M_1\eta_{(00)(00)}\\
&\quad+\eta^{-1}[0,0,2,0]p\sum_{k, l \in \mathcal{S}}\tilde{g}^{kl}\eta_{(0k)(0l)}\\
&\quad+\eta^{-1}[0,0,2,0]p\sum_{i, j \in \mathcal{S}}\tilde{g}^{ij}\eta_{(i0)(j0)}\\
&\quad+\sum_{i,j,k,l\in \mathcal{S}} \tilde{g}^{ij}\tilde{g}^{kl}\eta_{(ik)(jl)}.\\
L15&=\eta^{-2}[0,0,2,0]M_1\eta_{(00)(00)}\\
&\quad+\eta^{-1}[0,0,2,0]p\sum_{k, l \in \mathcal{S}}\tilde{g}^{kl}\eta_{(00)(kl)}\\
&\quad+\eta^{-1}[0,0,2,0]p\sum_{i, j \in \mathcal{S}}\tilde{g}^{ij}\eta_{(ij)(00)}\\
&\quad+\sum_{i,j,k,l\in \mathcal{S}} \tilde{g}^{ij}\tilde{g}^{kl}\eta_{(ij)(kl)}.\\
\end{align*}
\section{Programming with Mathematica}
\label{Prog_Mathematica}

For the error term distributions, $N(0,1)$, $t(\nu)$, we theoretically derived the function $\eta[i, j, k, l]$ respectively as in \eqref{eta_normal}, \eqref{eta_t}. When the error term distribution is a skew-normal, we derived its values numerically by Monte Carlo simulation. We will call this process of the Mathematica programming "eta part". In the next step (called ``main part'' in the programming), we first calculated $\tilde{g}^{ij} (i, j \in \mathcal{I} \cup \mathcal{S})$ and \eqref{eta_(ij)k} to \eqref{eta_ssss}, and then \eqref{def_L} with another input $m[i, j, k], m[i, j, k, l]$ (instead $m_{4}, m_{22}, m_{3}, m_{21}, m_{111}$ when $x$ is homogeneous), The actual process for the calculation of \eqref{def_L} is stated in the last part of Appendix \ref{detailed_cal}.
Here we put the ``main part'' of the program of Mathematica. 
{\scriptsize
\begin{verbatim}
Main Part 
Inputs: 
1. eta[i, j, k, l] (which is calculated in Eta Part for each error distribution)
2-a. m4, m22, m3, m21, m111 ( m_{4}, m_{22}, m_{3}, m{21}, m{111} in the text) for homogeneous x
2-b. m[i, j, k], m[i, j, k, l] for nonhomogeneous x

Outputs:
eedn[a_,n_,p_] (\overset{\alpha}{E\!D} in the text). 

Note: In this program, "a" is used instead of  "alpha" in the text.

Delta
delta=eta[0,0,2,0](1+2eta[0,0,1,1]+eta[0,0,2,2])-(eta[0,1,0,1])^2

\tilde{g}^{ij} denoted by tgin[i,j]
tgin[0,0]=(1+2eta[0,0,1,1]+eta[0,0,2,2])/delta
tgin[0,sigma]=eta[0,1,0,1]/delta
tgin[sigma,0]=eta[0,1,0,1]/delta
tgin[sigma,sigma]=eta[0,0,2,0]/delta

eta_{(i,j),k} denoted by etau[{i_,j_},k_]
etau[{i_,j_}, k_]:=-eta[0,1,1,0]/;i==0&&j==0 &&k==0 
etau[{i_,j_}, k_]:=-(eta[0,1,1,1]+eta[0,0,2,0])/;i==0 && j==sigma && k==0 
etau[{i_,j_}, k_]:=-(eta[0,1,1,1]+eta[0,0,2,0])/;j==0 && i==sigma && k==0 
etau[{i_,j_}, k_]:=-(eta[0,1,0,0]+eta[0,1,1,1])/;i==0&&j==0 && k==sigma 
etau[{i_,j_}, k_]:=-(eta[0,1,0,1]+eta[0,1,1,2]+eta[0,0,2,1])/;i==0 && j==sigma && k==sigma 
etau[{i_,j_}, k_]:=-(eta[0,1,0,1]+eta[0,1,1,2]+eta[0,0,2,1])/;j==0 && i==sigma && k==sigma 
etau[{i_,j_}, k_]:=-(eta[0,1,1,2]+2*eta[0,0,2,1])/;i==sigma && j==sigma && k==0 
etau[{i_,j_}, k_]:=-(1+3*eta[0,0,1,1]+eta[0,1,0,2]+2*eta[0,0,2,2]+eta[0,1,1,3])
/;i==sigma && j==sigma && k==sigma

eta_{i,j,k} denoted by etau[i_,j_,k_]
etau[i_,j_,k_]:=-eta[0,0,3,0]/;i==0&&j==0 &&k==0 
etau[i_,j_,k_]:=-(eta[0,0,2,0]+eta[0,0,3,1])/;i==0&&j==0 &&k==sigma 
etau[i_,j_,k_]:=-(eta[0,0,2,0]+eta[0,0,3,1])/;i==0&&j==sigma &&k==0
etau[i_,j_,k_]:=-(eta[0,0,2,0]+eta[0,0,3,1])/;i==sigma&&j==0 &&k==0
etau[i_,j_,k_]:=-(eta[0,0,1,0]+2*eta[0,0,2,1]+eta[0,0,3,2])/;i==0&&j==sigma &&k==sigma 
etau[i_,j_,k_]:=-(eta[0,0,1,0]+2*eta[0,0,2,1]+eta[0,0,3,2])/;i==sigma&&j==0 &&k==sigma
etau[i_,j_,k_]:=-(eta[0,0,1,0]+2*eta[0,0,2,1]+eta[0,0,3,2])/;i==sigma&&j==sigma &&k==0
etau[i_,j_,k_]:=-(1+3*eta[0,0,1,1]+3*eta[0,0,2,2]+eta[0,0,3,3])
/;i==sigma&&j==sigma 
&&k==sigma 

eta_{(i,j)(k,l)} denoted by etau[{i_,j_},{k_,l_}]
etau[{i_,j_},{k_,l_}]:=eta[0,2,0,0]/;i==0&&j==0 &&k==0 &&l==0
etau[{i_,j_},{k_,l_}]:=eta[0,2,0,1]+eta[0,1,1,0]/;i==sigma&&j==0 &&k==0 &&l==0
etau[{i_,j_},{k_,l_}]:=eta[0,2,0,1]+eta[0,1,1,0]/;i==0&&j==sigma &&k==0 &&l==0
etau[{i_,j_},{k_,l_}]:=eta[0,2,0,1]+eta[0,1,1,0]/;i==0&&j==0 &&k==sigma &&l==0
etau[{i_,j_},{k_,l_}]:=eta[0,2,0,1]+eta[0,1,1,0]/;i==0&&j==0&&k==0 &&l==sigma
etau[{i_,j_},{k_,l_}]:=eta[0,2,0,2]+2*eta[0,1,1,1]+eta[0,0,2,0]
/;i==0&&j==sigma&&k==0 &&l==sigma
etau[{i_,j_},{k_,l_}]:=eta[0,2,0,2]+2*eta[0,1,1,1]+eta[0,0,2,0]
/;i==0&&j==sigma&&k==sigma &&l==0
etau[{i_,j_},{k_,l_}]:=eta[0,2,0,2]+2*eta[0,1,1,1]+eta[0,0,2,0]
/;i==sigma&&j==0&&k==0 &&l==sigma
etau[{i_,j_},{k_,l_}]:=eta[0,2,0,2]+2*eta[0,1,1,1]+eta[0,0,2,0]
/;i==sigma&&j==0&&k==sigma &&l==0
etau[{i_,j_},{k_,l_}]:=eta[0,1,0,0]+eta[0,2,0,2]+2*eta[0,1,1,1]
/;i==0&&j==0&&k==sigma &&l==sigma
etau[{i_,j_},{k_,l_}]:=eta[0,1,0,0]+eta[0,2,0,2]+2*eta[0,1,1,1]
/;i==sigma&&j==sigma&&k==0&&l==0
etau[{i_,j_},{k_,l_}]:=eta[0,1,0,1]+eta[0,2,0,3]+3*eta[0,1,1,2]+2*eta[0,0,2,1]
/;i==0&&j==sigma&&k==sigma &&l==sigma
etau[{i_,j_},{k_,l_}]:=eta[0,1,0,1]+eta[0,2,0,3]+3*eta[0,1,1,2]+2*eta[0,0,2,1]
/;i==sigma&&j==0&&k==sigma &&l==sigma
etau[{i_,j_},{k_,l_}]:=eta[0,1,0,1]+eta[0,2,0,3]+3*eta[0,1,1,2]+2*eta[0,0,2,1]
/;i==sigma&&j==sigma&&k==0 &&l==sigma
etau[{i_,j_},{k_,l_}]:=eta[0,1,0,1]+eta[0,2,0,3]+3*eta[0,1,1,2]+2*eta[0,0,2,1]
/;i==sigma&&j==sigma&&k==sigma &&l==0
etau[{i_,j_},{k_,l_}]:=1+eta[0,2,0,4]+4*eta[0,0,2,2]+2*eta[0,1,0,2]
+4*eta[0,0,1,1]+4*eta[0,1,1,3]
/;i==sigma&&j==sigma&&k==sigma &&l==sigma

eta_{(i,j,k),l} denoted by etau[{i_,j_,k_},l_]
etau[{i_,j_,k_},l_]:=eta[1,0,1,0]/;i==0&&j==0 &&k==0 &&l==0
etau[{i_,j_,k_},l_]:=eta[1,0,0,0]+eta[1,0,1,1]/;i==0&&j==0 &&k==0 &&l==sigma
etau[{i_,j_,k_},l_]:=2*eta[0,1,1,0]+eta[1,0,1,1]/;i==0&&j==0 &&k==sigma &&l==0
etau[{i_,j_,k_},l_]:=2*eta[0,1,1,0]+eta[1,0,1,1]/;i==0&&j==sigma &&k==0 &&l==0
etau[{i_,j_,k_},l_]:=2*eta[0,1,1,0]+eta[1,0,1,1]/;i==sigma&&j==0 &&k==0 &&l==0
etau[{i_,j_,k_},l_]:=4*eta[0,1,1,1]+2*eta[0,0,2,0]+eta[1,0,1,2]
/;i==0&&j==sigma &&k==sigma &&l==0
etau[{i_,j_,k_},l_]:=4*eta[0,1,1,1]+2*eta[0,0,2,0]+eta[1,0,1,2]
/;i==sigma&&j==0&&k==sigma &&l==0
etau[{i_,j_,k_},l_]:=4*eta[0,1,1,1]+2*eta[0,0,2,0]+eta[1,0,1,2]
/;i==sigma&&j==sigma &&k==0&&l==0
etau[{i_,j_,k_},l_]:=2*eta[0,1,0,0]+eta[1,0,0,1]+2*eta[0,1,1,1]+eta[1,0,1,2]
/;i==0&&j==0 &&k==sigma &&l==sigma
etau[{i_,j_,k_},l_]:=2*eta[0,1,0,0]+eta[1,0,0,1]+2*eta[0,1,1,1]+eta[1,0,1,2]
/;i==0&&j==sigma &&k==0 &&l==sigma
etau[{i_,j_,k_},l_]:=2*eta[0,1,0,0]+eta[1,0,0,1]+2*eta[0,1,1,1]+eta[1,0,1,2]
/;i==sigma&&j==0 &&k==0 &&l==sigma
etau[{i_,j_,k_},l_]:=4*eta[0,1,0,1]+eta[1,0,0,2]+4*eta[0,1,1,2]
+2*eta[0,0,2,1]+eta[1,0,1,3]/;i==0&&j==sigma &&k==sigma &&l==sigma
etau[{i_,j_,k_},l_]:=4*eta[0,1,0,1]+eta[1,0,0,2]+4*eta[0,1,1,2]
+2*eta[0,0,2,1]+eta[1,0,1,3]/;i==sigma&&j==0 &&k==sigma &&l==sigma
etau[{i_,j_,k_},l_]:=4*eta[0,1,0,1]+eta[1,0,0,2]+4*eta[0,1,1,2]
+2*eta[0,0,2,1]+eta[1,0,1,3]/;i==sigma&&j==sigma &&k==0 &&l==sigma
etau[{i_,j_,k_},l_]:=6*eta[0,1,1,2]+6*eta[0,0,2,1]+eta[1,0,1,3]
/;i==sigma&&j==sigma &&k==sigma &&l==0
etau[{i_,j_,k_},l_]:=2+6*eta[0,1,0,2]+6*eta[0,0,1,1]+eta[1,0,0,3]
+2*eta[0,0,1,1]+6*eta[0,1,1,3]+6*eta[0,0,2,2]+eta[1,0,1,4]
/;i==sigma&&j==sigma &&k==sigma &&l==sigma

eta_{(i,j),k,l} denoted by etau[{i_,j_},k_,l_]
etau[{i_,j_},k_,l_]:=eta[0,1,2,0]/;i==0&&j==0 &&k==0 &&l==0
etau[{i_,j_},k_,l_]:=eta[0,1,1,0]+eta[0,1,2,1]/;i==0&&j==0 &&k==0 &&l==sigma
etau[{i_,j_},k_,l_]:=eta[0,1,1,0]+eta[0,1,2,1]/;i==0&&j==0 &&k==sigma &&l==0
etau[{i_,j_},k_,l_]:=eta[0,1,2,1]+eta[0,0,3,0]/;i==0&&j==sigma &&k==0 &&l==0
etau[{i_,j_},k_,l_]:=eta[0,1,2,1]+eta[0,0,3,0]/;i==sigma&&j==0 &&k==0 &&l==0
etau[{i_,j_},k_,l_]:=eta[0,1,0,0]+2*eta[0,1,1,1]+eta[0,1,2,2]
/;i==0&&j==0 &&k==sigma &&l==sigma
etau[{i_,j_},k_,l_]:=eta[0,1,1,1]+eta[0,0,2,0]+eta[0,1,2,2]+eta[0,0,3,1]
/;i==0&&j==sigma &&k==0 &&l==sigma
etau[{i_,j_},k_,l_]:=eta[0,1,1,1]+eta[0,0,2,0]+eta[0,1,2,2]+eta[0,0,3,1]
/;i==0&&j==sigma &&k==sigma &&l==0
etau[{i_,j_},k_,l_]:=eta[0,1,1,1]+eta[0,0,2,0]+eta[0,1,2,2]+eta[0,0,3,1]
/;i==sigma&&j==0 &&k==0 &&l==sigma
etau[{i_,j_},k_,l_]:=eta[0,1,1,1]+eta[0,0,2,0]+eta[0,1,2,2]+eta[0,0,3,1]
/;i==sigma&&j==0 &&k==sigma &&l==0
etau[{i_,j_},k_,l_]:=eta[0,0,2,0]+2*eta[0,0,3,1]+eta[0,1,2,2]
/;i==sigma&&j==sigma &&k==0 &&l==0
etau[{i_,j_},k_,l_]:=eta[0,1,0,1]+2*eta[0,1,1,2]+2*eta[0,0,2,1]
+eta[0,1,2,3]+eta[0,0,3,2]/;i==0&&j==sigma &&k==sigma &&l==sigma
etau[{i_,j_},k_,l_]:=eta[0,1,0,1]+2*eta[0,1,1,2]+2*eta[0,0,2,1]
+eta[0,1,2,3]+eta[0,0,3,2]/;i==sigma&&j==0 &&k==sigma &&l==sigma
etau[{i_,j_},k_,l_]:=2*eta[0,0,2,1]+eta[0,1,1,2]+eta[0,0,2,1]
+2*eta[0,0,3,2]+eta[0,1,2,3]/;i==sigma&&j==sigma &&k==0 &&l==sigma
etau[{i_,j_},k_,l_]:=2*eta[0,0,2,1
]+eta[0,1,1,2]+eta[0,0,2,1]+2*eta[0,0,3,2]+eta[0,1,2,3]
/;i==sigma&&j==sigma &&k==sigma &&l==0
etau[{i_,j_},k_,l_]:=1+4*eta[0,0,1,1]+eta[0,1,0,2]+5*eta[0,0,2,2]
+2*eta[0,1,1,3]+2*eta[0,0,3,3]+eta[0,1,2,4]
/;i==sigma&&j==sigma &&k==sigma &&l==sigma

eta_{i,j,k,l} denoted by etau[i_,j_,k_,l_]
etau[i_,j_,k_,l_]:=eta[0,0,4,0]/;i==0&&j==0 &&k==0 &&l==0
etau[i_,j_,k_,l_]:=eta[0,0,3,0]+eta[0,0,4,1]/;i==0&&j==0 &&k==0 &&l==sigma
etau[i_,j_,k_,l_]:=eta[0,0,3,0]+eta[0,0,4,1]/;i==0&&j==0 &&k==sigma &&l==0
etau[i_,j_,k_,l_]:=eta[0,0,3,0]+eta[0,0,4,1]/;i==0&&j==sigma &&k==0 &&l==0
etau[i_,j_,k_,l_]:=eta[0,0,3,0]+eta[0,0,4,1]/;i==sigma&&j==0 &&k==0 &&l==0
etau[i_,j_,k_,l_]:=eta[0,0,2,0]+2*eta[0,0,3,1]+eta[0,0,4,2]
/;i==sigma&&j==sigma &&k==0 &&l==0
etau[i_,j_,k_,l_]:=eta[0,0,2,0]+2*eta[0,0,3,1]+eta[0,0,4,2]
/;i==sigma&&j==0 &&k==sigma &&l==0
etau[i_,j_,k_,l_]:=eta[0,0,2,0]+2*eta[0,0,3,1]+eta[0,0,4,2]
/;i==sigma&&j==0 &&k==0 &&l==sigma
etau[i_,j_,k_,l_]:=eta[0,0,2,0]+2*eta[0,0,3,1]+eta[0,0,4,2]
/;i==0&&j==sigma &&k==sigma &&l==0
etau[i_,j_,k_,l_]:=eta[0,0,2,0]+2*eta[0,0,3,1]+eta[0,0,4,2]
/;i==0&&j==sigma &&k==0&&l==sigma
etau[i_,j_,k_,l_]:=eta[0,0,2,0]+2*eta[0,0,3,1]+eta[0,0,4,2]
/;i==0&&j==0 &&k==sigma &&l==sigma
etau[i_,j_,k_,l_]:=3*eta[0,0,2,1]+3*eta[0,0,3,2]+eta[0,0,4,3]
/;i==0&&j==sigma &&k==sigma &&l==sigma
etau[i_,j_,k_,l_]:=3*eta[0,0,2,1]+3*eta[0,0,3,2]+eta[0,0,4,3]
/;i==sigma&&j==0 &&k==sigma &&l==sigma
etau[i_,j_,k_,l_]:=3*eta[0,0,2,1]+3*eta[0,0,3,2]+eta[0,0,4,3]
/;i==sigma&&j==sigma &&k==0 &&l==sigma
etau[i_,j_,k_,l_]:=3*eta[0,0,2,1]+3*eta[0,0,3,2]+eta[0,0,4,3]
/;i==sigma&&j==sigma &&k==sigma &&l==0
etau[i_,j_,k_,l_]:=1+4*eta[0,0,1,1]+6*eta[0,0,2,2]+4*eta[0,0,3,3]+eta[0,0,4,4]
/;i==sigma&&j==sigma &&k==sigma &&l==sigma

M_{2a}(p), M_{2b}(p), M_1(p) denoted respectively by mmu2a[p_], mmu2b[p_],mmu1[p_]
(*mmu2a[p_]:=Sum[m[i, j, k]^2, {i, 1,p},{j,1,p},{k,1,p}]*)
(*mmu2b[p_]:=Sum[m[i, i, k]*m[j,j,k], {i, 1,p},{j,1,p},{k,1,p}]*)
(*mmu1[p_]:=Sum[m[i,i,j,j],{i, 1,p},{j,1,p}]*)
M_{2a}(p), M_{2b}(p), M_1(p) for the special case when $x$ is homogeneous
mmu2a[p_]:=p*m3^2+3*p*(p-1)*m21^2+p*(p-1)*(p-2)*m111^2
mmu2b[p_]:=p*m3^2+p*(p-1)^2*m21^2+2*p*(p-1)*m3*m21
mmu1[p_]:=p*m4+p*(p-1)*m22

L21,..., L26, L11,...,L15 denoted by ll21,...,ll26,ll11,...,ll15
ll21[p_]:=eta[0,0,2,0]^{-3}mmu2a[p]etau[{0,0},0]etau[0,0,0]
+eta[0,0,2,0]^{-2}p*Sum[tgin[s,u]etau[{0,0},s]etau[0,0,u],{s,{0,sigma}},{u,{0,sigma}}]
+eta[0,0,2,0]^{-2}p*Sum[tgin[k,l]etau[{0,k},0]etau[0,l,0],{k,{0,sigma}},{l,{0,sigma}}]+
eta[0,0,2,0]^{-2}p*Sum[tgin[i,j]etau[{i,0},0]etau[j,0,0],{i,{0,sigma}},{j,{0,sigma}}]+
Sum[tgin[i,j]tgin[k,l]tgin[s,u]etau[{i,k},s]etau[j,l,u],{i,{0,sigma}},{j,{0,sigma}},
{k,{0,sigma}},{l,{0,sigma}},{s,{0,sigma}},{u,{0,sigma}}]

ll22[p_]:=eta[0,0,2,0]^{-3}mmu2b[p]etau[{0,0},0]etau[0,0,0]
+eta[0,0,2,0]^{-2}p^2*Sum[tgin[k,l]etau[{0,0},k]etau[l,0,0],{k,{0,sigma}},{l,{0,sigma}}]
+eta[0,0,2,0]^{-1}p*Sum[tgin[k,l]tgin[s,u]etau[{0,0},k]etau[l,s,u],{k,{0,sigma}},
{l,{0,sigma}},{s,{0,sigma}},{u,{0,sigma}}]+eta[0,0,2,0]^{-1}p*Sum[tgin[i,j]tgin[k,l]etau[{i,j},k]
etau[l,0,0],{i,{0,sigma}},{j,{0,sigma}},{k,{0,sigma}},{l,{0,sigma}}]+Sum[tgin[i,j]tgin[k,l]tgin[s,u]
etau[{i,j},k]etau[l,s,u],{i,{0,sigma}},{j,{0,sigma}},{k,{0,sigma}},{l,{0,sigma}},{s,{0,sigma}},{u,{0,sigma}}]

ll23[p_]:=eta[0,0,2,0]^{-3}mmu2a[p]etau[0,0,0]etau[0,0,0]
+eta[0,0,2,0]^{-2}p*Sum[tgin[s,u]etau[0,0,s]etau[0,0,u],{s,{0,sigma}},{u,{0,sigma}}]
+eta[0,0,2,0]^{-2}p*Sum[tgin[k,l]etau[0,k,0]etau[0,l,0],{k,{0,sigma}},{l,{0,sigma}}]+
eta[0,0,2,0]^{-2}p*Sum[tgin[i,j]etau[i,0,0]etau[j,0,0],{i,{0,sigma}},{j,{0,sigma}}]+
Sum[tgin[i,j]tgin[k,l]tgin[s,u]etau[i,k,s]etau[j,l,u],{i,{0,sigma}},{j,{0,sigma}},{k,{0,sigma}},
{l,{0,sigma}},{s,{0,sigma}},{u,{0,sigma}}]

ll24[p_]:=eta[0,0,2,0]^{-3}mmu2b[p]etau[0,0,0]etau[0,0,0]
+eta[0,0,2,0]^{-2}p^2*Sum[tgin[k,l]etau[0,0,k]etau[l,0,0],{k,{0,sigma}},{l,{0,sigma}}]
+eta[0,0,2,0]^{-1}p*Sum[tgin[k,l]tgin[s,u]etau[0,0,k]etau[l,s,u],{k,{0,sigma}},{l,{0,sigma}},
{s,{0,sigma}},{u,{0,sigma}}]+eta[0,0,2,0]^{-1}p*Sum[tgin[i,j]tgin[k,l]etau[i,j,k]etau[l,0,0],
{i,{0,sigma}},{j,{0,sigma}},{k,{0,sigma}},{l,{0,sigma}}]
+Sum[tgin[i,j]tgin[k,l]tgin[s,u]etau[i,j,k]etau[l,s,u],{i,{0,sigma}},{j,{0,sigma}},
{k,{0,sigma}},{l,{0,sigma}},{s,{0,sigma}},{u,{0,sigma}}]

ll25[p_]:=eta[0,0,2,0]^{-3}mmu2a[p]etau[{0,0},0]etau[{0,0},0]
+eta[0,0,2,0]^{-2}p*Sum[tgin[s,u]etau[{0,0},s]etau[{0,0},u],{s,{0,sigma}},{u,{0,sigma}}]
+eta[0,0,2,0]^{-2}p*Sum[tgin[k,l]etau[{0,k},0]etau[{0,l},0],{k,{0,sigma}},{l,{0,sigma}}]
+eta[0,0,2,0]^{-2}p*Sum[tgin[i,j]etau[{i,0},0]etau[{j,0},0],{i,{0,sigma}},{j,{0,sigma}}]
+Sum[tgin[i,j]tgin[k,l]tgin[s,u]etau[{i,k},s]etau[{j,l},u],{i,{0,sigma}},{j,{0,sigma}},{k,{0,sigma}},
{l,{0,sigma}},{s,{0,sigma}},{u,{0,sigma}}]

ll26[p_]:=eta[0,0,2,0]^{-3}mmu2b[p]etau[{0,0},0]etau[{0,0},0]
+eta[0,0,2,0]^{-2}p^2*Sum[tgin[k,l]etau[{0,0},k]etau[{0,0},l],{k,{0,sigma}},{l,{0,sigma}}]
+eta[0,0,2,0]^{-1}p*Sum[tgin[k,l]tgin[s,u]etau[{0,0},k]etau[{s,u},l],{k,{0,sigma}},
{l,{0,sigma}},{s,{0,sigma}},{u,{0,sigma}}]+eta[0,0,2,0]^{-1}p*Sum[tgin[i,j]tgin[k,l]
etau[{i,j},k]etau[{0,0},l],{i,{0,sigma}},{j,{0,sigma}},{k,{0,sigma}},{l,{0,sigma}}]
+Sum[tgin[i,j]tgin[k,l]tgin[s,u]etau[{i,j},k]etau[{s,u},l],{i,{0,sigma}},{j,{0,sigma}},
{k,{0,sigma}},{l,{0,sigma}},{s,{0,sigma}},{u,{0,sigma}}]

ll11[p_]:=eta[0,0,2,0]^{-2}mmu1[p]etau[{0,0},0,0]
+eta[0,0,2,0]^{-1}p*Sum[tgin[k,l]etau[{0,l},0,k],{k,{0,sigma}},{l,{0,sigma}}]
+eta[0,0,2,0]^{-1}p*Sum[tgin[i,j]etau[{i,0},j,0],{i,{0,sigma}},{j,{0,sigma}}]+Sum[tgin[i,j]tgin[k,l]etau[{i,l},j,k],
{i,{0,sigma}},{j,{0,sigma}},{k,{0,sigma}},{l,{0,sigma}}]

ll12[p_]:=eta[0,0,2,0]^{-2}mmu1[p]etau[{0,0},{0,0}]
+eta[0,0,2,0]^{-1}p*Sum[tgin[k,l]etau[{0,0},k,l],{k,{0,sigma}},{l,{0,sigma}}]
+eta[0,0,2,0]^{-1}p*Sum[tgin[i,j]etau[{i,j},0,0],{i,{0,sigma}},{j,{0,sigma}}]
+Sum[tgin[i,j]tgin[k,l]etau[{i,j},k,l],{i,{0,sigma}},{j,{0,sigma}},{k,{0,sigma}},{l,{0,sigma}}]

ll13[p_]:=eta[0,0,2,0]^{-2}mmu1[p]etau[0,0,0,0]
+eta[0,0,2,0]^{-1}p*Sum[tgin[k,l]etau[0,0,k,l],{k,{0,sigma}},{l,{0,sigma}}]
+eta[0,0,2,0]^{-1}p*Sum[tgin[i,j]etau[i,j,0,0],{i,{0,sigma}},{j,{0,sigma}}]
+Sum[tgin[i,j]tgin[k,l]etau[i,j,k,l],{i,{0,sigma}},{j,{0,sigma}},{k,{0,sigma}},{l,{0,sigma}}]

ll14[p_]:=eta[0,0,2,0]^{-2}mmu1[p]etau[{0,0},{0,0}]
+eta[0,0,2,0]^{-1}p*Sum[tgin[k,l]etau[{0,l},{0,k}],{k,{0,sigma}},{l,{0,sigma}}]
+eta[0,0,2,0]^{-1}p*Sum[tgin[i,j]etau[{i,0},{j,0}],{i,{0,sigma}},{j,{0,sigma}}]
+Sum[tgin[i,j]tgin[k,l]etau[{i,k},{j,l}],{i,{0,sigma}},{j,{0,sigma}},{k,{0,sigma}},{l,{0,sigma}}]

ll15[p_]:=eta[0,0,2,0]^{-2}mmu1[p]etau[{0,0},{0,0}]
+eta[0,0,2,0]^{-1}p*Sum[tgin[k,l]etau[{0,0},{k,l}],{k,{0,sigma}},{l,{0,sigma}}]
+eta[0,0,2,0]^{-1}p*Sum[tgin[i,j]etau[{i,j},{0,0}],{i,{0,sigma}},{j,{0,sigma}}]
+Sum[tgin[i,j]tgin[k,l]etau[{i,j},{k,l}],{i,{0,sigma}},{j,{0,sigma}},{k,{0,sigma}},{l,{0,sigma}}]

F^e denoted by  ffe
ffe[p_]:=2ll11[p]+ll12[p]+ll13[p]-2ll21[p]-ll23[p]-ll22[p]

T_{ijk}T^{ijk} denoted by tt1
tt1[p_]:=ll23[p]

T_{is}^iT_j^{js} denoted by tt2
tt2[p_]:=ll24[p]

R^e_{ij}^{ij} denoted by rre
rre[p_]:=ll14[p]-ll15[p]+ll11[p]-ll12[p]-ll25[p]+ll26[p]+ll22[p]-ll21[p]

<A^e_i^j, A^e_j^i> denoted by aaee1 
aaee1[p_]:=ll14[p]-ll25[p]-p

<A^e_i^i, A^e_j^j> denoted by aaee2 
aaee2[p_]:=ll15[p]-ll26[p]-p^2

<A^e_i^j, A^m_j^i> denoted by aaem1
aaem1[p_]:=ll11[p]+ll14[p]-ll25[p]-ll21[p]

<A^e_i^i, A^m_j^j> denoted by aaem2
aaem2[p_]:=ll12[p]+ll15[p]-ll26[p]-ll22[p]

\overset{\alpha}{E\!D} denoted by eedn[a_,n_,p_]
eedn[a_,n_,p_]:=(p+2)/(2n)+1/(24n^2)*(((1-a)/2)^2
*(3ffe[p]+3tt1[p]-6aaem1[p]+6aaee1[p]-3aaem2[p]+3aaee2[p]+3p^2+6p)
+((1-a)/2)*(3ffe[p]-5tt1[p]-6tt2[p]+6aaem1[p]-6aaee1[p]+3aaem2[p]-3aaee2[p]-3p^2-6p)
+12aaee1[p]-2aaem1[p]-aaem2[p]+tt1[p]+9tt2[p]+8rre[p]-9ffe[p])
Simplify[eedn[-1,n,p]]
Simplify[eedn[1,n,p]]
Simplify[eedn[0,n,p]]
Simplify[eedn[2,n,p]]
\end{verbatim}
}

\section*{Acknowledgment}
This research was supported by JSPS KAKENHI Grant Number 25380265.


\begin{thebibliography}{99}
\bibitem{Amari3}
S. Amari.
Alpha divergence is unique, belonging to both classes of f-divergence and Bregman divergence
\textit{IEEE Trans. Information Theory},  55:4925-4931, 2009.

\bibitem{Amari&Cichocki}
S. Amari and A. Chichocki.
Information geometry of divergence function.
\textit{Bulletin of the Polish Academy of  Sciences : Technical Sciences}, 58:183-195, 2010.


\bibitem{Amari&Nagaoka}
S. Amari and H. Nagaoka.
\textit{Methods of Information Geometry}. 
Translations of Mathematical Monographs 191. American Mathematical Society, 2000.

\bibitem{Cortez_et_al}
P. Cortez, A. Cerdeira, F. Almeida, T. Matos and J. Reis. 
Modeling wine preferences by data mining from physicochemical properties. 
\textit{Decision Support Systems}, 47(4): 547-553, 2009.

\bibitem{Eguchi1}
S. Eguchi.
A differential geometric approach to statistical inference on the basis of contrast functionals.
\textit{Hiroshima Mathematical Journal}, 15: 341-391, 1985.

\bibitem{Mathematica}
Wolfram Research, Inc.
\textit{Mathematica 10.4}, 2016.

\bibitem{Sheena}
Y. Sheena.
\textit{Asymptotic expansion of the risk of maximum likelihood estimator with respect to $\alpha$-divergence as a measure of the difficulty of specifying a parametric model}, arXiv:1510.08226, 2016.


\bibitem{Vajda}
I. Vajda.
\textit{Theory of Statistical Inference and Information}, Kluwer Academic Publishers, 1989.



\end{thebibliography}
\end{document}